\documentclass[review,3p,times,12pt]{elsarticle}
\usepackage{lineno,hyperref,amsmath,amsthm}
\usepackage[noend]{algpseudocode}
\usepackage{algorithmicx,algorithm}
\usepackage{mathrsfs}
\usepackage{booktabs}
\usepackage{color}
\usepackage{graphics}
\usepackage{subfig}
\usepackage{amsfonts,amssymb}
\usepackage{bm}
\allowdisplaybreaks[4]
\usepackage{threeparttable}

\renewcommand{\eqref}[1]{\textup{{\normalfont Eq.~(\ref{#1}}\normalfont)}}
\newcommand\figref{Fig.~\ref}

\modulolinenumbers[1]

\journal{arXiv}

\begin{document}
	
	\begin{frontmatter}
		
		\title{Stochastic virtual element methods for uncertainty propagation of stochastic linear elasticity}

        \author{Zhibao Zheng}
		
		\author{Udo Nackenhorst}
		
		\address{Leibniz Universit{\"a}t Hannover, Institute of Mechanics and Computational Mechanics \& International Research Training Group 2657, Appelstra{\ss}e 9a, 30167, Hannover, Germany}
	
		\begin{abstract}
            This paper presents stochastic virtual element methods for propagating uncertainty in linear elastic stochastic problems.
            We first derive stochastic virtual element equations for 2D and 3D linear elastic problems that may involve uncertainties in material properties, external forces, boundary conditions, etc.
            A stochastic virtual element space that couples the deterministic virtual element space and the stochastic space is constructed for this purpose and used to approximate the unknown stochastic solution.
            Two numerical frameworks are then developed to solve the derived stochastic virtual element equations, including a Polynomial Chaos approximation based approach and a weakly intrusive approximation based approach.
            In the PC based framework, the stochastic solution is approximated using the Polynomial Chaos basis and solved via an augmented deterministic virtual element equation that is generated by applying the stochastic Galerkin procedure to the original stochastic virtual element equation.
            In the weakly intrusive approximation based framework, the stochastic solution is approximated by a summation of a set of products of random variables and deterministic vectors, where the deterministic vectors are solved via converting the original stochastic problem to deterministic virtual element equations by the stochastic Galerkin approach, and the random variables are solved via converting the original stochastic problem to one-dimensional stochastic algebraic equations by the classical Galerkin procedure.
            This method avoids the curse of dimensionality of high-dimensional stochastic problems successfully since all random inputs are embedded into one-dimensional stochastic algebraic equations whose computational effort weakly depends on the stochastic dimension.
    		Numerical results on 2D and 3D problems with low- and high-dimensional random inputs demonstrate the good performance of the proposed methods.
		\end{abstract}
		
		
	\end{frontmatter}
	
	
\section{Introduction} \label{sec:Intro}

    Numerical techniques for solving complex partial differential equations are continuously developing at an incredible rate, e.g. finite difference methods, finite volume methods, finite element methods and spectral element methods etc \cite{quarteroni2008numerical}.
    As a generalization of the finite element method, the virtual element method (VEM) has been proposed and received a lot of attention in the last decade \cite{beirao2013basic, beirao2014hitchhiker, mengolini2019engineering, antonietti2022virtual}.
    Compared to the classical finite element method, VEM can discretize 2D or 3D geometric domains utilizing arbitrary polygons or polyhedrons and is not limited to the regular elements used in the finite element method, which is thus highly flexible and mesh insensitive. 
    In VEM, the shape functions can be non-polynomial. 
    It does not require constructing explicit shape functions on elements since all numerical integration is transferred to edges rather than performed on the elements.
    Several limitations of finite element methods are also avoided, such as convex elements and element degradation caused by small edges and interior angles.
    Applications of VEM to various problems have been extensively studied, such as linear elastic problems \cite{da2013virtual, gain2014virtual, artioli2017arbitrary}, large deformation problems \cite{chi2017some, wriggers2021taylor}, contact problems \cite{wriggers2016virtual, wriggers2019virtual, aldakheel2020curvilinear}, fracture and crack propagation problems \cite{benedetto2014virtual, aldakheel2018phase, hussein2020combined}, topology optimization \cite{antonietti2017virtual, chi2020virtual} etc.
    Although extensive studies have been conducted, there is still a gap in using VEM to deal with problems with uncertainties.
    In many practical engineering problems, the inherent or epistemic uncertainty of systems are unavoidable.
    Predicting uncertainty propagation on the physical models has become an important part of the analysis of systems \cite{smith2013uncertainty}, which leads to the development of dedicated numerical methods for uncertainty analysis. 
  
    In this paper, we focus on extending the deterministic VEM to stochastic VEMs (SVEMs) for the uncertainty analysis of 2D and 3D linear elastic stochastic problems that may involve random material properties and stochastic external forces, etc., which is currently still lacking.
    Our main contributions in this paper consist of two parts: the first contribution is to extend the deterministic virtual element discretization to stochastic cases and derive corresponding stochastic virtual element equations (SVEEs), and the second contribution is to present two numerical methods to solve the derived SVEEs efficiently and accurately.
    For the first contribution, we extend the deterministic virtual element space to a stochastic virtual element space that couples the classical virtual element space and the stochastic space, which can provide (stochastic) virtual element approximations for both deterministic and stochastic functions in the space.
    The constructed stochastic virtual element space can be considered as a deterministic virtual element space parameterized by random inputs.
    For each sample realization of the random input, it degenerates into a deterministic virtual element space and inherits all properties of classical virtual element spaces.
    Thus, we can simply approximate the stochastic solution using a linear combination of deterministic virtual basis functions with random coefficients (i.e. the unknown stochastic solution).
    Further, numerical techniques for calculating the gradients of virtual basis functions and the stabilization term in the deterministic VEM \cite{beirao2014hitchhiker, gain2014virtual, artioli2017arbitrary, chi2020virtual} can still be used to the stochastic discretization with slight modifications.
    In this way we can obtain SVEEs by assembling the stochastic stiffness matrix and the stochastic force vector with a complexity similar to the deterministic VEM.

    Similar to VEM being a generalization of the finite element method, SVEM is also considered to be a generalization of the stochastic finite element method \cite{stefanou2009stochastic}.
    Thus, numerical algorithms for solving the derived SVEEs can benefit from stochastic finite element solution algorithms, such as Monte Carlo simulation (MCS) and its improvements \cite{papadrakakis1996robust, graham2011quasi}, spectral stochastic finite element methods \cite{ghanem2003stochastic, xiu2002wiener}, stochastic collocation methods \cite{babuvska2007stochastic, xiu2010numerical}, response surface and kriging methods \cite{khuri2010response, fuhg2021state}, etc.
    For the second contribution, we develop two numerical methods to solve the derived SVEE, including a Polynomial Chaos expansion based SVEM (PC-SVEM) and a Weakly INtrusive approximation based SVEM (WIN-SVEM).
    The other methods mentioned above can also be extended to SVEMs in a similar way as in this paper.
    The PC-SVEM is a natural extension of the spectral stochastic finite element method \cite{ghanem2003stochastic}.
    In this method, the stochastic solution is decomposed into a summation of a set of products of PC basis and deterministic vectors.
    By the use of stochastic Galerkin procedure, the original SVEE is transformed into an augmented deterministic equation whose size is much larger than the original SVEE.
    Also, the size increases dramatically as the degree of freedom of physical models, the stochastic dimension and the expansion order of PC basis increase, which leads to the curse of dimensionality when dealing with large-scale and/or high-dimensional stochastic problems.
    To address this issue, we further present a WIN-SVEM, which is an extension of our previous work for solving stochastic finite element equations \cite{zheng2022weak, zheng2023stochastic}.
    In this method, the stochastic solution is approximated by a summation of a set of products of random variables and deterministic vectors.
    Different from the PC-SVEM, both random variables and deterministic vectors are unknown a priori.
    To this end, we solve them using a dedicated iteration.
    The deterministic vectors are solved via a few number of deterministic equations that are obtained by a similar stochastic Galerkin process used to PC-SVEM.
    The random variables are solved via one-dimensional stochastic algebraic equations that are obtained by applying the classical Galerkin procedure to the original SVEE.
    In this way, all random inputs are embedded into these one-dimensional stochastic algebraic equations, and their efficient solutions are achieved using a non-intrusive sampling method with weak dimension dependence.
    The proposed WIN-SVEM thus avoids the curse of dimensionality of high-dimensional stochastic problems successfully.

	The paper is organized as follows: 
    Section \ref{sec:SVEE} presents the stochastic virtual element discretization for linear elastic stochastic problems and stochastic virtual element equations are then derived.
    In Section \ref{sec:PC_SVEM}, the PC-SVEM is developed to solve the derived stochastic systems.
    Following that, the WIN-SVEM is proposed in Section \ref{sec:WIN_SVEM} to solve the derived stochastic systems efficiently, with special emphasis on high-dimensional stochastic problems.
    2D and 3D numerical examples involving low- and high-dimensional random inputs are given in Section \ref{sec:Ex} to demonstrate the performance of the proposed methods.
    Conclusions and outlook follow in Section \ref{sec:Con}.



\section{Stochastic virtual element equations} \label{sec:SVEE}

\subsection{Stochastic elastic equations} \label{subsection:EP}
	Let $( {\Theta, \Xi, \cal{P}} )$ be a suitable probability space , where $\Theta$ denotes the space of elementary events, $\Xi$ is a $\sigma$-algebra defined on $\Theta$ and $\cal{P}$ is a probability measure. 
	In this paper, we consider the following elastic stochastic equation

	\vspace{-0.5cm}
	\begin{equation} \label{eq:EP_probelm}
    	\left \{  
    	\begin{aligned}
        	- \nabla \cdot {\bm \sigma} \left( {\bf x},\theta \right) &= {\bm f}\left( {\bf x},\theta \right) ~~&{\rm{in}}~~&\Omega \\
        	{\bm \sigma} \left( {\bf x},\theta \right) \cdot {\bf n} &= {\bm g}\left( {\bf x},\theta \right)~~&{\rm on}~~&{\Gamma _N}\\
        	{\bm u}\left( {\bf x},\theta \right) &= {\bm u}_D\left( {\bf x},\theta \right)~~&{\rm on}~~&{\Gamma _D}
    	\end{aligned}  
    	\right. ,
	\end{equation}
	where the deterministic domain $\Omega \subset \mathbb{R}^{d}$ with the boundary $\partial \Omega$, the spatial dimension may be $d = 2,3$ and the $d$-dimensional spatial coordinate is given by ${\bf x} = \left( x_1, \cdots, x_d \right) \in \Omega$, $\nabla \cdot \left( \cdot \right)$ denotes the divergence operator, ${\bm \sigma} \left( {{\bf{x}},\theta } \right)$ is the stochastic stress tensor, 
    the vector-valued displacement field ${\bm u}\left( {\bf x}, \theta \right) = \left[ u_1\left( {\bf x}, \theta \right), \cdots, u_d\left( {\bf x}, \theta \right) \right]^{\rm T} \in \mathbb{R}^d$ is the unknown stochastic solution to be solved,
    the vector-valued field ${\bm f}\left( {{\bf{x}},\theta } \right) = \left[ f_1\left( {\bf x}, \theta \right), \cdots, f_d\left( {\bf x}, \theta \right) \right]^{\rm T} \in \mathbb{R}^d$ is associated with stochastic external forces, and ${\Gamma _N}$ and ${\Gamma _D}$ are boundary segments associated with the Neumann boundary condition ${\bm g}\left( {{\bf{x}},\theta } \right) = \left[ g_1\left( {\bf x}, \theta \right), \cdots, g_d\left( {\bf x}, \theta \right) \right]^{\rm T} \in \mathbb{R}^d$ and the Dirichlet boundary condition ${\bm u}_D\left( {{\bf{x}},\theta } \right) = \left[ u_{D,1}\left( {\bf x}, \theta \right), \cdots, u_{D,d}\left( {\bf x}, \theta \right) \right]^{\rm T} \in \mathbb{R}^d$.
	In this paper, we only consider linear elastic stochastic problems.
    The following linear stochastic strain tensor and linear elastic constitutive relation are adopted
		
	\vspace{-0.5cm}
	\begin{equation} \label{eq:sto_strain_stress}
        {\bm \varepsilon} \left( {\bm u}\left( {{\bf{x}},\theta } \right) \right) = \frac{1}{2}\left( {\nabla {\bm u}\left( {{\bf{x}},\theta } \right) + \left( {\nabla}{\bm u}\left( {{\bf{x}},\theta } \right) \right)^{\rm T} } \right) \in \mathbb{R}^{d \times d}, \quad {\bm \sigma} \left( {\bm u}\left( {{\bf{x}},\theta } \right),\theta \right) = {\bm C}\left( {\bf x}, \theta \right) {\bm \varepsilon} \left( {\bm u}\left( {{\bf{x}},\theta } \right) \right) \in \mathbb{R}^{d \times d},
	\end{equation}
    where ${\bm C}\left( {\bf x}, \theta \right)$ is a fourth order elastic tensor that may be related to stochastic material properties, e,g. stochastic Young's modulus and stochastic Poisson ratio.
    
    To solve \eqref{eq:EP_probelm}, let us consider its weak form written as follows: find a stochastic function ${\bm u}\left( {{\bf{x}},\theta } \right) \in {\mathscr V}: \Omega \times \Theta \rightarrow \mathbb{R}^d$ such that the following equation holds for ${\cal P}$-almost surely $\theta \in \Theta$,
		
	\vspace{-0.5cm}
	\begin{equation} \label{eq:weak_form}
	   {\mathscr W} \left( {\bm u}\left( {{\bf{x}},\theta } \right), {\bm v}\left( {\bf{x}} \right),\theta \right) = {\mathscr F} \left( {\bm v}\left( {\bf{x}} \right),\theta \right), \quad \forall {\bm v}\left( {\bf{x}} \right) \in {\mathscr V},
	\end{equation} 
	where the functional space is defined as ${\mathscr V} = \left\{ {\bm v} \in \left[ {\bm{\mathscr H}}^1\left( \Omega \right) \right]^d : {\bm v} = {\bf 0}~{\rm on}~{\Gamma}_D \right\}$, and ${\mathscr H}^1\left( \Omega \right)$ is the subspace of the space of square integrable scalar functions on $\Omega$ that contains both the function and its weak derivatives \cite{yosida2012functional}.
    The terms ${\mathscr W} \left( {\bm u}\left( {{\bf{x}},\theta } \right), {\bm v}\left( {\bf{x}} \right),\theta \right)$ and ${\mathscr F} \left( {\bm v}\left( {\bf{x}} \right),\theta \right)$ are given by

	\vspace{-0.5cm}
	\begin{align}
    	{\mathscr W} \left( {\bm u}\left( {{\bf{x}},\theta } \right), {\bm v}\left( {\bf{x}} \right),\theta \right) &= \int_\Omega {{\bm \sigma} \left( {{\bm u}\left( {{\bf{x}},\theta } \right),\theta } \right) \colon {\bm \varepsilon} \left( {{\bm v}\left( {\bf{x}} \right)} \right){\rm d}{\bf{x}}} \label{eq:weak_W_1} \\
    	&= \int_\Omega \left[ {\bm C}\left( {\bf x}, \theta \right) {\bm \varepsilon} \left( {{\bm u}\left( {{\bf{x}},\theta } \right) } \right) \right] \colon {\bm \varepsilon} \left( {{\bm v}\left( {\bf{x}} \right)} \right){\rm d}{\bf{x}}, \label{eq:weak_W_2} \\
    	{\mathscr F} \left( {\bm v}\left( {\bf{x}} \right),\theta \right) &= \int_\Omega  {{\bm f}\left( {{\bf{x}},\theta } \right) \cdot {\bm v}\left( {\bf{x}} \right){\rm d}{\bf{x}}}  + \int_{{\Gamma _N}} {{\bm g}\left( {{\bf{x}},\theta } \right) \cdot {\bm v}\left( {\bf{x}} \right){\rm d}{\bf{s}}}. \label{eq:weak_F}
	\end{align}

\subsection{Stochastic virtual element discretization}
    We adopt the stochastic virtual element discretization for the weak form \eqref{eq:weak_form}.
    Specifically, the domain $\Omega$ is partitioned into $n_e$ non-overlapping polygonal elements $\Omega = \overline {\bigcup\nolimits_{{\rm e} = 1}^{{n_e}} \Omega^{\left( {\rm e} \right)} }$, and each element $\Omega^{\left( {\rm e} \right)}$, ${\rm e} = 1, \dots, n_e$ includes $n^{\left( {\rm e} \right)}$ vertices and $m^{\left( {\rm e} \right)}$ edges.
    In this paper, we only consider the low-order virtual element, but the proposed method can be extended to higher-order virtual elements \cite{da2017high}.
    We give the following approximate discretized virtual space ${\mathscr V}_h\left( \Omega^{\left( {\rm e} \right)} \right) \subset {\mathscr V}$ of the element $\Omega^{\left( {\rm e} \right)}$
    
	\vspace{-0.5cm}
    \begin{align} \label{eq:Vh}
        {\mathscr V}_h\left( \Omega^{\left( {\rm e} \right)} \right) = \left\{ {\bm v}_h \in \left[ {\bm{\mathscr H}}^1\left( \Omega^{\left( {\rm e} \right)} \right) \bigcap {\bm{\mathscr C}}^0\left( \Omega^{\left( {\rm e} \right)} \right) \right]^d: {\bm v}_{h,\partial \Omega^{\left( {\rm e} \right)}_i} \in \left[ {\mathscr P}_1\left( \partial \Omega^{\left( {\rm e} \right)}_i \right) \right]^d, \forall \partial \Omega^{\left( {\rm e} \right)}_i \in \partial \Omega^{\left( {\rm e} \right)}, \right. \qquad \nonumber \\
        \left. \nabla \cdot \left[ {\bm C}\left( \theta \right) {\bm \varepsilon} \left( {\bm v}_h \right) \right] = {\bf 0}, ~\forall \theta \in \Theta ~~{\rm on}~~ \Omega^{\left( {\rm e} \right)} \right\},
    \end{align}
    where ${\mathscr P}_1$ represents the space of polynomials of degree up to 1, ${\bm v}_h$ is a polynomial on each edge $\partial \Omega^{\left( {\rm e} \right)}_i$, $i = 1, \cdots, m^{\left( {\rm e} \right)}$ of $\Omega^{\left( {\rm e} \right)}$ and ${\bm{\mathscr C}}^0$-continuity on the element $\Omega^{\left( {\rm e} \right)}$, and $\nabla \cdot \left[ {\bm C}\left( \theta \right) {\bm \varepsilon} \left( {\bm v}_h \right) \right]$ vanishes on the element $\Omega^{\left( {\rm e} \right)}$ for all $\theta \in \Theta$.
    In this way, we couple the classical virtual element space and the random input $\theta$.
    The function ${\bm v}_h$ is not known on the element $\Omega^{\left( {\rm e} \right)}$ but explicitly known on the edge $\partial \Omega^{\left( {\rm e} \right)}$.
    We now consider the weak form \eqref{eq:weak_form} on the discretized space: find ${\bm u}_h\left( {{\bf{x}},\theta } \right) \in {\mathscr V}_h: \Omega \times \Theta \rightarrow \mathbb{R}^d$ such that the following equation holds for ${\cal P}$-almost surely $\theta \in \Theta$,

	\vspace{-0.5cm}
	\begin{equation} \label{eq:weak_form_h}
	   {\mathscr W}_h \left( {\bm u}_h\left( {{\bf{x}},\theta } \right), {\bm v}_h\left( {\bf{x}} \right),\theta \right) = {\mathscr F}_h \left( {\bm v}_h\left( {\bf{x}} \right),\theta \right), \quad \forall {\bm v}_h\left( {\bf{x}} \right) \in {\mathscr V}_h,
	\end{equation}
    where ${\mathscr W}_h$ and ${\mathscr F}_h$ are assembled by looping through all elements $\left\{ \Omega^{\left( {\rm e} \right)} \right\}_{{\rm e} = 1}^{n_e}$

	\vspace{-0.5cm}
	\begin{equation} \label{eq:WF_e}
        {\mathscr W}_h \left( {\bm u}_h\left( {{\bf{x}},\theta } \right), {\bm v}_h\left( {\bf{x}} \right),\theta \right) = \sum\limits_{{\rm e} = 1}^{n_e} {\mathscr W}_h^{\left( {\rm e} \right)} \left( {\bm u}_h\left( {{\bf{x}},\theta } \right), {\bm v}_h\left( {\bf{x}} \right),\theta \right), \quad
        {\mathscr F}_h \left( {\bm v}_h\left( {\bf{x}} \right),\theta \right) = \sum\limits_{{\rm e} = 1}^{n_e} {\mathscr F}_h^{\left( {\rm e} \right)} \left( {\bm v}_h\left( {\bf{x}} \right),\theta \right).
	\end{equation}
    It is noted that ${\mathscr W}_h^{\left( {\rm e} \right)} \left( {\bm u}_h\left( {{\bf{x}},\theta } \right), {\bm v}_h\left( {\bf{x}} \right),\theta \right)$ and ${\mathscr F}_h^{\left( {\rm e} \right)} \left( {\bm v}_h\left( {\bf{x}} \right),\theta \right)$ cannot be evaluated in a similar way to the classical finite element method since the function ${\bm v}_h$ is unknown on the element $\Omega^{\left( {\rm e} \right)}$.
    To this end, a projection operator ${{\bf \Pi} ^{\left( {\rm e} \right)}}: {\mathscr V}_h\left( \Omega^{\left( {\rm e} \right)} \right) \rightarrow \left[ {\mathscr P}_1\left( \Omega^{\left( {\rm e} \right)} \right) \right]^d$ is defined similar to the deterministic VEM such that for $\forall {\bm v}_h \in {\mathscr V}_h\left( \Omega^{\left( {\rm e} \right)} \right)$, ${\bm p}_1 \in \left[ {\mathscr P}_1\left( \Omega^{\left( {\rm e} \right)} \right) \right]^d$, $\theta \in \Theta$,
    
	\vspace{-0.5cm} 
    \begin{equation} \label{eq:VEM_ortho}
        {\mathscr W}_h^{\left( {\rm e} \right)} \left( {\bm v}_h\left( {\bf{x}} \right) - {{\bf \Pi} ^{\left( {\rm e} \right)}}{\bm v}_h\left( {\bf{x}} \right), {\bm p}_1\left( {\bf x} \right),\theta \right) = 0
    \end{equation}
    holds.
    We can thus reformulate the term ${\mathscr W}_h^{\left( {\rm e} \right)} \left( {\bm u}_h\left( {{\bf{x}},\theta } \right), {\bm v}_h\left( {\bf{x}} \right),\theta \right)$ as
    
	\vspace{-0.5cm}
	\begin{align}
    	&{\mathscr W}_h^{\left( {\rm e} \right)} \left( {\bm u}_h\left( {{\bf{x}},\theta } \right), {\bm v}_h\left( {\bf{x}} \right),\theta \right) \nonumber \\
        = &{\mathscr W}_h^{\left( {\rm e} \right)} \left( {{\bf \Pi} ^{\left( {\rm e} \right)}}{\bm u}_h\left( {{\bf{x}},\theta } \right) + {\bm u}_h\left( {{\bf{x}},\theta } \right) - {{\bf \Pi} ^{\left( {\rm e} \right)}}{\bm u}_h\left( {{\bf{x}},\theta } \right), {{\bf \Pi} ^{\left( {\rm e} \right)}}{\bm v}_h\left( {\bf{x}} \right) + {\bm v}_h\left( {\bf{x}} \right) - {{\bf \Pi} ^{\left( {\rm e} \right)}}{\bm v}_h\left( {\bf{x}} \right),\theta \right), \label{eq:weak_W_h_1} \\
    	= &\underbrace{{\mathscr W}_h^{\left( {\rm e} \right)} \left( {{\bf \Pi} ^{\left( {\rm e} \right)}}{\bm u}_h\left( {{\bf{x}},\theta } \right), {{\bf \Pi} ^{\left( {\rm e} \right)}}{\bm v}_h\left( {\bf{x}} \right), \theta \right)}_{={\mathscr W}_C^{\left( {\rm e} \right)}\left( \theta \right)}
        + \underbrace{{\mathscr W}_h^{\left( {\rm e} \right)} \left( {\bm u}_h\left( {{\bf{x}},\theta } \right) - {{\bf \Pi} ^{\left( {\rm e} \right)}}{\bm u}_h\left( {{\bf{x}},\theta } \right), {\bm v}_h\left( {\bf{x}} \right) - {{\bf \Pi} ^{\left( {\rm e} \right)}}{\bm v}_h\left( {\bf{x}} \right),\theta \right)}_{={\mathscr W}_S^{\left( {\rm e} \right)}\left( \theta \right)} \nonumber \\
        & + \underbrace{{\mathscr W}_h^{\left( {\rm e} \right)} \left( {{\bf \Pi} ^{\left( {\rm e} \right)}}{\bm u}_h\left( {{\bf{x}},\theta } \right), {\bm v}_h\left( {\bf{x}} \right) - {{\bf \Pi} ^{\left( {\rm e} \right)}}{\bm v}_h\left( {\bf{x}} \right),\theta \right)}_{=0} 
        + \underbrace{{\mathscr W}_h^{\left( {\rm e} \right)} \left( {\bm u}_h\left( {{\bf{x}},\theta } \right) - {{\bf \Pi} ^{\left( {\rm e} \right)}}{\bm u}_h\left( {{\bf{x}},\theta } \right), {{\bf \Pi} ^{\left( {\rm e} \right)}}{\bm v}_h\left( {\bf{x}} \right),\theta \right)}_{=0}, \label{eq:weak_W_h_2}
	\end{align}
    where ${\mathscr W}_C^{\left( {\rm e} \right)}\left( \theta \right)$ and ${\mathscr W}_S^{\left( {\rm e} \right)}\left( \theta \right)$ are called the stochastic consistent term and the stochastic stabilization term, respectively.
    The last two terms are obtained according to  \eqref{eq:VEM_ortho} since ${{\bf \Pi} ^{\left( {\rm e} \right)}}v_h\left( {\bf{x}} \right)$, ${{\bf \Pi} ^{\left( {\rm e} \right)}}{\bm u}_h\left( {{\bf{x}},\theta } \right) \in \left[ {\mathscr P}_1\left( \Omega^{\left( {\rm e} \right)} \right) \right]^d$.
    Furthermore, ${\bm \varepsilon}^{\left( {\rm e} \right)} \left( {{\bf \Pi} ^{\left( {\rm e} \right)}}{\bm u}_h\left( {{\bf{x}},\theta } \right) \right)$ and ${\bm \varepsilon}^{\left( {\rm e} \right)} \left( {{{\bf \Pi} ^{\left( {\rm e} \right)}}{\bm v}_h\left( {\bf{x}} \right)} \right)$ are (stochastic) constant strain tensors.
    Therefore, the stochastic consistent term ${\mathscr W}_C^{\left( {\rm e} \right)}\left( \theta \right)$ is evaluated via
    
	\vspace{-0.5cm} 
    \begin{align}\label{eq:W_C_theta}
        {\mathscr W}_C^{\left( {\rm e} \right)}\left( \theta \right) &= \int_{\Omega^{\left( {\rm e} \right)}} \left[ {\bm C}^{\left( {\rm e} \right)}\left( {\bf x}, \theta \right) {\bm \varepsilon}^{\left( {\rm e} \right)} \left( {{\bf \Pi} ^{\left( {\rm e} \right)}}{\bm u}_h\left( {{\bf{x}},\theta } \right) \right) \right] \colon {\bm \varepsilon}^{\left( {\rm e} \right)} \left( {{{\bf \Pi} ^{\left( {\rm e} \right)}}{\bm v}_h\left( {\bf{x}} \right)} \right){\rm d}{\bf{x}} \nonumber \\
        &= a^{\left( {\rm e} \right)} \left[ {\bm C}^{\left( {\rm e} \right)}\left( {\bf x}, \theta \right) {\bm \varepsilon}^{\left( {\rm e} \right)} \left( {{\bf \Pi} ^{\left( {\rm e} \right)}}{\bm u}_h\left( {{\bf{x}},\theta } \right) \right) \right] \colon {\bm \varepsilon}^{\left( {\rm e} \right)} \left( {{{\bf \Pi} ^{\left( {\rm e} \right)}}{\bm v}_h\left( {\bf{x}} \right)} \right),
    \end{align}
    where $a^{\left( {\rm e} \right)}$ is the area (for 2D polygonal element) or the volume (for 3D polygonal element) of the element $\Omega^{\left( {\rm e} \right)}$.

    Further, we let the stochastic solution ${\bm u}_h\left( {\bf x},\theta \right)$ and the function $v_h\left( {\bf x} \right)$ on the element $\Omega^{\left( {\rm e} \right)}$ be approximated using a set of virtual basis functions $\left\{ \varphi_i\left( {\bf x} \right) \in {\mathscr V}_h\left( \Omega^{\left( {\rm e} \right)} \right) \right\}_{i=1}^{n^{\left( {\rm e} \right)}}$
    
	\vspace{-0.5cm} 
    \begin{equation}\label{eq:basis_uv}
        {\bm u}_h\left( {\bf x},\theta \right) = \sum\limits_{i=1}^{n^{\left( {\rm e} \right)}} \varphi_i\left( {\bf x} \right) {\bf u}_i^{\left( {\rm e} \right)}\left( \theta \right), \quad
        {\bm v}_h\left( {\bf x} \right) = \sum\limits_{i=1}^{n^{\left( {\rm e} \right)}} \varphi_i\left( {\bf x} \right) {\bf v}_i^{\left( {\rm e} \right)},
    \end{equation}
    where the solution vector ${\bf u}_i^{\left( {\rm e} \right)}\left( \theta \right)$ of the element $\Omega^{\left( {\rm e} \right)}$ is ${\bf u}^{\left( {\rm e} \right)}\left( \theta \right) = \left[ {\bf u}_1^{\left( {\rm e} \right){\rm T}}\left( \theta \right), \cdots, {\bf u}_{n^{\left( {\rm e} \right)}}^{\left( {\rm e} \right){\rm T}}\left( \theta \right) \right]^{\rm T} \in \mathbb{R}^{n^{\left( {\rm e} \right)}d}$,
    and ${\bf u}_i^{\left( {\rm e} \right)}\left( \theta \right) = \left[ u_{i,1}^{\left( {\rm e} \right)}\left( \theta \right), \cdots, u_{i,d}^{\left( {\rm e} \right)}\left( \theta \right) \right]^{\rm T} \in \mathbb{R}^d$ is the solution vector of the $i$-th vertex.
    The vector ${\bf v}_i^{\left( {\rm e} \right)} \in \mathbb{R}^{n^{\left( {\rm e} \right)}d}$ has a similar expression but does not involve the random input $\theta$.
    Applying the projection operator ${{\bf \Pi} ^{\left( {\rm e} \right)}}$ to ${\bm u}_h\left( {\bf x},\theta \right)$ and ${\bm v}_h\left( {\bf x} \right)$ we have
    
	\vspace{-0.5cm} 
    \begin{equation}\label{eq:basis_uv_II}
        {{\bf \Pi} ^{\left( {\rm e} \right)}}{\bm u}_h\left( {\bf x},\theta \right) = \sum\limits_{i=1}^{n^{\left( {\rm e} \right)}} {{\bf \Pi} ^{\left( {\rm e} \right)}}\varphi_i\left( {\bf x} \right) {\bf u}_i^{\left( {\rm e} \right)}\left( \theta \right), \quad
        {{\bf \Pi} ^{\left( {\rm e} \right)}}{\bm v}_h\left( {\bf x} \right) = \sum\limits_{i=1}^{n^{\left( {\rm e} \right)}} {{\bf \Pi} ^{\left( {\rm e} \right)}}\varphi_i\left( {\bf x} \right) {\bf v}_i^{\left( {\rm e} \right)}.
    \end{equation}
    On the basis of this we rewrite the strain tensors ${\bm \varepsilon}^{\left( {\rm e} \right)} \left( {{\bf \Pi} ^{\left( {\rm e} \right)}}{\bm u}_h\left( {{\bf{x}},\theta } \right) \right)$ and ${\bm \varepsilon}^{\left( {\rm e} \right)} \left( {{{\bf \Pi} ^{\left( {\rm e} \right)}}{\bm v}_h\left( {\bf{x}} \right)} \right)$ as the following vector forms
    
	\vspace{-0.5cm} 
    \begin{equation}\label{eq:e_u}
        {\bm \varepsilon}_{{\bf \Pi},{\rm vec}}^{\left( {\rm e} \right)} \left( {\bf u}^{\left( {\rm e} \right)}\left( \theta \right) \right) = {\bf B}^{\left( {\rm e} \right)} {\bf u}^{\left( {\rm e} \right)}\left( \theta \right) = \left[ {\bf B}_1^{\left( {\rm e} \right)}, \cdots, {\bf B}_{n^{\left( {\rm e} \right)}}^{\left( {\rm e} \right)} \right] {\bf u}^{\left( {\rm e} \right)}\left( \theta \right) \in \mathbb{R}^{\frac{d\left( d+1 \right)}{2}}, \quad
        {\bm \varepsilon}_{{\bf \Pi},{\rm vec}}^{\left( {\rm e} \right)} \left( {\bf v}^{\left( {\rm e} \right)} \right) = {\bf B}^{\left( {\rm e} \right)} {\bf v}^{\left( {\rm e} \right)} \in \mathbb{R}^{\frac{d\left( d+1 \right)}{2}},
    \end{equation}
    where the matrices ${\bf B}_i^{\left( {\rm e} \right)} \in \mathbb{R}^{\frac{d\left( d+1 \right)}{2} \times d}$, $i = 1, \cdots, n^{\left( {\rm e} \right)}$ are given by
    
	\vspace{-0.5cm} 
    \begin{equation} \label{eq:B_i}
        {\bf B}_i^{\left( {\rm e} \right)} = \left[ {\begin{array}{*{20}{c}}
        {\frac{{\partial {{\bf \Pi} ^{\left( {\rm e} \right)}} {\varphi _i}\left( {\bf{x}} \right)}}{{\partial x_1}}}&0\\
        0&{\frac{{\partial {{\bf \Pi} ^{\left( {\rm e} \right)}} {\varphi _i}\left( {\bf{x}} \right)}}{{\partial x_2}}}\\
        {\frac{{\partial {{\bf \Pi} ^{\left( {\rm e} \right)}} {\varphi _i}\left( {\bf{x}} \right)}}{{\partial x_2}}}&{\frac{{\partial {{\bf \Pi} ^{\left( {\rm e} \right)}} {\varphi _i}\left( {\bf{x}} \right)}}{{\partial x_1}}}
        \end{array}} \right]~(d=2)~~{\rm or}~~
        \left[ {\begin{array}{*{20}{c}}
        {\frac{{\partial {{\bf \Pi} ^{\left( {\rm e} \right)}} {\varphi _i}\left( {\bf{x}} \right)}}{{\partial x_1}}}&0&0\\
        0&{\frac{{\partial {{\bf \Pi} ^{\left( {\rm e} \right)}} {\varphi _i}\left( {\bf{x}} \right)}}{{\partial x_2}}}&0\\
        0&0&{\frac{{\partial {{\bf \Pi} ^{\left( {\rm e} \right)}} {\varphi _i}\left( {\bf{x}} \right)}}{{\partial x_3}}}\\
        {\frac{{\partial {{\bf \Pi} ^{\left( {\rm e} \right)}} {\varphi _i}\left( {\bf{x}} \right)}}{{\partial x_2}}}&{\frac{{\partial {{\bf \Pi} ^{\left( {\rm e} \right)}} {\varphi _i}\left( {\bf{x}} \right)}}{{\partial x_1}}}&0\\
        0&{\frac{{\partial {{\bf \Pi} ^{\left( {\rm e} \right)}} {\varphi _i}\left( {\bf{x}} \right)}}{{\partial x_3}}}&{\frac{{\partial {{\bf \Pi} ^{\left( {\rm e} \right)}} {\varphi _i}\left( {\bf{x}} \right)}}{{\partial x_2}}}\\
        {\frac{{\partial {{\bf \Pi} ^{\left( {\rm e} \right)}} {\varphi _i}\left( {\bf{x}} \right)}}{{\partial x_3}}}&0&{\frac{{\partial {{\bf \Pi} ^{\left( {\rm e} \right)}} {\varphi _i}\left( {\bf{x}} \right)}}{{\partial x_1}}}
        \end{array}} \right]~(d=3),
    \end{equation}
    where the components ${\frac{{\partial {{\bf \Pi} ^{\left( {\rm e} \right)}} {\varphi _i}\left( {\bf{x}} \right)}}{{\partial x_j}}}$, $i = 1, \cdots, n^{\left( {\rm e} \right)}$, $j = 1, \cdots, d$ are evaluated via transferring the calculations to edges and faces of the element $\Omega^{\left( {\rm e} \right)}$ based on \eqref{eq:VEM_ortho}, which is the same as the deterministic virtual element method, see \cite{beirao2013basic, beirao2014hitchhiker, gain2014virtual, artioli2017arbitrary} for details.
    In this way, it does not require knowing explicit representations of the functions $\left\{ \varphi_i\left( {\bf x} \right) \right\}_{i=1}^{n^{\left( {\rm e} \right)}}$ and only needs to know their traces on edges.
    The Lagrangian linear basis functions similar to those used in the classical finite element method can be adopted for the purpose.
    Substituting \eqref{eq:e_u} into \eqref{eq:W_C_theta} we have
    
	\vspace{-0.5cm} 
    \begin{equation}\label{eq:W_C_mat}
        {\mathscr W}_C^{\left( {\rm e} \right)}\left( \theta \right) = a^{\left( {\rm e} \right)} {\bf v}^{\left( {\rm e} \right){\rm T}} {\bf B}^{\left( {\rm e} \right){\rm T}} {\bf G}^{\left( {\rm e} \right)}\left( \theta \right) {\bf B}^{\left( {\rm e} \right)} {\bf u}^{\left( {\rm e} \right)}\left( \theta \right),
    \end{equation}
    where ${\bf G}^{\left( {\rm e} \right)}\left( \theta \right) \in \mathbb{R}^{\frac{d\left( d+1 \right)}{2} \times \frac{d\left( d+1 \right)}{2}}$ is the matrix form of the tensor ${\bm C}\left( {\bf x}, \theta \right)$ of the element $\Omega^{\left( {\rm e} \right)}$.
    Two detailed representations of ${\bf G}\left( \theta \right)$ for 2D and 3D problems can be found in the numerical example section.
    Hence, the stochastic element stiffness matrix corresponding to the stochastic consistent term ${\mathscr W}_C^{\left( {\rm e} \right)}\left( \theta \right)$ is given by
    
	\vspace{-0.5cm} 
    \begin{equation}\label{eq:k_C_ele}
        {\bf k}_C^{\left( {\rm e} \right)}\left( \theta \right) = a^{\left( {\rm e} \right)} {\bf B}^{\left( {\rm e} \right){\rm T}} {\bf G}^{\left( {\rm e} \right)}\left( \theta \right) {\bf B}^{\left( {\rm e} \right)} \in \mathbb{R}^{n^{\left( {\rm e} \right)}d \times n^{\left( {\rm e} \right)}d}.
    \end{equation}
    
    Further, let us consider the calculation of the stochastic stabilization term ${\mathscr W}_S^{\left( {\rm e} \right)}\left( \theta \right)$ in \eqref{eq:weak_W_h_2}, which can be achieved taking advantage of several numerical strategies \cite{beirao2013basic, gain2014virtual, artioli2017arbitrary}.
    Here we adopt the approach presented in \cite{beirao2014hitchhiker, chi2020virtual}, which corresponds to

	\vspace{-0.5cm} 
    \begin{align}\label{eq:W_S_0}
        {\mathscr W}_S^{\left( {\rm e} \right)}\left( \theta \right)
        = \gamma_S ^{\left( {\rm e} \right)}\left( \theta \right) \sum\limits_{i=1}^{n^{\left( {\rm e} \right)}} \left[ {\bm u}_h\left( {{\bf{x}}_i,\theta } \right) - {{\bf \Pi} ^{\left( {\rm e} \right)}}{\bm u}_h\left( {{\bf{x}}_i,\theta } \right) \right] \cdot \left[ {\bm v}_h\left( {\bf{x}}_i \right) - {{\bf \Pi} ^{\left( {\rm e} \right)}}{\bm v}_h\left( {\bf{x}}_i \right) \right],
    \end{align}
    where ${\bf x}_i = \left( x_{i,1}, \cdots, x_{i,d} \right)$, $x_{i,j}$ denotes the $i$-th coordinate value of the $j$-th vertex of the element $\Omega^{\left( {\rm e} \right)}$.
    In the practical implementation, it is calculated as \cite{chi2020virtual}
    
	\vspace{-0.5cm} 
    \begin{align}\label{eq:W_S_theta}
        {\mathscr W}_S^{\left( {\rm e} \right)}\left( \theta \right)
        = \gamma_S ^{\left( {\rm e} \right)}\left( \theta \right) {\bf v}^{\left( {\rm e} \right){\rm T}} 
        \left( {\bf I}_{n^{\left( {\rm e} \right)}d} - {\bf S}^ {\left( {\rm e} \right)} \right)^{\rm T} \left( {\bf I}_{n^{\left( {\rm e} \right)}d} - {\bf S}^ {\left( {\rm e} \right)} \right)
        {\bf u}^{\left( {\rm e} \right)}\left( \theta \right),
    \end{align}
    where the coefficient $\gamma_S ^{\left( {\rm e} \right)}\left( \theta \right) = \frac{1}{2}{\rm Tr} \left( {\bf G}^{\left( {\rm e} \right)}\left( \theta \right) \right)$, ${\rm{Tr}}\left( \cdot \right)$ is the trace operator of matrices, ${\bf I}_{n^{\left( {\rm e} \right)}d} \in \mathbb{R}^{n^{\left( {\rm e} \right)}d \times n^{\left( {\rm e} \right)}d}$ is the identity matrix,
    and the deterministic matrix ${\bf S}^ {\left( {\rm e} \right)} \in \mathbb{R}^{n^{\left( {\rm e} \right)}d \times n^{\left( {\rm e} \right)}d}$ is given by
    
	\vspace{-0.5cm} 
    \begin{equation}\label{eq:S_e}
        {\bf S}^ {\left( {\rm e} \right)} = \left[ {\begin{array}{*{20}{c}}
        h_{11}^ {\left( {\rm e} \right)} {\bf I}_d & \cdots & h_{1n^{\left( {\rm e} \right)}}^ {\left( {\rm e} \right)} {\bf I}_d \\
        \vdots & \ddots & \vdots \\
        h_{n^{\left( {\rm e} \right)}1}^ {\left( {\rm e} \right)} {\bf I}_d & \cdots & h_{n^{\left( {\rm e} \right)}n^{\left( {\rm e} \right)}}^ {\left( {\rm e} \right)} {\bf I}_d
        \end{array}} \right], \quad
        h_{ij}^ {\left( {\rm e} \right)} = {\bf X}_i^{\left( {\rm e} \right){\rm T}} {\bf A}_j^{\left( {\rm e} \right)} + \frac{1}{n^{\left( {\rm e} \right)}},
    \end{equation}
    where the vectors ${\bf X}_i^{\left( {\rm e} \right)}$ and ${\bf A}_i^{\left( {\rm e} \right)}$ are given by
    
	\vspace{-0.5cm} 
    \begin{equation}\label{eq:X_A}
        {\bf X}_i^{\left( {\rm e} \right)} = \left[ x_{1,i} - \frac{1}{n^{\left( {\rm e} \right)}} \sum\limits_{j=1}^{n^{\left( {\rm e} \right)}}x_{1,j}, \cdots, x_{d,i} - \frac{1}{n^{\left( {\rm e} \right)}} \sum\limits_{j=1}^{n^{\left( {\rm e} \right)}}x_{d,j} \right]^{\rm T} \in \mathbb{R}^d, \quad
        {\bf A}_i^{\left( {\rm e} \right)} = \left[ {\frac{{\partial {\bf \Pi} {\varphi _i}\left( {\bf{x}} \right)}}{{\partial x_1}}}, \cdots, {\frac{{\partial {\bf \Pi} {\varphi _i}\left( {\bf{x}} \right)}}{{\partial x_d}}} \right]^{\rm T} \in \mathbb{R}^d.
    \end{equation}
    Hence, the stochastic element stiffness matrix corresponding to the stochastic stabilization term ${\mathscr W}_S^{\left( {\rm e} \right)}\left( \theta \right)$ is given by
    
	\vspace{-0.5cm} 
    \begin{align} \label{eq:k_S_ele}
        {\bf k}_S^{\left( {\rm e} \right)}\left( \theta \right) = \gamma_S ^{\left( {\rm e} \right)}\left( \theta \right)
        \left( {\bf I}_{n^{\left( {\rm e} \right)}d} - {\bf S}^ {\left( {\rm e} \right)} \right)^{\rm T} \left( {\bf I}_{n^{\left( {\rm e} \right)}d} - {\bf S}^ {\left( {\rm e} \right)} \right) \in \mathbb{R}^{n^{\left( {\rm e} \right)}d \times n^{\left( {\rm e} \right)}d},
    \end{align}
    which is very close to that in the deterministic VEM, but the coefficient $\gamma_S ^{\left( {\rm e} \right)}\left( \theta \right)$ involves the random input $\theta$.
    The calculation procedure of deterministic virtual element matrices can thus be inherited for ${\bf k}_S^{\left( {\rm e} \right)}\left( \theta \right)$.
    We only need to pay a little attention to the calculation of the random coefficient $\gamma_S ^{\left( {\rm e} \right)}\left( \theta \right)$.

    In the last step, let us consider the stochastic term ${\mathscr F}_h^{\left( {\rm e} \right)} \left( {\bm v}_h\left( {\bf{x}} \right),\theta \right)$ calculated exactly using the one-point integration rule on the edges and the face (for 2D case) or the faces and the element (for 3D case)
    
    \vspace{-0.5cm} 
    \begin{align}\label{eq:F_weak}
        {\mathscr F}_h^{\left( {\rm e} \right)} \left( {\bm v}_h\left( {\bf{x}} \right),\theta \right) &= \int_{\Omega^{\left( {\rm e} \right)}}  {f^{\left( {\rm e} \right)}\left( {{\bf{x}},\theta } \right) \cdot {\bm v}_h\left( {\bf{x}} \right){\rm d}{\bf{x}}}  + \int_{\Gamma _N^{\left( {\rm e} \right)}} {g^{\left( {\rm e} \right)}\left( {{\bf{x}},\theta } \right) \cdot {\bm v}_h\left( {\bf{x}} \right){\rm d}{\bf{s}}} \nonumber \\
        &= \frac{a^{\left( {\rm e} \right)}}{\left( n^{\left( {\rm e} \right)} \right)^2} {\bf v}^{\left( {\rm e} \right){\rm T}} {\bf Z}_1^{\rm T}{\bf Z}_1 {\bm f}^{\left( {\rm e} \right)}\left( \theta \right) +
        {\bf v}^{\left( {\rm e} \right){\rm T}} \sum\limits_{j=1}^{n_{\Gamma_N}^{\left( {\rm e} \right)}} \frac{b_j^{\left( {\rm e} \right)}}{ \left( n_{\Gamma_{N,j}}^{\left( {\rm e} \right)} \right)^2 } {\bf Z}_{2,j}^{\rm T} {\bf Z}_{2,j} {\bm g}^{\left( {\rm e} \right)}\left( \theta \right),
    \end{align}
    where $n_{\Gamma_N}^{\left( {\rm e} \right)}$ is the number of edges (for 2D case) or faces (for 3D case) of the element $\Omega^{\left( {\rm e} \right)}$, 
    $\Gamma_{N,j}^{\left( {\rm e} \right)}$ is the $j$-th edge or face and includes $n_{\Gamma_{N,j}}^{\left( {\rm e} \right)}$ vertices,
    and $b_j^{\left( {\rm e} \right)}$ is the length (for 2D element) or the area (for 3D element) of $\Gamma_{N,j}^{\left( {\rm e} \right)}$.
    The vectors consisting of the values on each vertex are ${\bm f}^{\left( {\rm e} \right)}\left( \theta \right) = \left[ {\bm f}_1^{\left( {\rm e} \right){\rm T}}\left( \theta \right), \cdots, {\bm f}_{n^{\left( {\rm e} \right)}}^{\left( {\rm e} \right){\rm T}}\left( \theta \right) \right]^{\rm T} \in \mathbb{R}^{n^{\left( {\rm e} \right)}d}$,
    ${\bm f}_i^{\left( {\rm e} \right)}\left( \theta \right) = \left[ f_{i,1}^{\left( {\rm e} \right)}\left( \theta \right), \cdots, f_{i,d}^{\left( {\rm e} \right)}\left( \theta \right) \right]^{\rm T} \in \mathbb{R}^d$ and
    ${\bm g}^{\left( {\rm e} \right)}\left( \theta \right) = \left[ {\bm g}_1^{\left( {\rm e} \right){\rm T}}\left( \theta \right), \cdots, {\bm g}_{n^{\left( {\rm e} \right)}}^{\left( {\rm e} \right){\rm T}}\left( \theta \right) \right]^{\rm T} \in \mathbb{R}^{n^{\left( {\rm e} \right)}d}$,
    ${\bm g}_i^{\left( {\rm e} \right)}\left( \theta \right) = \left[ g_{i,1}^{\left( {\rm e} \right)}\left( \theta \right), \cdots, g_{i,d}^{\left( {\rm e} \right)}\left( \theta \right) \right]^{\rm T} \in \mathbb{R}^d$.
    The deterministic matrices are given by
    ${\bf Z}_1 = \left[ \underbrace{ {\bf I}_d, \cdots, {\bf I}_d }_{n^{\left( {\rm e} \right)}} \right] \in \mathbb{R}^{d \times n^{\left( {\rm e} \right)}d}$ and
    ${\bf Z}_{2,j} = \left[ \delta_{1,\Gamma_{N,j}^{\left( {\rm e} \right)}}^* {\bf I}_d, \cdots, \delta_{n^{\left( {\rm e} \right)},\Gamma_{N,j}^{\left( {\rm e} \right)}}^* {\bf I}_d \right] \in \mathbb{R}^{d \times n^{\left( {\rm e} \right)}d}$, where $\delta_{i,\Gamma_{N,j}^{\left( {\rm e} \right)}}^* = 1$ if the vertex $i \in \Gamma_{N,j}^{\left( {\rm e} \right)}$ and 0 otherwise.
    Thus, the stochastic element force vector is calculated as

    \vspace{-0.5cm} 
    \begin{equation}\label{eq:f_e}
        {\bf f}^{\left( {\rm e} \right)} \left( \theta \right) = \frac{a^{\left( {\rm e} \right)}}{\left( n^{\left( {\rm e} \right)} \right)^2} {\bf Z}_1^{\rm T}{\bf Z}_1 {\bm f}^{\left( {\rm e} \right)}\left( \theta \right) +
        \sum\limits_{j=1}^{n_{\Gamma_N}^{\left( {\rm e} \right)}} \frac{b^{\left( {\rm e} \right)}}{ \left( n_{\Gamma_{N,j}}^{\left( {\rm e} \right)} \right)^2 } {\bf Z}_{2,j}^{\rm T} {\bf Z}_{2,j} {\bm g}^{\left( {\rm e} \right)}\left( \theta \right) \in \mathbb{R}^{n^{\left( {\rm e} \right)}d}.
    \end{equation}

    Assembling the above stochastic element stiffness matrices and stochastic element force vector we obtain the following SVEE

    \vspace{-0.5cm}
    \begin{equation}\label{eq:SVEE}
        {\bf K}\left( \theta \right) {\bf u}\left( \theta \right) = {\bf F}\left( \theta \right),
    \end{equation}
    where the global stochastic stiffness matrix ${\bf K}\left( \theta \right)$ and the the global stochastic force vector ${\bf F}\left( \theta \right)$ are assembled by

    \vspace{-0.5cm}
    \begin{equation}\label{eq:SVEE_ass}
        {\bf K}\left( \theta \right) = \bigcup \limits_{{\rm e} = 1}^{n_e} \left( {\bf k}_C^{\left( {\rm e} \right)}\left( \theta \right) + {\bf k}_S^{\left( {\rm e} \right)}\left( \theta \right) \right) \in \mathbb{R}^{n \times n}, \quad
        {\bf F}\left( \theta \right) = \bigcup \limits_{{\rm e} = 1}^{n_e} {\bf f}^{\left( {\rm e} \right)}\left( \theta \right) \in \mathbb{R}^n,
    \end{equation}
    where $\bigcup \left( \cdot \right)$ represents the assembly operator for all stochastic element matrices and vectors, and the total degree of freedom (DoF) is given by $n = n_ed$.

    Further, the stochastic matrix ${\bf G}\left( {\bf x}, \theta \right)$ (or the tensor ${\bm C}\left( {\bf x}, \theta \right)$) can be approximated using the following separated form in many cases
    
	\vspace{-0.5cm}
	\begin{equation} \label{eq:G_i}
        {\bf G}\left( {\bf x}, \theta \right) = \sum\limits_{i=0}^m \xi_i\left( \theta \right) {\bf G}_i,
	\end{equation}
    where $\xi_0\left( \theta \right) \equiv 1$, $\left\{ \xi_i\left( \theta \right) \right\}_{i=1}^m$ are scalar random variables,
    $\left\{ {\bf G}_i \in \mathbb{R}^{\frac{d\left( d+1 \right)}{2} \times \frac{d\left( d+1 \right)}{2}} \right\}_{i=0}^m$ are deterministic matrices.
    For the non-separated stochastic matrix ${\bf G}\left( {\bf x}, \theta \right)$, the approaches for simulating random fields can be adopted to achieve \eqref{eq:G_i}-like approximations for both Gaussian and non-Gaussian random inputs, e.g. Karhunen-Lo{\`e}ve expansion and Polynomial Chaos expansion \cite{sakamoto2002polynomial, zheng2017simulation, zheng2021sample}.
    In this way, we can reformulate the SVEE (\ref{eq:SVEE}) as

    \vspace{-0.5cm}
    \begin{equation} \label{eq:SVEE_decom}
        \left( \sum\limits_{i=0}^m \xi_i\left( \theta \right) {\bf K}_i \right) {\bf u}\left( \theta \right) = {\bf F}\left( \theta \right),
    \end{equation}
    where the deterministic matrices $\left\{ {\bf K}_i \right\}_{i=0}^m$ are assembled via

    \vspace{-0.5cm}
    \begin{equation} \label{eq:SVEE_decom_Ki}
        {\bf K}_i = \bigcup \limits_{{\rm e}=1}^{n_e} \left( {\bf k}_C^{\left( {\rm e} \right)}\left( {\bf G}_i \right) + {\bf k}_S^{\left( {\rm e} \right)}\left( {\bf G}_i \right) \right) \in \mathbb{R}^{n \times n}
    \end{equation}
    since the stochastic element matrices ${\bf k}_C^{\left( {\rm e} \right)}\left( \theta \right)$ in \eqref{eq:k_C_ele} and ${\bf k}_S^{\left( {\rm e} \right)}\left( \theta \right)$ in \eqref{eq:k_S_ele} depend linearly on the matrix components $\left\{ {\bf G}_i \right\}_{i=0}^m$.

\section{PC-SVEM: Polynomial Chaos based spectral stochastic virtual element method} \label{sec:PC_SVEM}
    The PC-based methods have been well developed and widely applied to solve a variety of stochastic problems.
    In this section, we present a PC-SVEM to solve the SVEE (\ref{eq:SVEE}) (or the separated form \eqref{eq:SVEE_decom}).
    In this method, the stochastic solution ${\bf u} \left( \theta \right)$ is expanded using (generalized) PC basis as follows

    \vspace{-0.5cm}
    \begin{equation} \label{eq:sol_PC}
        {\bf u}_{{\rm PC},k} \left( \theta \right) = \sum\limits_{i=1}^k \Gamma_i \left( \theta \right) {\bf d}_{{\rm PC},i} , 
    \end{equation}
    where $\left\{ {\bf d}_{{\rm PC},i} \right\}_{i=1}^k$ are the corresponding deterministic vectors to be solved, $\left\{ \Gamma_i \left( \theta \right) \right\}_{i=1}^k$ are the PC basis.
    In practice, we can choose different PC basis for different types of random inputs, such as the Hermite PC basis for Gaussian random variables and the Legendre PC basis for uniform random variables \cite{ghanem2003stochastic, xiu2002wiener}.
    The stochastic Galerkin approach is then used to transform SVEE (\ref{eq:SVEE}) into the following deterministic equation \cite{ghanem2003stochastic}

    \vspace{-0.5cm}
    \begin{equation}\label{eq:PC_equ}
        \int_{\Theta} \left[ {\bf K}\left( \theta \right) \sum\limits_{i=1}^k \Gamma_i \left( \theta \right) {\bf d}_{{\rm PC},i} - {\bf F}\left( \theta \right) \right] \Gamma_j \left( \theta \right) {\rm d} {\cal P}\left( \theta \right) = 0, \quad j = 1, \cdots, k,
    \end{equation}
    where ${\cal P}\left( \theta \right)$ is the probability measurement of the random input $\theta$.
    Further, the above equation can be rewritten as a compact form

    \vspace{-0.5cm}
    \begin{equation} \label{eq:Sys_PC}
        {\bf K}_{\rm PC} {\bf d}_{\rm PC} = {\bf F}_{\rm PC},
    \end{equation}
    where the augmented deterministic matrix ${\bf K}_{\rm PC} \in \mathbb{R}^{nk \times nk}$ and the augmented deterministic vectors ${\bf d}_{\rm PC}$, ${\bf F}_{\rm PC} \in \mathbb{R}^{nk}$ are assembled by
    
    \vspace{-0.5cm}
    \begin{equation} \label{eq:PC_KUF}
        {\bf K}_{\rm PC} = \left[ {\begin{array}{*{20}{c}}
        {\bf K}_{{\rm PC},11} & \cdots  &  {\bf K}_{{\rm PC},1k} \\
          \vdots & \ddots & \vdots \\
          {\bf K}_{{\rm PC},k1} & \cdots  & {\bf K}_{{\rm PC},kk}
        \end{array}} \right], \quad
        {\bf d}_{\rm PC} = \left[ {\begin{array}{*{20}{c}}
        {\bf d}_{{\rm PC},1}\\
          \vdots\\
          {\bf d}_{{\rm PC},k}
        \end{array}} \right],
        \quad
        {\bf F}_{\rm PC} = \left[ {\begin{array}{*{20}{c}}
        {\bf F}_{{\rm PC},1}\\
          \vdots\\
          {\bf F}_{{\rm PC},k}
        \end{array}} \right],
    \end{equation}
    where the matrix and vector components ${\bf K}_{{\rm PC},ji} \in \mathbb{R}^{n \times n}$ and ${\bf F}_{{\rm PC},j} \in \mathbb{R}^n$, $i,j = 1, \cdots, k$ are
    
    \vspace{-0.5cm}
    \begin{equation}\label{eq:PC_KF_ij}
        {\bf K}_{{\rm PC},ji} = \int_{\Theta} {\bf K}\left( \theta \right) \Gamma_i \left( \theta \right) \Gamma_j \left( \theta \right) {\rm d} {\cal P}\left( \theta \right), \quad
        {\bf F}_{{\rm PC},j} = \int_{\Theta} {\bf F}\left( \theta \right) \Gamma_j \left( \theta \right) {\rm d} {\cal P}\left( \theta \right).
    \end{equation}
    Further, if the separated form \eqref{eq:SVEE_decom} is considered, the above calculation of the matrices ${\bf K}_{{\rm PC},ji}$ is simplified as 
        
    \vspace{-0.5cm}
    \begin{equation}\label{eq:PC_K_ij}
        {\bf K}_{{\rm PC},ji} = \sum\limits_{l=0}^m \left[ \int_{\Theta} \xi_l\left( \theta \right) \Gamma_i \left( \theta \right) \Gamma_j \left( \theta \right) {\rm d} {\cal P}\left( \theta \right) \right] {\bf K}_i,
    \end{equation}
    which only involves the numerical integration $\int_{\Theta} \xi_l\left( \theta \right) \Gamma_i \left( \theta \right) \Gamma_j \left( \theta \right) {\rm d} {\cal P}\left( \theta \right)$.
    For low-dimensional stochastic problems, the calculation is cheap enough benefiting from efficient numerical integration strategies on stochastic spaces \cite{xiu2010numerical}.
    
    It is noted that the total number of PC basis is $k = \frac{\left(m+r\right)!}{m!r!}$, where $\left( \cdot \right)!$ is the factorial operator, $r$ is the expansion order of the PC basis.
    Thus, similar to classical PC-based methods, the proposed PC-SVEM still suffers from the curse of dimensionality since the matrix/vector size $nk$ in \eqref{eq:Sys_PC} increases sharply as the spatial DoF $n$ of the physical model, the stochastic dimension $m$ of the random input and the expansion order $r$ of the PC basis increase.
    For instance, the size is about $nk = 1 \times 10^6$ when $n = 1 \times 10^3$, $m = 10$ and $r=4$, which requires extremely expensive computational effort.
    Although several methods are developed to alleviate the computational burden, e.g. dedicated iterative algorithms and sparse PC approximations \cite{pellissetti2000iterative, keese2005hierarchical, blatman2010adaptive}.
    It is still challenging to solve very high-dimensional stochastic problems using the PC-SVEM.

\section{WIN-SVEM: weakly intrusive stochastic virtual element method} \label{sec:WIN_SVEM}
\subsection{A weakly intrusive stochastic virtual element method}
    To avoid the curse of dimensionality arising in the above PC-SVEM, we present a WIN-SVEM in this section, which can be considered as an extension of our previous work \cite{zheng2022weak, zheng2023stochastic} on stochastic finite element methods to SVEM.
    To this end, we consider the stochastic solution ${\bf u} \left( \theta \right)$ approximated by the following series expansion

    \vspace{-0.5cm}
    \begin{equation} \label{eq:sol_WIN}
        {\bf u}_{{\rm WIN},k} \left( \theta \right) = \sum\limits_{i=1}^k \lambda_i \left( \theta \right) {\bf d}_{{\rm WIN},i} , 
    \end{equation}
    where $\left\{ \lambda_i \left( \theta \right) \in \mathbb{R} \right\}_{i=1}^k$ are scalar random variables, $\left\{ {\bf d}_{{\rm WIN},i} \in \mathbb{R}^n \right\}_{i=1}^k$ are deterministic vectors, and $k$ is the number of retained terms.
    It is noted that the number $k$ and all pairs $\left\{ \lambda_i\left(\theta\right), {\bf d}_{{\rm WIN},i} \right\}_{i=1}^k$ are not known a priori.
    An iterative algorithm is presented to solve the pairs $\left\{ \lambda_i\left(\theta\right), {\bf d}_{{\rm WIN},i} \right\}$ one by one.
    Specifically, we assume that the $\left( k-1 \right)$-th approximation ${\bf u}_{{\rm WIN},k-1} \left( \theta \right) = \sum_{i=1}^{k-1} \lambda_i \left( \theta \right) {\bf d}_{{\rm WIN},i}$ has been known and the goal is to solve the $k$-th pair $\left\{ \lambda_k\left(\theta\right), {\bf d}_{{\rm WIN},k} \right\}$.
    The original SVEE (\ref{eq:SVEE}) can be rewritten as 

    \vspace{-0.5cm}
    \begin{equation} \label{eq:SVEE_k}
        {\bf K}\left( \theta \right) \lambda_k \left( \theta \right) {\bf d}_{{\rm WIN},k} = {\bf F}_k\left( \theta \right),
    \end{equation}
    where the stochastic vector ${\bf F}_k\left( \theta \right) = {\bf F}\left( \theta \right) - {\bf K}\left( \theta \right) \sum_{i=1}^{k-1} \lambda_i \left( \theta \right) {\bf d}_{{\rm WIN},i}$.
    In this way, \eqref{eq:SVEE_k} only involves one unknown pair $\left\{ \lambda_k\left(\theta\right), {\bf d}_{{\rm WIN},k} \right\}$.
    However, different from that the random basis (i.e. PC basis) has been known in PC-SVEM, both the random variable $\lambda_k\left(\theta\right)$ and the deterministic vector ${\bf d}_{{\rm WIN},k}$ are unknown in this case.
    To avoid this issue, an alternating iteration is adopted to solve them.
    Specifically, if the random variable $\lambda_k\left(\theta\right)$ has been known (or given an initial value), \eqref{eq:SVEE_k} is transformed into the following deterministic virtual element equation by taking advantage of the stochastic Galerkin procedure \cite{ghanem2003stochastic} similar to that in \eqref{eq:PC_equ}

    \vspace{-0.5cm}
    \begin{equation} \label{eq:WIN_d}
        \int_{\Theta} \left[ {\bf K}\left( \theta \right) \lambda_k \left( \theta \right) {\bf d}_{{\rm WIN},k} - {\bf F}_k\left( \theta \right) \right] \lambda_k \left( \theta \right) {\rm d} {\cal P}\left( \theta \right) = 0,
    \end{equation}
    which is equivalent to 

    \vspace{-0.5cm}
    \begin{equation} \label{eq:WIN_d_k}
        {\bf K}_{{\rm WIN},k}{\bf d}_{{\rm WIN},k} = {\bf F}_{{\rm WIN},k},
    \end{equation}
    where the deterministic matrix ${\bf K}_{{\rm WIN},k} = \int_{\Theta} {\bf K}\left( \theta \right) \lambda_k^2 \left( \theta \right) {\rm d} {\cal P}\left( \theta \right) \in \mathbb{R}^{n \times n}$ and the deterministic vector ${\bf F}_{{\rm WIN},k} = \int_{\Theta} {\bf F}_k\left( \theta \right) \lambda_k \left( \theta \right) {\rm d} {\cal P}\left( \theta \right) \in \mathbb{R}^n$.
    Existing numerical solvers can be adopted to solve it efficiently and accurately \cite{young2014iterative}.
    Note that the size of \eqref{eq:WIN_d_k} is the same as the original SVEE (\ref{eq:SVEE}), which is different from the augmented size of PC-based derived equation (\ref{eq:Sys_PC}) and can thus save a lot of computational effort.
    In practical implementations, we let the vector ${\bf d}_{{\rm WIN},k}$ orthogonal to the obtained vectors $\left\{ {\bf d}_{{\rm WIN},i} \right\}_{i=1}^{k-1}$ to speed up the convergence, which is achieved by using the Gram-Schmidt orthonormalization
	
	\vspace{-0.5cm}
	\begin{equation} \label{eq:d_or}
        {\bf d}_{{\rm WIN},k} = {\bf d}_{{\rm WIN},k} - \sum_{i=1}^{k-1}\left( {\bf d}_{{\rm WIN},k}^{\rm T} {\bf d}_{{\rm WIN},i} \right) {\bf d}_{{\rm WIN},i}, \quad {\bf d}_{{\rm WIN},k}^{\rm T} {\bf d}_{{\rm WIN},k} = 1,
	\end{equation}
	where $\left\{ {\bf d}_{{\rm WIN},i} \right\}_{i=1}^{k-1}$ are normalized orthogonal vectors that meet ${\bf d}_{{\rm WIN},i}^{\rm T} {\bf d}_{{\rm WIN},j}=\delta_{ij}$, where $\delta_{ij}$ is the Kronecker delta.

    With the deterministic vector ${\bf d}_{{\rm WIN},k}$ solved using \eqref{eq:WIN_d_k}, the random variable $\lambda_k\left( \theta \right)$ is then recalculated taking advantage of the following classical Galerkin procedure

    \vspace{-0.5cm}
    \begin{equation} \label{eq:WIN_lam}
        {\bf d}_{{\rm WIN},k}^{\rm T} \left[ {\bf K}\left( \theta \right) \lambda_k \left( \theta \right) {\bf d}_{{\rm WIN},k} - {\bf F}_k\left( \theta \right) \right] = 0.
    \end{equation}
    Since the stochastic matrix ${\bf K}\left( \theta \right)$ is positive definite and ${\bf z}^{\rm T} {\bf K}\left( \theta \right) {\bf z} > 0$, $\forall {\bf z} \ne {\bf 0} \in \mathbb{R}^n$, $\forall \theta \in \Theta$ holds, the above equation can be rewritten as

    \vspace{-0.5cm}
    \begin{equation} \label{eq:WIN_lam_k}
        \lambda_k \left( \theta \right) = \frac{ {\bf d}_{{\rm WIN},k}^{\rm T}{\bf F}_k\left( \theta \right) }{ {\bf d}_{{\rm WIN},k}^{\rm T}{\bf K}\left( \theta \right){\bf d}_{{\rm WIN},k} }.
    \end{equation}
    To avoid the curse of dimensionality arising in the high-dimensional problems, we adopt a non-intrusive sampling approach \cite{zheng2023stochastic} to solve \eqref{eq:WIN_lam_k} instead of the PC-based approximation, which corresponds to

    \vspace{-0.5cm}
    \begin{equation} \label{eq:WIN_lam_k_sample}
        \lambda_k \left( {\widehat{\bm \theta}} \right) = \left[ {\bf d}_{{\rm WIN},k}^{\rm T}{\bf F}_k\left( {\widehat{\bm \theta}} \right) \right] \oslash \left[ {\bf d}_{{\rm WIN},k}^{\rm T}{\bf K}\left( {\widehat{\bm \theta}} \right){\bf d}_{{\rm WIN},k} \right] \in \mathbb{R}^{n_s},
    \end{equation}
    where ${\Box}\left( {\widehat{\bm \theta}} \right)$ represents $n_s$ sample realizations of ${\Box}\left( \theta \right)$, 
    ${\bf d}_{{\rm WIN},k}^{\rm T}{\bf F}_k\left( {\widehat{\bm \theta}} \right) \in \mathbb{R}^{n_s}$ and ${\bf d}_{{\rm WIN},k}^{\rm T}{\bf K}\left( {\widehat{\bm \theta}} \right){\bf d}_{{\rm WIN},k} \in \mathbb{R}^{n_s}$ are the sample vectors of the random variables ${\bf d}_{{\rm WIN},k}^{\rm T}{\bf F}_k\left( \theta \right)$ and ${\bf d}_{{\rm WIN},k}^{\rm T}{\bf K}\left( \theta \right){\bf d}_{{\rm WIN},k}$, respectively, and $\oslash$ represents the element-wise division of two sample vectors, also known as Hadamard division operator.
    In this way, all random inputs are embedded into the sample realization vectors ${\bf d}_{{\rm WIN},k}^{\rm T}{\bf F}_k\left( {\widehat{\bm \theta}} \right)$ and ${\bf d}_{{\rm WIN},k}^{\rm T}{\bf K}\left( {\widehat{\bm \theta}} \right){\bf d}_{{\rm WIN},k}$, which is insensitive to the stochastic dimension of random inputs.
    The curse of dimensionality can thus be avoided successfully, which will be discussed in detail in the next section.

    We can calculate the $k$-th pair $\left\{ \lambda_k\left(\theta\right), {\bf d}_{{\rm WIN},k} \right\}$ by performing the iterative process of \eqref{eq:WIN_d_k} and \eqref{eq:WIN_lam_k_sample} until reaching a specified precision.
    A similar iteration is also adopted to calculate other pairs $\left\{ \lambda_{k+1}\left(\theta\right), {\bf d}_{{\rm WIN},k+1} \right\}$, $\cdots$ until a good approximation of the stochastic solution is achieved.
    However, it is noted that the stochastic solution ${\bf u}_{{\rm WIN},k} \left( \theta \right)$ in \eqref{eq:sol_WIN} is approximated in a sequential way and it does not exactly fulfill the original SVEE (\ref{eq:SVEE}).
    The approximation has low accuracy for some cases  \cite{zheng2022weak}.
    We introduce a recalculation process to avoid this problem.
    To this end, ${\bf D}_{\rm WIN} = \left[ {\bf d}_{{\rm WIN},1}, \cdots, {\bf d}_{{\rm WIN},k} \right] \in \mathbb{R}^{n \times k}$ is considered as a set of reduced basis functions and the random variable vector ${\bf \Lambda}\left( \theta \right) = \left[ \lambda_1 \left( \theta \right), \cdots, \lambda_k \left( \theta \right) \right]^{\rm T} \in \mathbb{R}^k$ is recalculated via the following reduced-order stochastic equation

    \vspace{-0.5cm}
    \begin{equation} \label{eq:SVEE_D}
        \left[ {\bf D}_{\rm WIN}^{\rm T}{\bf K}\left( \theta \right) {\bf D}_{\rm WIN} \right] {\bf \Lambda} \left( \theta \right) = {\bf D}_{\rm WIN} {\bf F}\left( \theta \right),
    \end{equation}
    which requires to be solved repeatedly for $n_s$ sample realizations to get the final solution ${\bf \Lambda}\left( \theta^{\left( i \right)} \right)$, $i = 1, \cdots, n_s$, but only very low computational effort is involved since the sizes of the reduced-order stochastic matrix ${\bf D}_{\rm WIN}^{\rm T}{\bf K}\left( \theta  \right){\bf D}_{\rm WIN} \in {{\mathbb{R}}^{k \times k}}$ and the reduced-order stochastic vector ${\bf D}_{\rm WIN}^{\rm T}{\bf F}\left( \theta  \right) \in {{\mathbb{R}}^k}$ are greatly reduced compared to the original SVEE (\ref{eq:SVEE}) in most cases.

    Let us highlight the weak intrusiveness of the proposed method.
    On one hand, \eqref{eq:sol_WIN} is considered as a kind of intrusive approximation of the stochastic solution that is very similar to the PC-based intrusive approximation (\ref{eq:sol_PC}).
    However, on the other hand, the implementation for solving ${\bf d}_{{\rm WIN},k}$ in \eqref{eq:WIN_d_k} only involves deterministic calculations and the matrix ${\bf K}_{{\rm WIN},k}$ keeps the same size and matrix properties as the original stochastic matrix ${\bf K}\left( \theta \right)$, which is weakly intrusive.
    Also, \eqref{eq:WIN_lam_k_sample} for calculating the random variable $\lambda_k\left( \theta \right)$ is fully non-intrusive.
    In these senses, we implement the intrusive stochastic solution approximation only in weakly intrusive and fully non-intrusive ways.
    The method combines the high efficiency of intrusive methods and the weak dimensionality dependence of non-intrusive methods.
    It can solve high-dimensional stochastic problems efficiently and accurately.

\subsection{High-dimensional stochastic problems}
    In this section, we will show that the proposed WIN-SVEM can be applied to high-dimensional stochastic problems without any modification.
    We explain this point from the perspective of the influence of high stochastic dimensions on solving \eqref{eq:WIN_d_k} and \eqref{eq:WIN_lam_k_sample}.
    We only consider \eqref{eq:SVEE_decom} in this section and
    a large number $m$ is truncated in \eqref{eq:SVEE_decom} to generate a high-dimensional stochastic problem.
    In this way, the deterministic matrix ${\bf K}_{{\rm WIN},k}$ and the deterministic vector ${\bf F}_{{\rm WIN},k}$ in \eqref{eq:WIN_d_k} are calculated via

    \vspace{-0.5cm}
    \begin{align} \label{eq:HD_d}
        {\bf K}_{{\rm WIN},k} &= \sum\limits_{i=0}^m \left[ \int_{\Theta} \xi_i\left( \theta \right) \lambda_k^2 \left( \theta \right) {\rm d} {\cal P}\left( \theta \right) \right] {\bf K}_i, \\
        {\bf F}_{{\rm WIN},k} &= \int_{\Theta} {\bf F}\left( \theta \right) \lambda_k \left( \theta \right) {\rm d} {\cal P}\left( \theta \right) - \sum\limits_{i=0}^m \sum\limits_{j=1}^{k-1} \left[ \int_{\Theta} \xi_i\left( \theta \right) \lambda_j \left( \theta \right) \lambda_k \left( \theta \right) {\rm d} {\cal P}\left( \theta \right) \right] {\bf K}_i {\bf d}_{{\rm WIN},j},
    \end{align}
    where the probability integrals are approximated using the following non-intrusive sampling approach
    
    \vspace{-0.5cm}
    \begin{equation} \label{eq:HD_E}
        \int_{\Theta} \xi_i\left( \theta \right) \lambda_j \left( \theta \right) \lambda_k \left( \theta \right) {\rm d} {\cal P}\left( \theta \right) = {\widehat{\mathbb{E}}}\left\{ \xi_i\left( {\widehat{\bm \theta}} \right) \odot \lambda_j \left( {\widehat{\bm \theta}} \right) \odot \lambda_k \left( {\widehat{\bm \theta}} \right) \right\} , \quad i = 0, \cdots, m, j = 1, \cdots, k,
    \end{equation}
    where $\left\{ \xi_i\left( {\widehat{\bm \theta}} \right) \in \mathbb{R}^{n_s} \right\}_{i = 0}^m$ are sample vectors of the random variables $\left\{ \xi_i\left( \theta \right) \right\}_{i = 0}^m$,
    the operator $\odot$ represents the element-wise multiplication of sample vectors, ${\widehat{\mathbb{E}}}\left\{ \cdot \right\}$ is the expectation operator of the sample vector.
    \eqref{eq:HD_E} takes a total of $k\left( m+1 \right)$ expectation operations, which is not sensitive to the stochastic dimension and has low computational effort even for very high stochastic dimensions.
    Note that although we only illustrate the high-dimensional input in the stochastic matrix ${\bf K}\left( \theta \right)$, the above calculation also works efficiently for $\int_{\Theta} {\bf F}\left( \theta \right) \lambda_k \left( \theta \right) {\rm d} {\cal P}\left( \theta \right)$ if the stochastic vector ${\bf F}\left( \theta \right)$ involves high-dimensional random inputs.
    
    Further, the sample vectors in right side of \eqref{eq:WIN_lam_k_sample} are calculated via

    \vspace{-0.5cm}
    \begin{align}\label{eq:HD_lam}
        {\bf d}_{{\rm WIN},k}^{\rm T}{\bf K}\left( {\widehat{\bm \theta}} \right){\bf d}_{{\rm WIN},k} &= \sum\limits_{i=0}^m \xi_i \left( {\widehat{\bm \theta}} \right) \left( {\bf d}_{{\rm WIN},k}^{\rm T}{\bf K}_i{\bf d}_{{\rm WIN},k} \right) \in \mathbb{R}^{n_s}, \\
        {\bf d}_{{\rm WIN},k}^{\rm T}{\bf F}_k\left( {\widehat{\bm \theta}} \right) &= {\bf d}_{{\rm WIN},k}^{\rm T} {\bf F}\left( \theta \right) - \sum\limits_{i=0}^m \sum\limits_{j=1}^{k-1} \left[ \xi_i\left( {\widehat{\bm \theta}} \right) \odot \lambda_j \left( {\widehat{\bm \theta}} \right) \right]  \left( {\bf d}_{{\rm WIN},k}^{\rm T}{\bf K}_i{\bf d}_{{\rm WIN},j} \right) \in \mathbb{R}^{n_s},
    \end{align}
    which requires a total of $k\left( m+1 \right)$ operations for ${\bf d}_{{\rm WIN},k}^{\rm T}{\bf K}_i{\bf d}_{{\rm WIN},j}$, $i = 0, \cdots, m$, $j = 1, \cdots, k$ and is also computationally cheap for high-dimensional stochastic problems.
    Therefore, both \eqref{eq:WIN_d_k} and \eqref{eq:WIN_lam_k_sample} are insensitive to the stochastic dimension.
    The proposed method can avoid the curse of dimensionality successfully.

\subsection{Algorithm implementation}

	The above proposed WIN-SVEM for solving SVEEs is summarized in Algorithm~\ref{alg_WIN_SVEM}, which includes two loop processes.
    The inner loop is from step \ref{Alg:step_04} to step \ref{Alg:step_09} and used to solve the $k$-th pair $\left\{ \lambda_k\left(\theta\right), {\bf d}_{{\rm WIN},k} \right\}$.
    The outer loop from step \ref{Alg:step_02} to step \ref{Alg:step_13} is to approximate the stochastic solution using a set of pairs $\left\{ \lambda_i\left(\theta\right), {\bf d}_{{\rm WIN},i} \right\}_{i=1}^k$.
    To execute the inner loop, a random sample vector ${\lambda_k ^{\left(0\right)}}\left( {\widehat{\bm{\theta}}} \right) \in \mathbb{R}^{n_s}$ is initialized in step \ref{Alg:step_03}.
    In the numerical implementation, any nonzero vectors of size $n_s$ can be chosen as the initialization since the initial samples have little influence on the computational accuracy and efficiency of the proposed method.
    Following each inner loop, we only need to update the stochastic force vector ${\bf F}_{k+1}\left( \theta \right)$ in step \ref{Alg:step_10} and store the reduced-order matrix ${\bf D}_{\rm WIN}$ in step \ref{Alg:step_11} in the outer loop.
    
	\begin{algorithm} [ht]
		\caption{WIN-SVEM for solving SVEEs} \label{alg_WIN_SVEM}
		\begin{algorithmic}[1]
			\State $k \leftarrow 1$
			\label{Alg:step_01}
			\While {$\epsilon _{{\bf u}, k} \ge \epsilon _{{\bf u}}$}
			\label{Alg:step_02}
			\State Initialize random samples ${\lambda_k ^{\left(0\right)}}\left( {\widehat{\bm{\theta}}} \right) = \left\{ {{\lambda_k ^{\left(0\right)}}\left( {{\theta ^{\left( i \right)}}} \right)} \right\}_{i = 1}^{{n_s}} \in \mathbb{R}^{n_s}$
			\label{Alg:step_03}
			\While {${\epsilon _{{\bf{d}},j}} > {\epsilon _{\bf{d}}}$}
			\label{Alg:step_04}
			\State Calculate the deterministic vector ${\bf d}_{{\rm WIN},k}^{\left( j \right)}$ by solving \eqref{eq:WIN_d_k}
			\label{Alg:step_05}
			\State Orthonormalize ${\bf{d}}_{{\rm WIN},k}^{\left( j \right)}$ using \eqref{eq:d_or}
			\label{Alg:step_06}
			\State Update the random sample vector $\lambda_k^{\left(j\right)}\left( {\widehat{\bm{\theta}}}\right) \in \mathbb{R}^{n_s}$ via \eqref{eq:WIN_lam_k_sample}
			\label{Alg:step_07}
			\State Compute the locally iterative error ${\epsilon _{{\bf{d}},j}}$, $j \leftarrow j+1$
            \label{Alg:step_08}
			\EndWhile
			\State {\bf{end}}
			\label{Alg:step_09}
			\State Update the stochastic force vector ${\bf F}_{k+1}\left( \theta \right) = {\bf F}_k\left( \theta \right) - {\bf K}\left( \theta \right) \lambda_k \left( \theta \right) {\bf d}_{{\rm WIN},k}$
			\label{Alg:step_10}
            \State Update the deterministic matrix ${\bf D}_{\rm WIN} = \left[ {\bf D}_{\rm WIN}, {\bf d}_{{\rm WIN},k} \right] \in \mathbb{R}^{n \times k}$
			\label{Alg:step_11}
			\State Compute the locally iterative error ${\epsilon _{{\bf{u}},k}}$, $k \leftarrow k + 1$
			\label{Alg:step_12}
			\EndWhile
			\State {\bf{end}}
			\label{Alg:step_13}
			\State Recalculate the random variable vector ${\bf{\Lambda}}\left( \theta \right) \in \mathbb{R}^k$ via \eqref{eq:SVEE_D}
			\label{Alg:step_14}
		\end{algorithmic}
	\end{algorithm}

    Two iterative criteria are involved in the above algorithm to check the convergence, i.e. ${\epsilon _{{\bf d},j}}$ in step \ref{Alg:step_08} for the inner loop and ${\epsilon _{{\bf u},k}}$ in step \ref{Alg:step_12} for the outer loop.
	The iterative error ${\epsilon _{{\bf d},j}}$ is defined as
	
	\vspace{-0.5cm}
	\begin{equation}\label{eq:loc_err_0}
        {\epsilon _{{\bf d},j}} = \frac{\left( {\bf{d}}^{\left( j \right)}_{{\rm WIN},k} - {\bf{d}}^{\left( j-1 \right)}_{{\rm WIN},k} \right)^{\rm T}\left( {\bf{d}}^{\left( j \right)}_{{\rm WIN},k} - {\bf{d}}^{\left( j-1 \right)}_{{\rm WIN},k} \right)}{{\bf{d}}^{\left( j \right){\rm T}}_{{\rm WIN},k} {\bf{d}}^{\left( j \right)}_{{\rm WIN},k}} = 2 - 2{\bf{d}}^{\left( j \right){\rm T}}_{{\rm WIN},k}{{\bf{d}}^{\left( j-1 \right)}_{{\rm WIN},k}},
	\end{equation}
	which measures the difference between the vectors ${\bf{d}}_{{\rm WIN},k}^{\left( j \right)}$ and ${\bf{d}}_{{\rm WIN},k}^{\left( j-1 \right)}$ and the calculation is stopped when ${\epsilon _{{\bf d},j}} < {\epsilon _{{\bf d}}}$ is met.
	Similarly, the iterative error ${\epsilon _{{\bf u},k}}$ is defined as
	
	\vspace{-0.5cm}
	\begin{align} \label{eq:In_It_err}
        {\epsilon _{{\bf u},k}} &= \frac{ \int_{\Theta} \left[ {{{\bf{u}}_{{\rm WIN},k}}\left( {\theta} \right) - {{\bf{u}}_{{\rm WIN},k-1}}\left( {\theta} \right)} \right]^{\rm T} \left[ {{{\bf{u}}_{{\rm WIN},k}}\left( {\theta} \right) - {{\bf{u}}_{{\rm WIN},k-1}}\left( {\theta} \right)} \right] {\rm d} {\cal P}\left( \theta \right) }{ \int_{\Theta} {{\bf{u}}_{{\rm WIN},k}^{\rm T}}\left( {\theta} \right) {{\bf{u}}_{{\rm WIN},k}}\left( {\theta} \right) {\rm d} {\cal P}\left( \theta \right) } \nonumber \\
        &= \frac{\int_{\Theta} \lambda_k^2\left(\theta\right) {\rm d} {\cal P}\left( \theta \right) {\bf d}_{{\rm WIN},k}^{\rm T} {\bf d}_{{\rm WIN},k} }{ \sum\limits_{i,j = 1}^k \int_{\Theta} \lambda_i\left(\theta\right)\lambda_j\left(\theta\right) {\rm d} {\cal P}\left( \theta \right) {\bf d}_{{\rm WIN},i}^{\rm T} {\bf d}_{{\rm WIN},j}} = \frac{\int_{\Theta} \lambda_k^2\left(\theta\right) {\rm d} {\cal P}\left( \theta \right)}{\sum\limits_{i = 1}^k \int_{\Theta} \lambda_i^2\left(\theta\right) {\rm d} {\cal P}\left( \theta \right)},
	\end{align}
    which measures the contribution of the $k$-th pair $\left\{ \lambda_k\left(\theta\right), {\bf d}_{{\rm WIN},k} \right\}$ to the stochastic solution ${\bf u}_{{\rm WIN},k}\left( \theta \right)$.
    However, \eqref{eq:In_It_err} may be not a good error indicator in many cases \cite{zheng2022weak} since the random variables $\left\{ \lambda_i\left( \theta \right) \right\}_{i=1}^k$ are calculated in a sequential way and $\int_{\Theta} \lambda_k^2\left(\theta\right) {\rm d} {\cal P}\left( \theta \right)$ may not keep decreasing.
    We avoid this problem by replacing $\left\{ \lambda_i\left( \theta \right) \right\}_{i=1}^k$ in \eqref{eq:In_It_err} with equivalent random variables $\left\{ {\widetilde \lambda}_i\left( \theta \right) \right\}_{i=1}^k$. 
    To this end, we calculate the autocorrelation function of the random variable vector ${\bf \Lambda}\left( \theta \right) = \left[ \lambda_1 \left( \theta \right), \cdots, \lambda_k \left( \theta \right) \right]^{\rm T}$ by
    
	\vspace{-0.5cm}
	\begin{equation}\label{eq:Lam_cov}
        {\bf C}_{{\bf{\Lambda}}{\bf{\Lambda}}} = {\widehat{\mathbb{E}}} \left\{ {\bf{\Lambda}}\left( {{\widehat{\bm{\theta}}} } \right) {\bf{\Lambda}}^{\rm T}\left( {{\widehat{\bm{\theta}}} } \right) \right\} \in \mathbb{R}^{k \times k},
	\end{equation}
	which is decomposed by eigendecomposition into
		
	\vspace{-0.5cm}
	\begin{equation}\label{eq:dom_cov}
	{\bf C}_{{\bf{\Lambda}}{\bf{\Lambda}}} = {\bf Q} {\bf Z} {\bf Q}^{\rm T},
	\end{equation}
	where ${\bf Q} \in \mathbb{R}^{k \times k}$ is an orthonormal matrix and ${\bf Z}$ is a diagonal matrix consisting of descending eigenvalues of the matrix ${\bf C}_{{\bf{\Lambda}}{\bf{\Lambda}}}$.
    We construct an equivalent random variable vector ${\bm{\widetilde \Lambda }}\left( {\theta } \right) = {\bf{Q}}^{\rm T} {\bm{\Lambda }}\left( {\theta } \right) = \left[ {{\widetilde \lambda }_1}\left( {\theta } \right), \cdots,{{\widetilde \lambda }_k}\left( {\theta } \right) \right]^{\rm T} \in \mathbb{R}^k$ whose autocorrelation function happens to be ${\bf C}_{{\bm{\widetilde \Lambda }}{\bm{\widetilde \Lambda }}} = {\bf Z}$.
    Substituting the equivalent random variables $\left\{ {\widetilde \lambda}_i\left( \theta \right) \right\}_{i=1}^k$ into \eqref{eq:In_It_err} we recalculate the iterative error $\epsilon_{{\bf u},k}$ as
		
	\vspace{-0.5cm}
	\begin{equation}\label{eq:Improve_It_err}
	\epsilon_{{\bf u},k} = \frac{\int_{\Theta} {\widetilde \lambda}_k^2\left(\theta\right) {\rm d} {\cal P}\left( \theta \right)}{\sum\limits_{i = 1}^k \int_{\Theta} {\widetilde \lambda}_i^2\left(\theta\right) {\rm d} {\cal P}\left( \theta \right)} = \frac{{\bf{Z}}_k}{{\rm{Tr}}\left( {\bf{Z}} \right)},  
	\end{equation}
    where ${\bf{Z}}_k$ is the $k$-th diagonal element of the matrix ${\bf Z}$.
    It is noted that the above reformulation does not improve the approximation accuracy of the stochastic solution and only provides an equivalent representation.
    In this way, the iterative error $\epsilon_{{\bf u},k}$ keeps decreasing as the retained item $k$ increases. 
    More details of the basic implementation and the comparison between \eqref{eq:In_It_err} and \eqref{eq:Improve_It_err} can be found in \cite{zheng2022weak}.

\section{Numerical examples} \label{sec:Ex}
	In this section, we test the proposed two methods with the aid of 2D and 3D numerical examples.
	For Algorithm \ref{alg_WIN_SVEM}, the convergence errors ${\epsilon _{{\bf d}}}$ in step \ref{Alg:step_02} for the inner loop and ${\epsilon _{{\bf u}}}$ in step \ref{Alg:step_04} for the outer loop are set as $1 \times 10^{-3}$ and $1 \times 10^{-6}$, respectively.
    For both examples, $n_s = 1 \times 10^4$ random samples are used for performing MC simulations and generating reference solutions.
    In our cases, $1 \times 10^4$ samples are enough to achieve convergent probabilistic solutions.
    The same $1 \times 10^4$ random samples are also used in the proposed WIN-SFEM to eliminate the influence caused by sampling processes.
    Further, the examples are performed on one core of a desktop computer (sixteen cores, Intel Core i7, 2.50GHz).

\subsection{Example 1: SVEM analysis of a 2D stochastic problem} \label{sec:ex_1}

\subsubsection{Simulation of random inputs} \label{sec:ex_1_RF}

 	\begin{figure}[!b]
		\centering
		\includegraphics[width=0.5\linewidth]{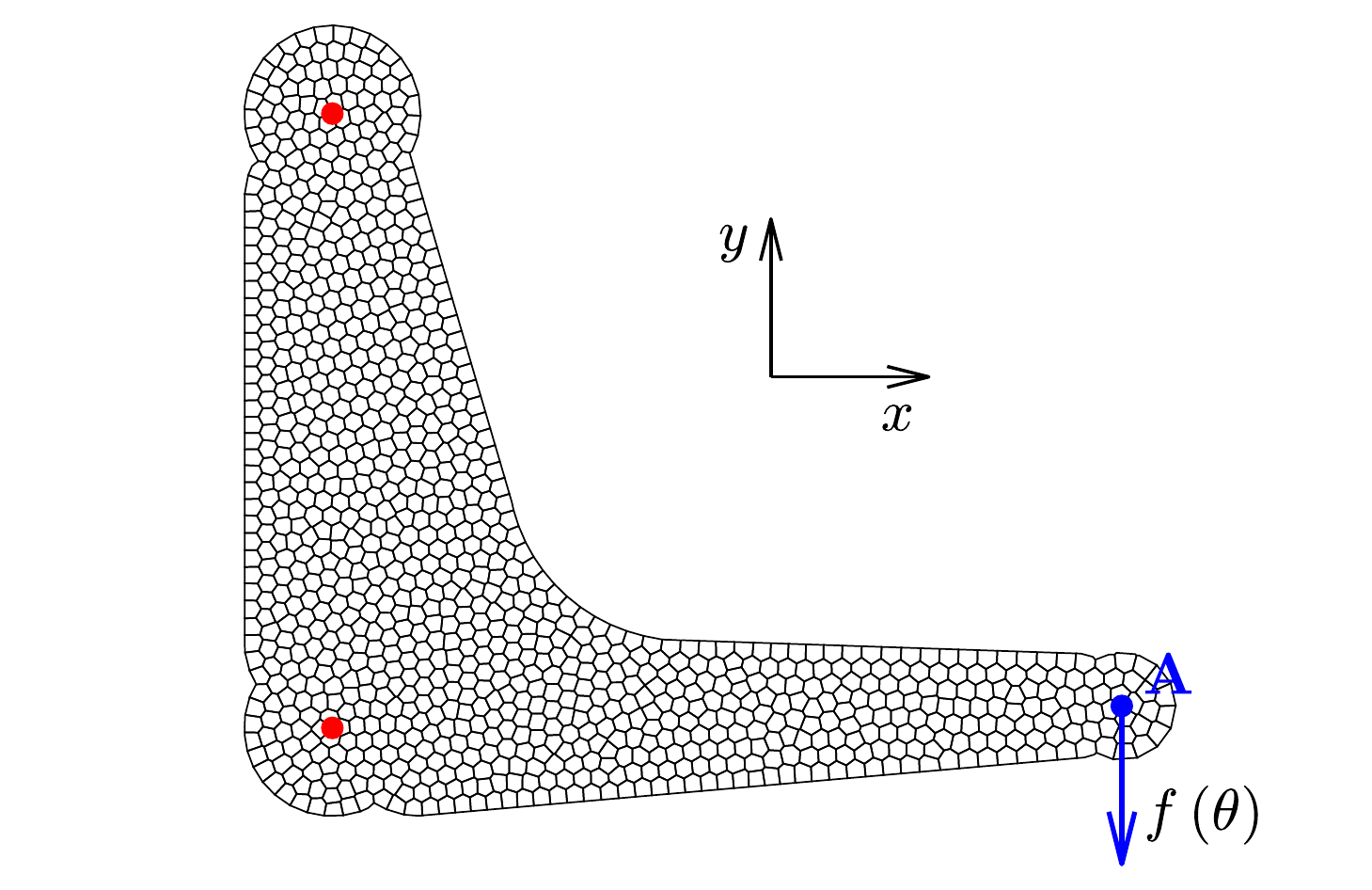}
		\caption{Geometry and Voronoi mesh of the 2D model.}\label{fig_e1_model}
	\end{figure}

    In this example, we consider the SVEM-based plane stress analysis of a 2D model shown in \figref{fig_e1_model} \cite{talischi2012polymesher}.
    The Voronoi mesh is adopted for the spatial discretization, including a total of 2018 vertices, 1000 Voronoi elements and 4036 DoFs.
    The model is fixed at the two red points as shown in \figref{fig_e1_model}.
    A stochastic force $f\left( \theta \right) = -1000 - 100 \xi_f\left( \theta \right)$ (unit: N) is applied to the blue point along the $y$ direction, where $\xi_f\left( \theta \right)$ is a standard Gaussian random variable.
    Further, the material property matrix ${\bf G}\left( x,y, \theta \right)$ is given by

	\vspace{-0.5cm}
	\begin{equation} \label{eq:Ex1_G}
        {\bf G}\left( x,y, \theta \right) = \frac{E\left( x,y, \theta \right)}{1 - \nu^2}
        \left[ {\begin{array}{*{20}{c}}
        1 & \nu & 0 \\
        \nu & 1 & 0\\
        0 & 0 & \frac{1}{2}\left( 1-\nu \right)
        \end{array}} \right] \in \mathbb{R}^{3 \times 3},
	\end{equation}
    where the Poisson ratio $\nu = 0.3$ and the Young's modulus $E\left( x,y, \theta \right)$ is modeled as a two-dimensional random field with the mean value $E_0\left( x,y \right) = 100$~MPa and the covariance function \cite{spanos2007karhunen}
    
    \vspace{-0.5cm}
    \begin{equation} \label{eq:Ex1_Cov}
        {\rm Cov}_{EE}\left(x_1,y_1;x_2,y_2\right) = \sigma_E^2\left( { 1 + \frac{{\left| x_1 - x_2 \right|}}{l_x} } \right) \left( { 1 + \frac{{\left| y_1 - y_2 \right|}}{l_y}} \right) \exp \left( { - \frac{{\left| x_1 - x_2 \right|}}{l_x} - \frac{{\left| y_1 - y_2 \right|}}{l_y}} \right),
    \end{equation}
    where the standard deviation $\sigma_E = 10$~MPa, and $l_x$ and $l_y$ are the correlation lengths in the $x$ and $y$ directions, respectively.
    The random field $E\left( x,y, \theta \right)$ is approximated via the following Karhunen-Lo{\`e}ve expansion \cite{ghanem2003stochastic, zheng2017simulation}
    
    \vspace{-0.5cm}
    \begin{equation} \label{eq:Ex1_E}
        E\left( x,y, \theta \right) = E_0\left( x,y \right) + \sum\limits_{i=1}^m \xi_i\left( \theta \right) \sqrt{\kappa_i} E_i\left( x,y \right),
    \end{equation}
    where $\left\{ \xi_i\left( \theta \right) \right\}_{i=1}^m$ are mutually independent standard Gaussian random variables and they are also independent of the random variable $\xi_f\left( \theta \right)$.
    $\left\{ \kappa_i, E_i\left( x,y \right) \right\}_{i=1}^m$ are eigenvalues and eigenvectors of the covariance function ${\rm Cov}_{EE}\left(x_1,y_1;x_2,y_2\right)$.
    They are solved by the following Fredholm integral equation of the second kind

    \vspace{-0.5cm}
    \begin{equation} \label{eq:Ex1_Cov_sol}
        \int_{\Omega} {\rm Cov}_{EE}\left(x_1,y_1; x_2,y_2\right) E_i\left( x_1,y_1 \right){\rm d}x_1{\rm d}y_1 = \kappa_i E_i\left( x_2,y_2 \right),
    \end{equation}
    which can be solved efficiently by taking advantage of existing eigenvalue solvers \cite{saad2011numerical}.
    To ensure $E\left( x,y, \theta \right) > 0$, the samples ${\theta}^{\left( i \right)}$ such that $\mathop{\min}\limits_{x,y \in \Omega}E\left( x,y,{\theta}^{\left( i \right)} \right) \le 1 \times 10^{-3}$ are dropped out in numerical implementations.
    In this way, the deterministic matrices $\left\{ {\bf G}_i\left( x,y \right) \right\}_{i=0}^m$ in \eqref{eq:G_i} are given by

	\vspace{-0.5cm}
	\begin{equation} \label{eq:Ex1_G_i}
        {\bf G}_i\left( x,y \right) = \frac{\sqrt{\kappa_i} E_i\left( x,y, \theta \right)}{1 - \nu^2}
        \left[ {\begin{array}{*{20}{c}}
        1 & \nu & 0 \\
        \nu & 1 & 0\\
        0 & 0 & \frac{1}{2}\left( 1-\nu \right)
        \end{array}} \right] \in \mathbb{R}^{3 \times 3}, ~~ {\rm where}~ \kappa_0 \equiv 1.
	\end{equation}

  	\begin{figure}[ht]
		\centering
		\includegraphics[width=0.8\linewidth]{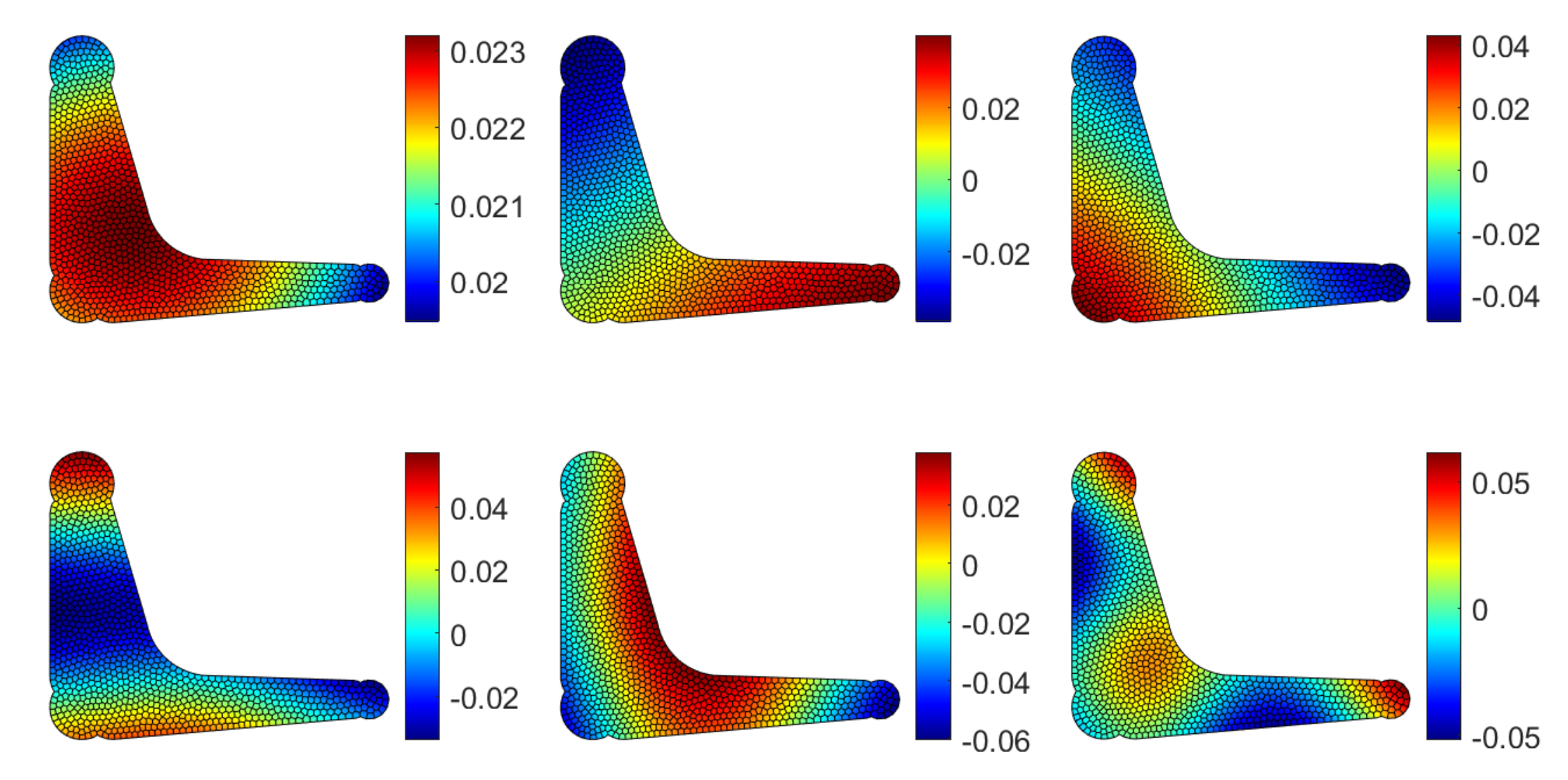}
		\caption{Eigenvectors $\left\{ E_i\left( x,y \right) \right\}_{i=1}^6$ of the covariance function ${\rm Cov}_{EE}\left(x_1,y_1; x_2,y_2\right)$ of the low-dimensional case.}\label{fig_e1_RF_ev}
	\end{figure}
 
 	\begin{figure}[!ht]
		\centering
		\includegraphics[width=0.5\linewidth]{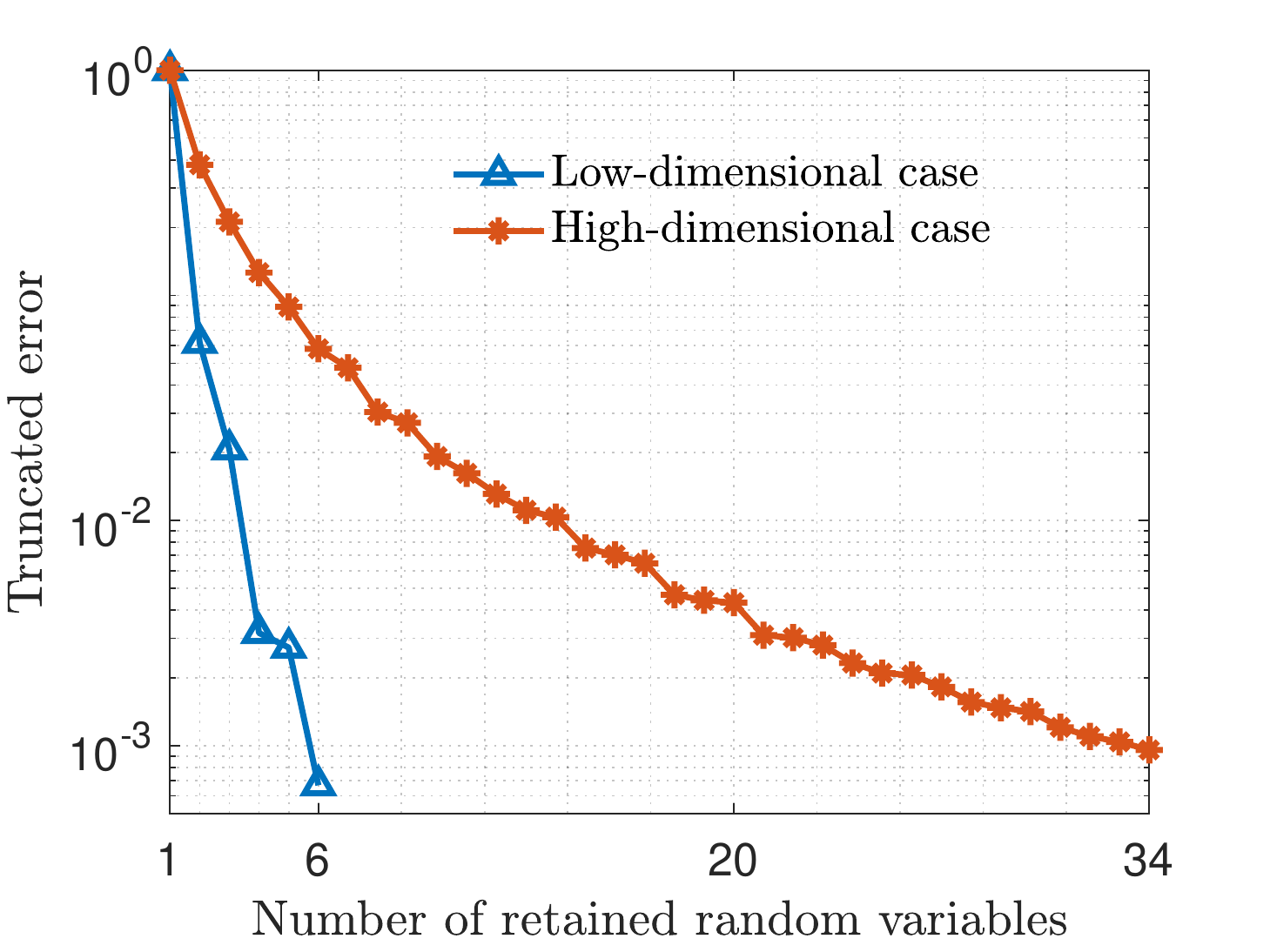}
		\caption{Truncated errors of different numbers of retained random variables.} \label{fig_e1_RF_err}
	\end{figure}

    In this example, we consider a low-dimensional case and a high-dimensional case via introducing different truncated items $m$ in \eqref{eq:Ex1_E}.
    For the low-dimensional case, we let the correlation lengths be $l_x = \max(x) - \min(x)$ and $l_y = \max(y) - \min(y)$.
    To achieve the truncated error ${\kappa_m} \mathord{\left/ {\vphantom{ \kappa_m \sum_{i=1}^m \kappa_i }} \right. \kern-\nulldelimiterspace} \sum_{i=1}^m \kappa_i < 1 \times 10^{-3}$, the number of truncated items is $m = 6$.
    Corresponding eigenvectors $\left\{ E_i\left( x,y \right) \right\}_{i=1}^6$ of the covariance function ${\rm Cov}_{EE}\left(x_1,y_1; x_2,y_2\right)$ are shown in \figref{fig_e1_RF_ev}.
    It is noted that we can solve \eqref{eq:Ex1_Cov_sol} using the discretized covariance matrix ${\bf Cov}_{EE} \in \mathbb{R}^{2018 \times 2018}$ of ${\rm Cov}_{EE}\left(x_1,y_1; x_2,y_2\right)$, which is only dependent of the vertices and independent of elements.
    In this way, the eigenvectors in \figref{fig_e1_RF_ev} are plotted on each vertex of the mesh.
    For the high-dimensional case, the correlation lengths are set as $l_x = \frac{1}{8}\left[ \max(x) - \min(x) \right]$ and $l_y = \frac{1}{8}\left[ \max(y) - \min(y) \right]$, and $m = 34$ truncated items are retained to achieve ${\kappa_{34}} \mathord{\left/ {\vphantom{ \kappa_{34} \sum_{i=1}^{34} \kappa_i }} \right. \kern-\nulldelimiterspace} \sum_{i=1}^{34} \kappa_i < 1 \times 10^{-3}$.
    The truncated errors of low- and high-dimensional cases are seen from \figref{fig_e1_RF_err}.
    Due to smaller correlation lengths in the high-dimensional case, it is slower to converge to the specified truncated error and more truncated items are required to capture the local correlation property.

\subsubsection{Low-dimensional case} \label{sec:ex_1_low_case}

    In this section, we solve the low-dimensional case using the proposed PC-SVEM and WIN-SVEM.
    For the PC-SVEM, the second order Hermite PC basis of seven standard Gaussian random variables ($\left\{ \xi_i\left( \theta \right) \right\}_{i=1}^6$ and $\xi_f\left( \theta \right)$) is adopted for $\left\{ \Gamma_i \left( \theta \right) \right\}_{i=1}^{36}$.
    The size of the augmented deterministic equation (\ref{eq:Sys_PC}) is 145296, which is much larger than the 4036 DoFs of the original stochastic problem.
    We do not perform the numerical implementations for higher order PC basis in this example. 
    On one hand, the second order PC basis is enough to achieve a good approximation of the stochastic solution.
    On the other hand, the size of the augmented equation (\ref{eq:Sys_PC}) is 484320 if the third order PC basis is adopted, which leads to too high computational burden in terms of storage and solution for a problem of such a small spatial scale.

 	\begin{figure}[ht]
		\centering
		\includegraphics[width=0.5\linewidth]{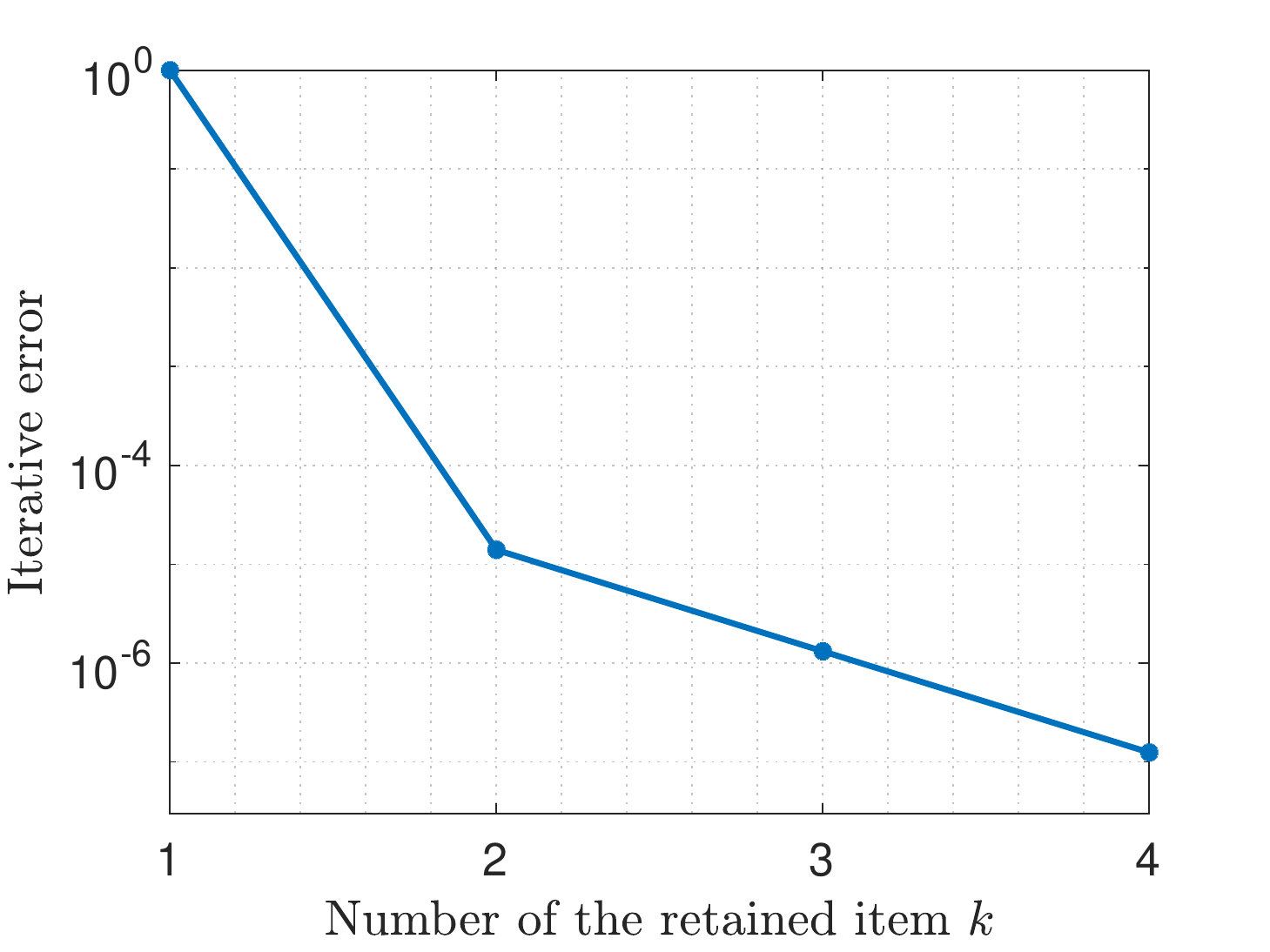}
		\caption{Iterative errors of different retained terms in the low-dimensional case.} \label{fig_e1_It_err_6}
	\end{figure}
  
\begin{figure}[ht]
		\centering
		\includegraphics[width=1.0\linewidth]{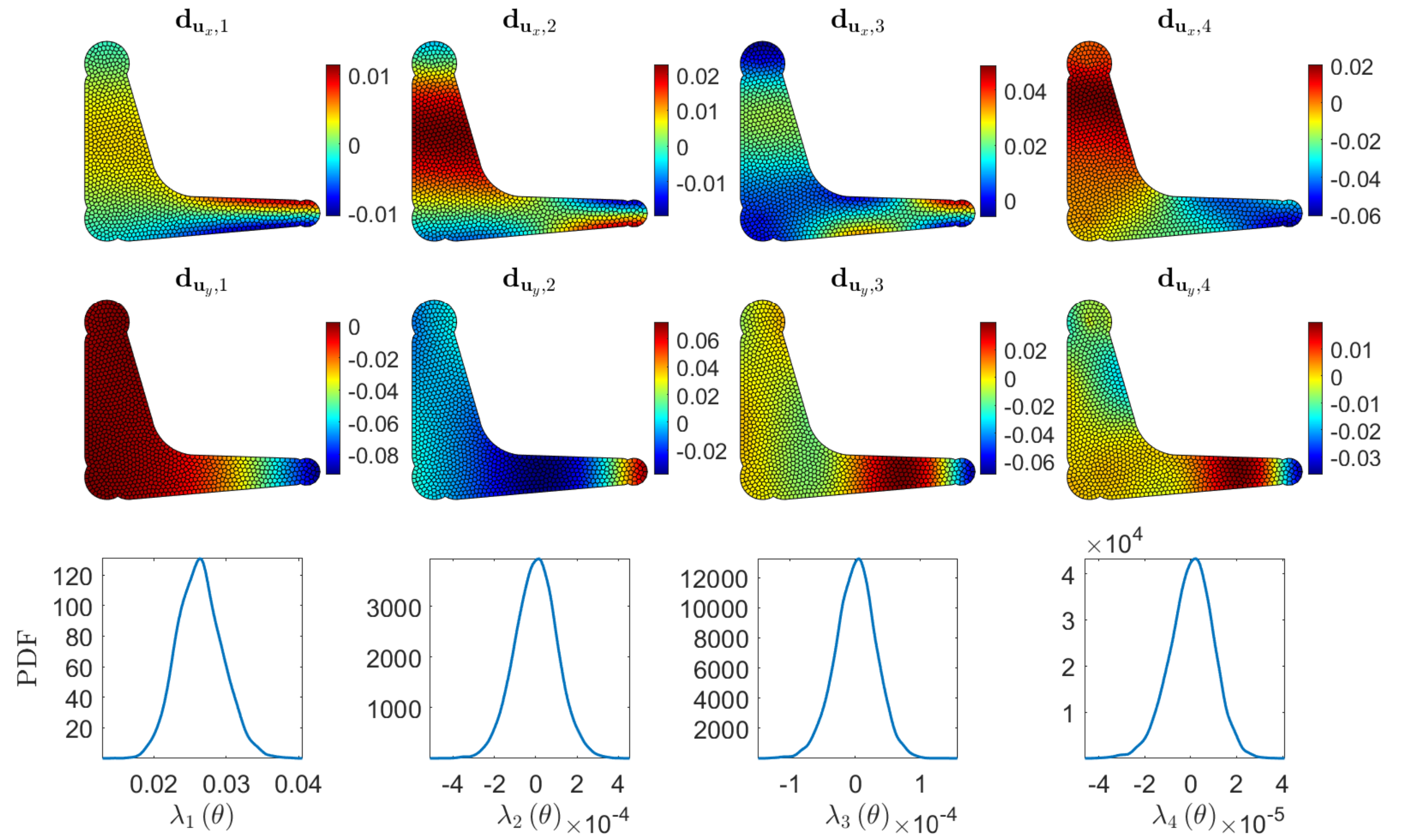}
		\caption{Components of the stochastic solution: the deterministic vectors $\left\{ {\bf d}_{{\bf u}_x,i} \right\}_{i=1}^4$ in the $x$ direction (the first line), the deterministic vectors $\left\{ {\bf d}_{{\bf u}_y,i} \right\}_{i=1}^4$ in the $y$ direction (the second line) and PDFs of the random variables $\left\{ \lambda_i\left( \theta \right) \right\}_{i=1}^4$ (the third line), respectively.}\label{fig_e1_DL_6}
	\end{figure}

    For the WIN-SVEM, Algorithm \ref{alg_WIN_SVEM} is performed and the iterative errors ${\epsilon _{{\bf u},k}}$ of different retained terms are shown in \figref{fig_e1_It_err_6}.
    Only four terms $k =4$ are retained for the stochastic solution approximation \eqref{eq:sol_WIN}, which demonstrates the good convergence of the proposed WIN-SVEM.
    Further, the components $\left\{ {\bf d}_{{\bf u}_x,i} \in \mathbb{R}^{2018} \right\}_{i=1}^4$ and $\left\{ {\bf d}_{{\bf u}_y,i} \in \mathbb{R}^{2018} \right\}_{i=1}^4$ of the deterministic vectors $\left\{ {\bf d}_{{\rm WIN},i} \right\}_{i=1}^4$ in the $x$ and $y$ directions are depicted in the first and second lines of \figref{fig_e1_DL_6}, respectively,
    and probability density functions (PDFs) of corresponding random variables $\left\{ \lambda_i\left( \theta \right) \right\}_{i=1}^4$ solved by  \eqref{eq:SVEE_D} are seen from the third line of \figref{fig_e1_DL_6}.
    As the retained term $k$ increases, the ranges of random variables $\left\{ \lambda_i\left( \theta \right) \right\}_{i=1}^4$ become smaller and more concentrated around zero, which indicates that the subsequent retained terms contribute less and less to the stochastic solution.

\begin{figure}[!ht]
    	\begin{minipage}{0.50\textwidth}
    		\centering
    		\subfloat[][Comparison of PDFs of the stochastic displacement $u_{A,x}\left( \theta \right)$ in the $x$ direction of the point A obtained by PC-, WIN-SVEMs and MCS, respectively.]{\label{fig_e1_u_x_PDF} \includegraphics[width=1.0\textwidth]{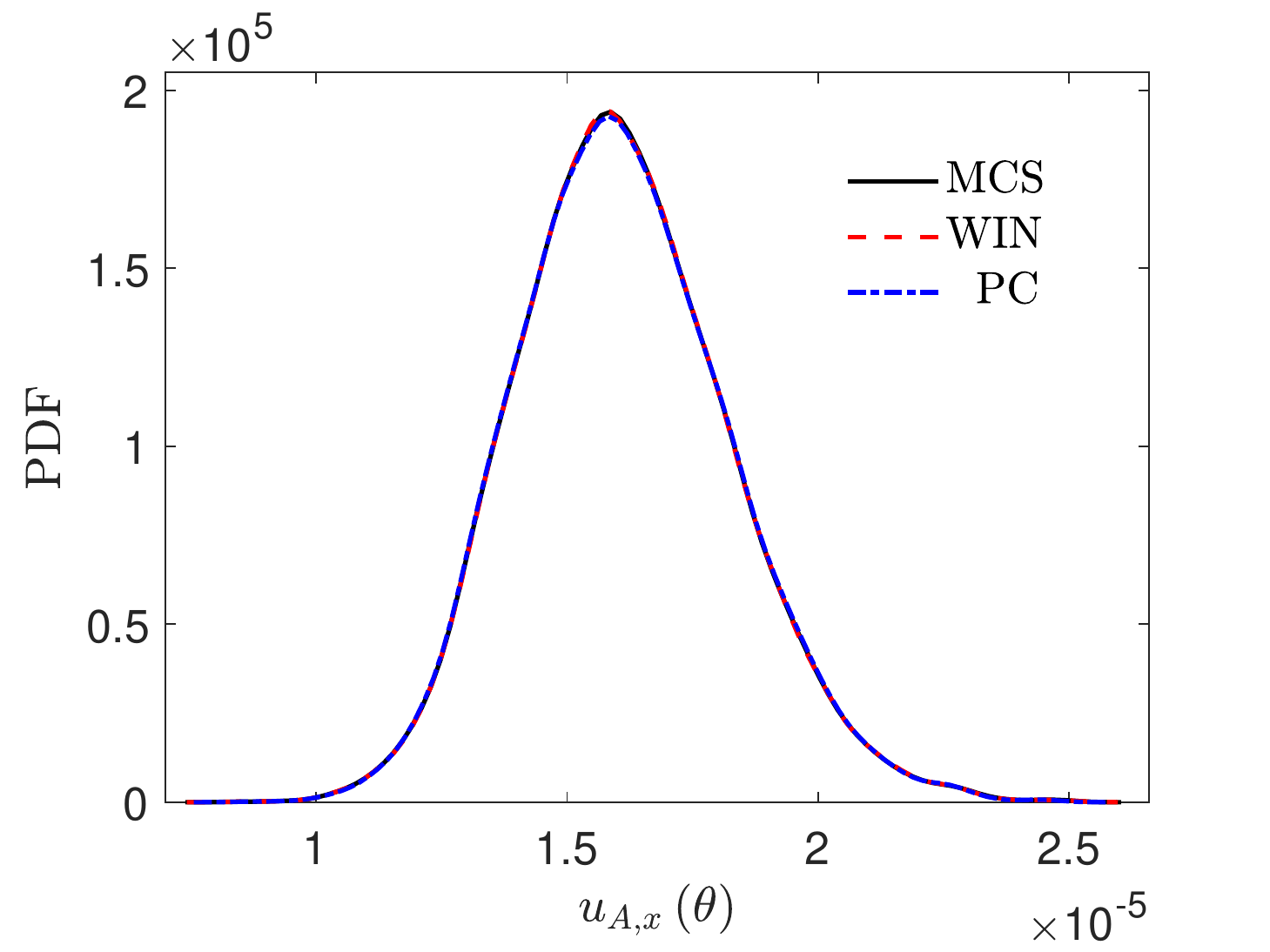}}
    	\end{minipage}
        \hspace{0.12cm}
    	\begin{minipage}{0.50\textwidth}
    		\centering
    		\subfloat[][Absolute errors of PDFs of the stochastic displacement $u_{A,x}\left( \theta \right)$ in the $x$ direction of the point A obtained by PC- and WIN-SVEMs to the MCS  reference PDF.]{\label{fig_e1_u_x_PDF_error} \includegraphics[width=1.0\textwidth]{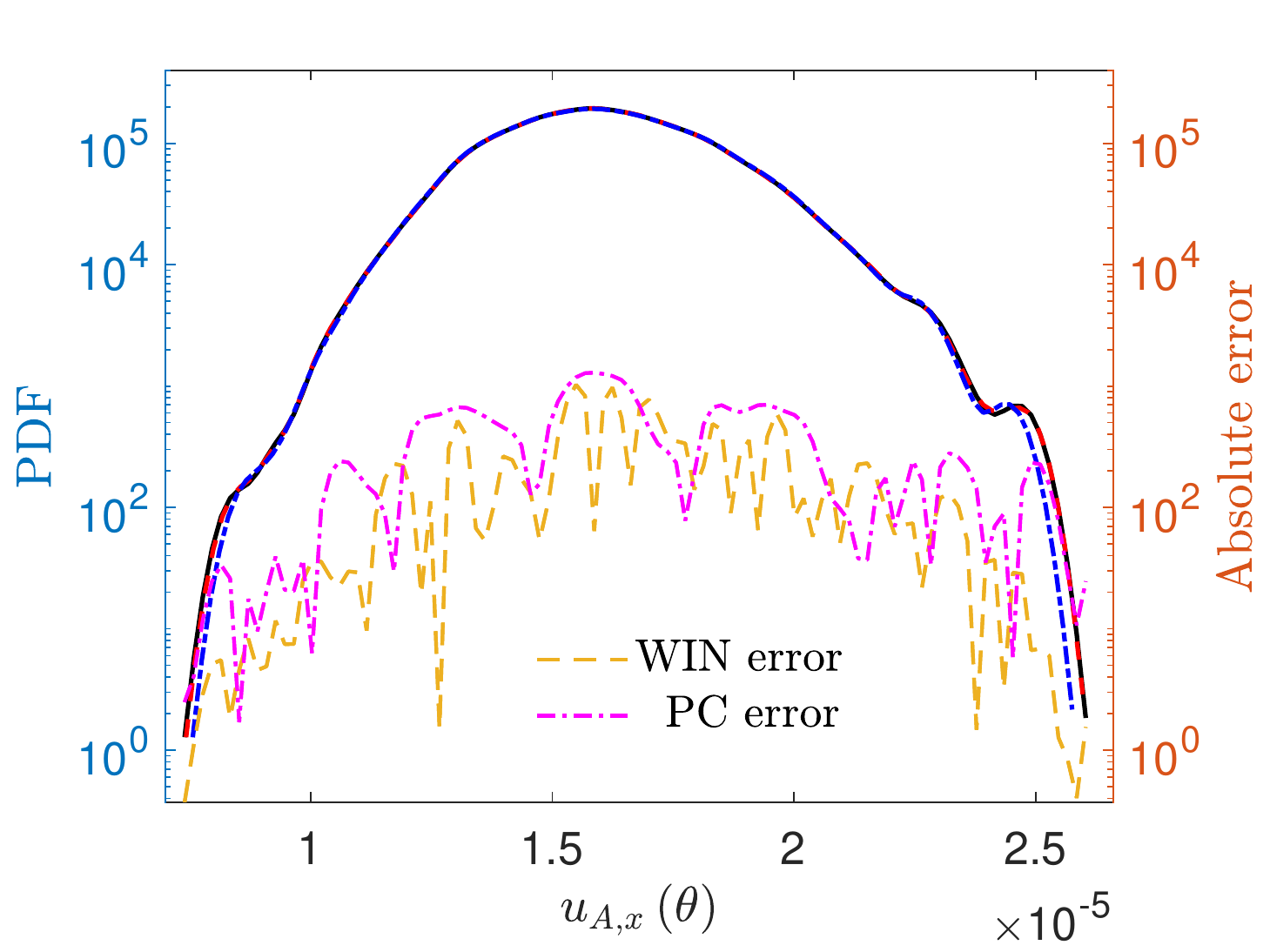}}
    	\end{minipage}
     	\begin{minipage}{0.50\textwidth}
    		\centering
    		\subfloat[][Comparison of PDFs of the stochastic displacement $u_{A,y}\left( \theta \right)$ in the $y$ direction of the point A obtained by PC-, WIN-SVEMs and MCS, respectively.]{\label{fig_e1_u_y_PDF} \includegraphics[width=1.0\textwidth]{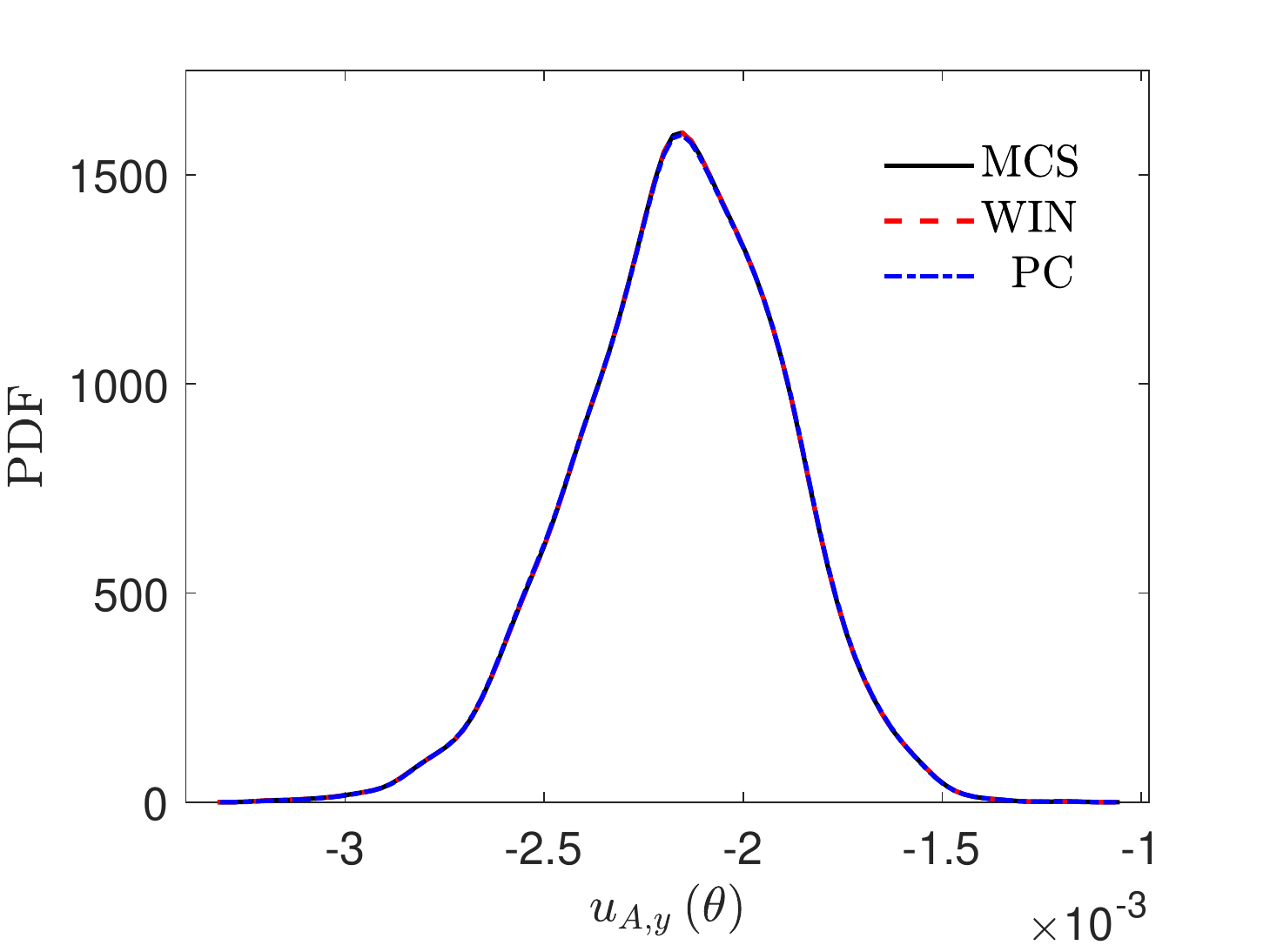}}
    	\end{minipage}
        \hspace{0.12cm}
    	\begin{minipage}{0.50\textwidth}
    		\centering
    		\subfloat[][Absolute errors of PDFs of the stochastic displacement $u_{A,y}\left( \theta \right)$ in the $y$ direction of the point A obtained by PC- and WIN-SVEMs to the MCS  reference PDF.]{\label{fig_e1_u_y_PDF_error} \includegraphics[width=1.0\textwidth]{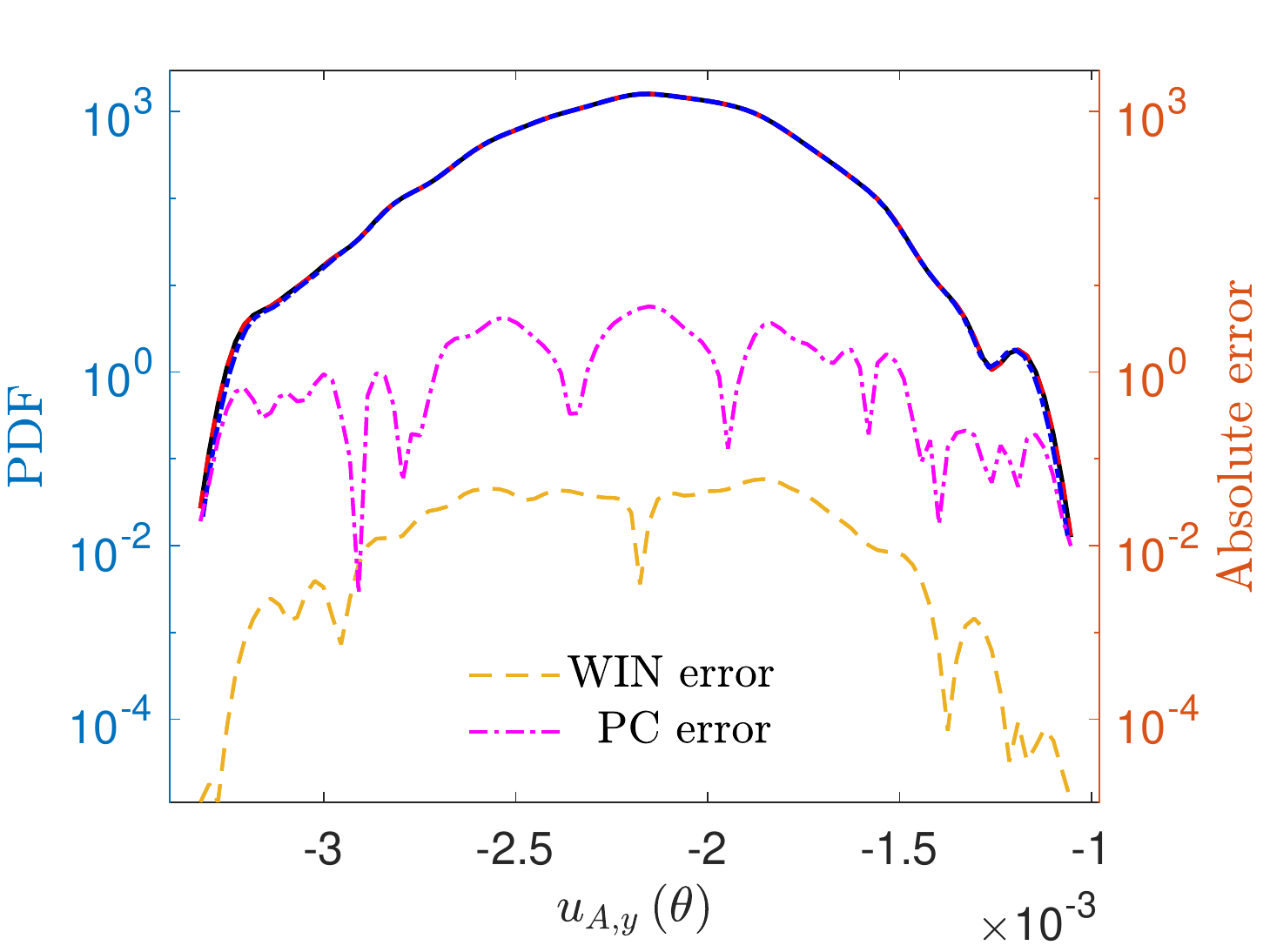}}
    	\end{minipage}
    	\caption{PDFs of the stochastic displacements $u_{A,x}\left( \theta \right)$ and $u_{A,y}\left( \theta \right)$ in the $x$ and $y$ directions of the point A obtained by PC-, WIN-SVEMs and MCS and their absolute errors.} \label{fig_e1_PDF}
    \end{figure}

    To show the computational accuracy of the proposed PC-SVEM and WIN-SVEM, we compare PDFs of the stochastic displacements $u_{A,x}\left( \theta \right)$ and $u_{A,y}\left( \theta \right)$ in the $x$ and $y$ directions of the point A (i.e. the blue point where the force $f\left( \theta \right)$ acts as shown in \figref{fig_e1_model}) obtained by PC-, WIN- and MCS-based SVEMs and their absolute errors in \figref{fig_e1_PDF}.
    It is seen from \figref{fig_e1_u_x_PDF} and \figref{fig_e1_u_y_PDF} that the PDFs of both stochastic displacements $u_{A,x}\left( \theta \right)$ and $u_{A,y}\left( \theta \right)$ obtained by PC- and WIN-SVEMs are in good accordance with those of MCS, which illustrates the good accuracy of the two proposed methods.
    The comparison in logarithmic scales shown in \figref{fig_e1_u_x_PDF_error} and \figref{fig_e1_u_y_PDF_error} demonstrates that WIN-SVEM can achieve smaller absolute errors than PC-SVEM, especially for the stochastic displacement $u_{A,y}\left( \theta \right)$ in the $y$ direction.
    Further, WIN-SVEM can capture tails of the PDFs more accurately than PC-SVEM, which is very useful for many uncertainty quantification problems, such as the simulation of physical phenomena with long tailed probability distributions and the estimation of small failure probability in structural reliability analysis.
    Therefore, WIN-SVEM is recommended for problems with such requirements.

    Further, let us focus on the computational efficiency of the proposed methods.
    Computational times (unit: second) of the numerical execution of PC-SVEM, WIN-SVEM and MCS are listed in the second to fourth columns of Table \ref{tab_TimeCost}, where the solving and recalculating times of WIN-SVEM are the computational times of step \ref{Alg:step_02} to step \ref{Alg:step_13} and the recalculation step \ref{Alg:step_14} of Algorithm \ref{alg_WIN_SVEM}, respectively.
    The cost of the recalculation process of WIN-SVEM is almost negligible since only deterministic equations with size 4 are solved for different sample realizations. 
    It is found that both PC-SVEM and WIN-SVEM are much cheaper than MCS.
    However, since a very large deterministic equation needs to be solved in PC-SVEM, it is more computationally intensive than WIN-SVEM.

	\begin{table*}[ht]
		\caption{Computational costs of the stochastic dimensions 7 and 35.} \label{tab_TimeCost}
		\centering
		\begin{tabular}{c|ccc|cc}
			\toprule
			Method & WIN & PC & MCS & WIN & MCS \\
            Stochastic dimension & \multicolumn{3}{c|}{7} & \multicolumn{2}{c}{35} \\
			\hline
			Solving time & 0.51 & & & 12.33 \\
			Recalculating time & 0.04 & & & 0.35 \\
			\hline
			Total time (second) & 0.55 & 26.09 & 159.13 & 12.68 & 221.52 \\
			\bottomrule
		\end{tabular}
	\end{table*}

\subsubsection{High-dimensional case} \label{sec:ex_1_high_case}

\begin{figure}[t]
		\centering
		\includegraphics[width=0.5\linewidth]{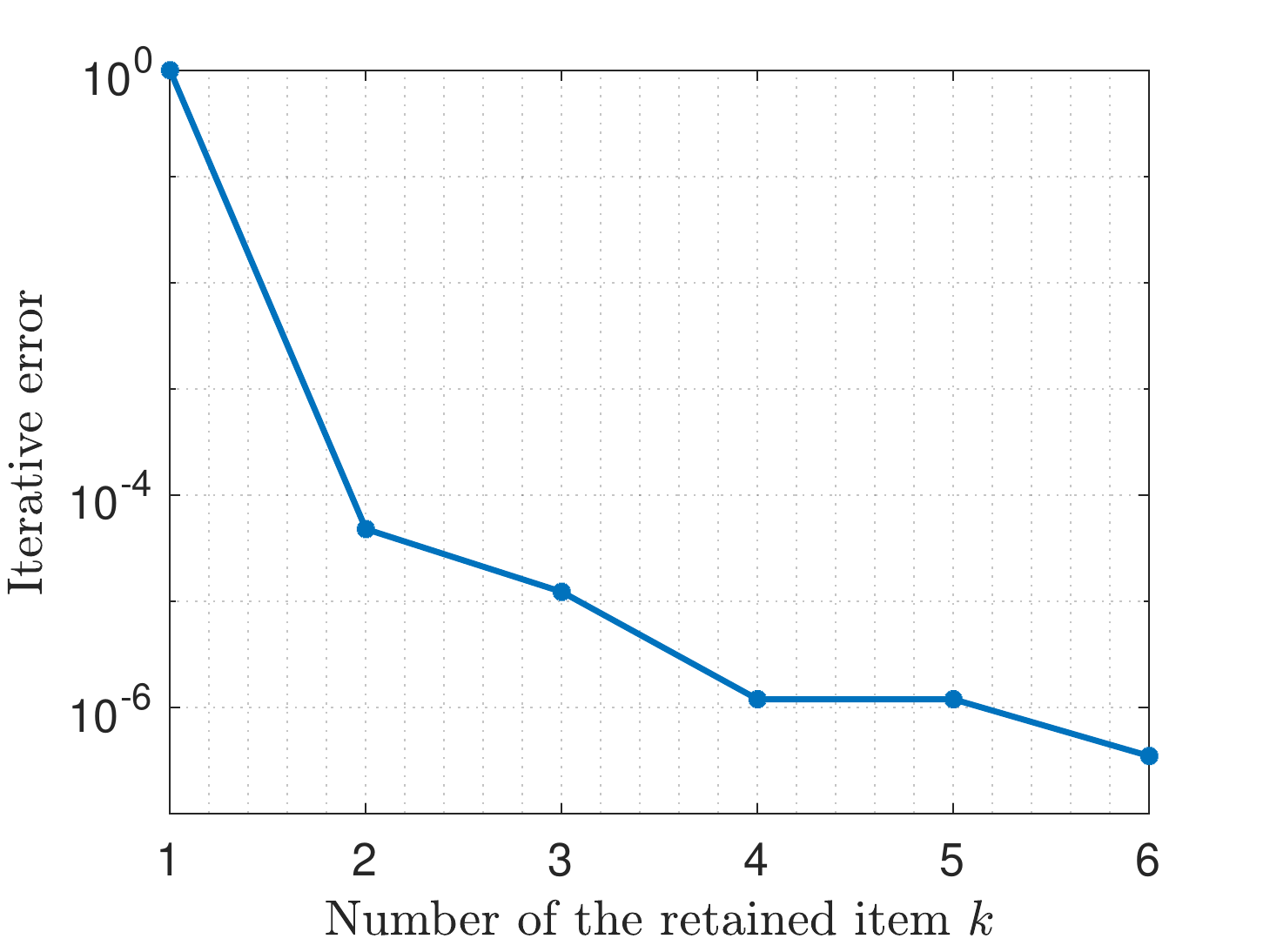}
		\caption{Iterative errors of different retained terms in the high-dimensional case.}\label{fig_e1_It_err_HD}
	\end{figure}

\begin{figure}[!ht]
    	\begin{minipage}{0.5\textwidth}
    		\centering
    		\subfloat[][Comparison of PDFs of the stochastic displacement $u_{A,x}\left( \theta \right)$ in the $x$ direction of the point A obtained by WIN-SVEM and MCS, respectively.]{\label{fig_e1_u_x_PDF_HD} \includegraphics[width=1.0\textwidth]{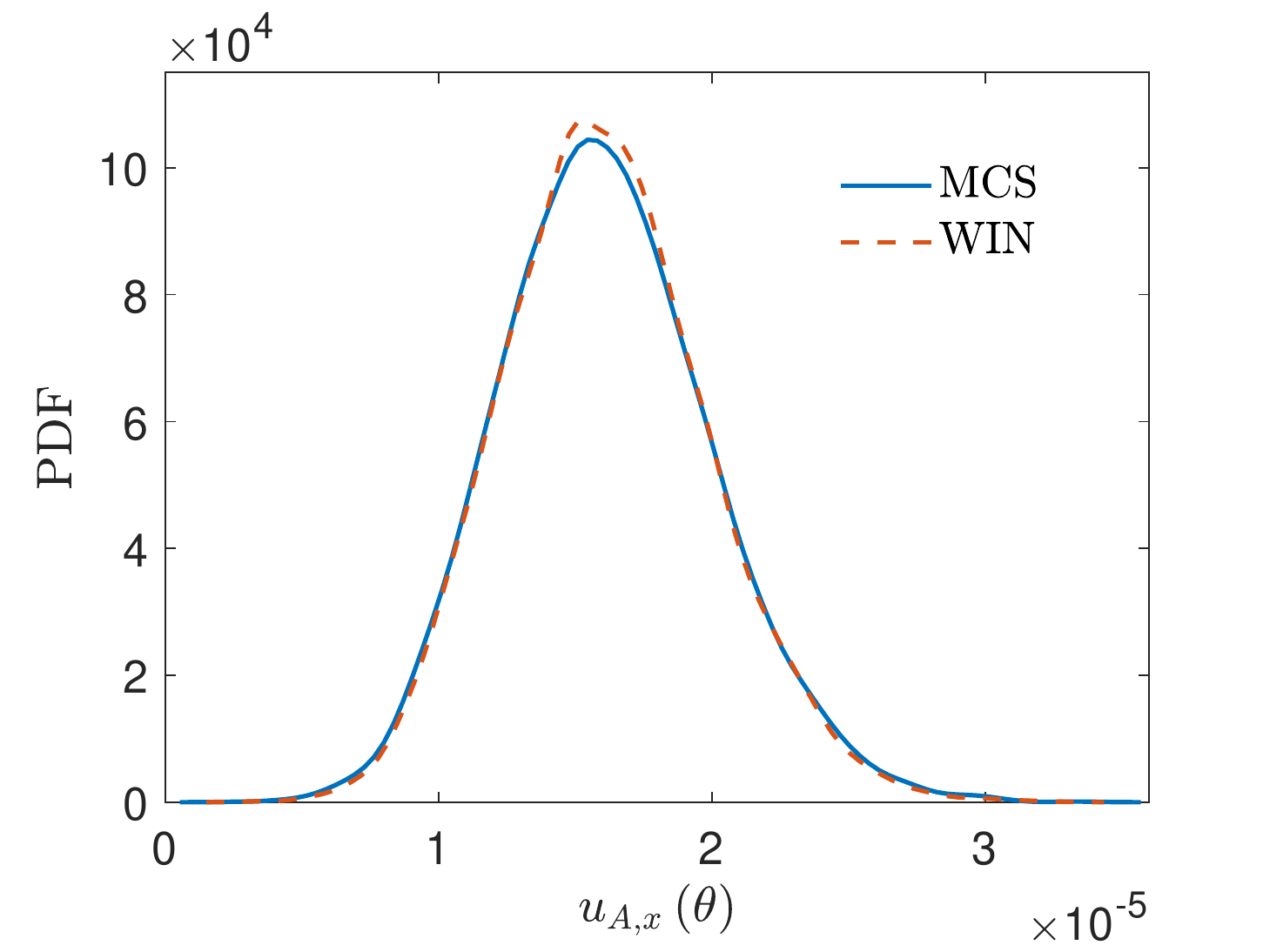}}
    	\end{minipage}
        \hspace{0.12cm}
    	\begin{minipage}{0.5\textwidth}
    		\centering
    		\subfloat[][Comparison of PDFs of the stochastic displacement $u_{A,y}\left( \theta \right)$ in the $y$ direction of the point A obtained by WIN-SVEM and MCS, respectively.]{\label{fig_e1_u_y_PDF_HD} \includegraphics[width=1.0\textwidth]{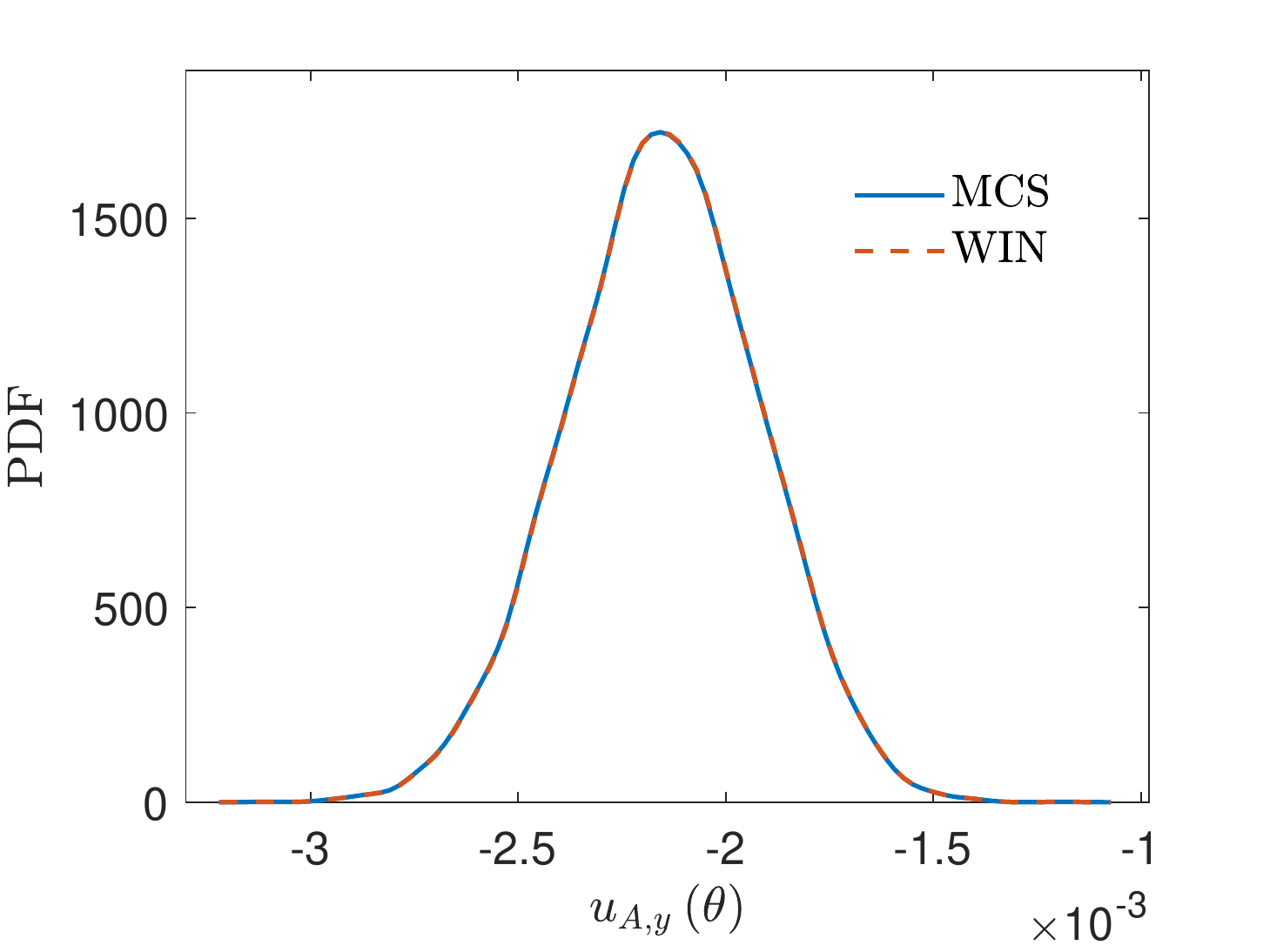}}
    	\end{minipage}
    	\caption{PDFs of the stochastic displacements $u_{A,x}\left( \theta \right)$ (left) and $u_{A,y}\left( \theta \right)$ (right) in the $x$ and $y$ directions of the point A obtained by WIN-SVEM and MCS, respectively.} \label{fig_e1_PDF_HD}
    \end{figure}

    In this section, we solve the high-dimensional case with a total of 35 stochastic dimensions ($\left\{ \xi_i\left( \theta \right) \right\}_{i=1}^{34}$ and $\xi_f\left( \theta \right)$).
    Only WIN-SVEM is adopted in this section.
    For PC-SVEM, even if we only use the second order PC basis, the size of the derived deterministic equation (\ref{eq:Sys_PC}) is about $2.69 \times 10^6$, which suffers from the curse of dimensionality.
    By using Algorithm \ref{alg_WIN_SVEM}, the iterative errors ${\epsilon _{{\bf u},k}}$ of different retained terms for the high-dimensional case are shown in \figref{fig_e1_It_err_HD}.
    The proposed WIN-SVEM still has good convergence for high-dimensional stochastic problems.
    $k = 6$ terms are retained in this case, which is slightly increased compared to the low-dimensional case. 
    Computational times for this case are seen from the fifth and sixth columns of Table \ref{tab_TimeCost}.
    WIN-SVEM has very low cost even for high-dimensional cases and is much cheaper than MCS.
    Compared to the low-dimensional case, the computational cost of the high-dimensional case does not increase dramatically as the stochastic dimension increases.
    In these senses, the proposed WIN-SVEM avoids the curse of dimensionality successfully.
    Further, PDFs of the stochastic displacements $u_{A,x}\left( \theta \right)$ and $u_{A,y}\left( \theta \right)$ in the $x$ and $y$ directions of the point A obtained by WIN-SVEM and MCS are compared in \figref{fig_e1_PDF_HD}.
    For both two stochastic displacements, the computational accuracy of WIN-SVEM is still comparable to MCS.

\vspace{-0.11cm}

\subsection{Example 2: SVEM analysis of a 3D mechanical part} \label{sec:ex_2}

    \begin{figure}[ht] 
    	\begin{minipage}{0.5\textwidth}
    		\centering
    		\subfloat[][Geometric model.]{\label{fig_e2_model} \includegraphics[width=0.9\textwidth]{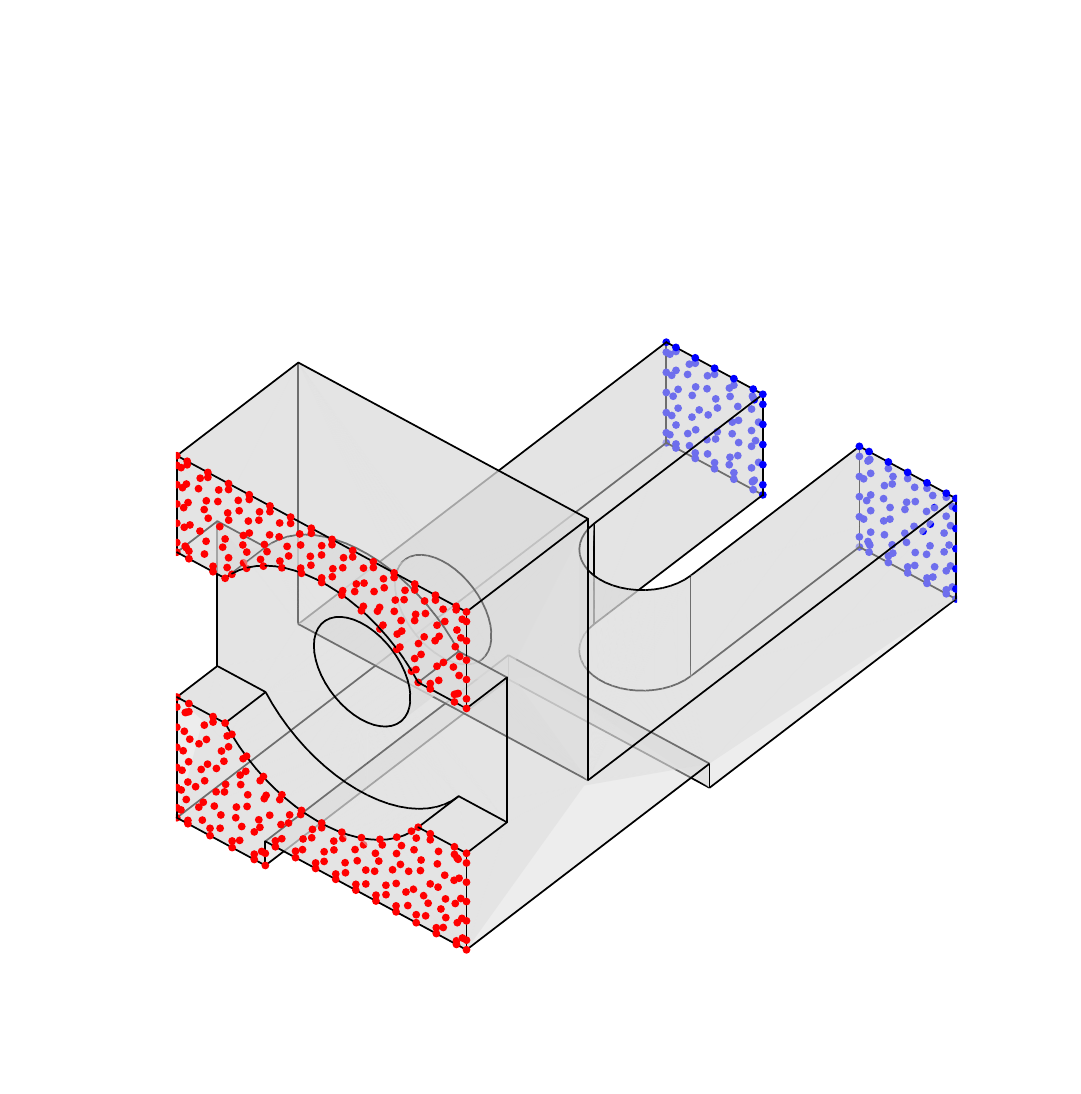}}
    	\end{minipage}
    	\begin{minipage}{0.5\textwidth}
    		\centering
    		\subfloat[][Voronoi mesh.]{\label{fig_e2_mesh} \includegraphics[width=0.9\textwidth]{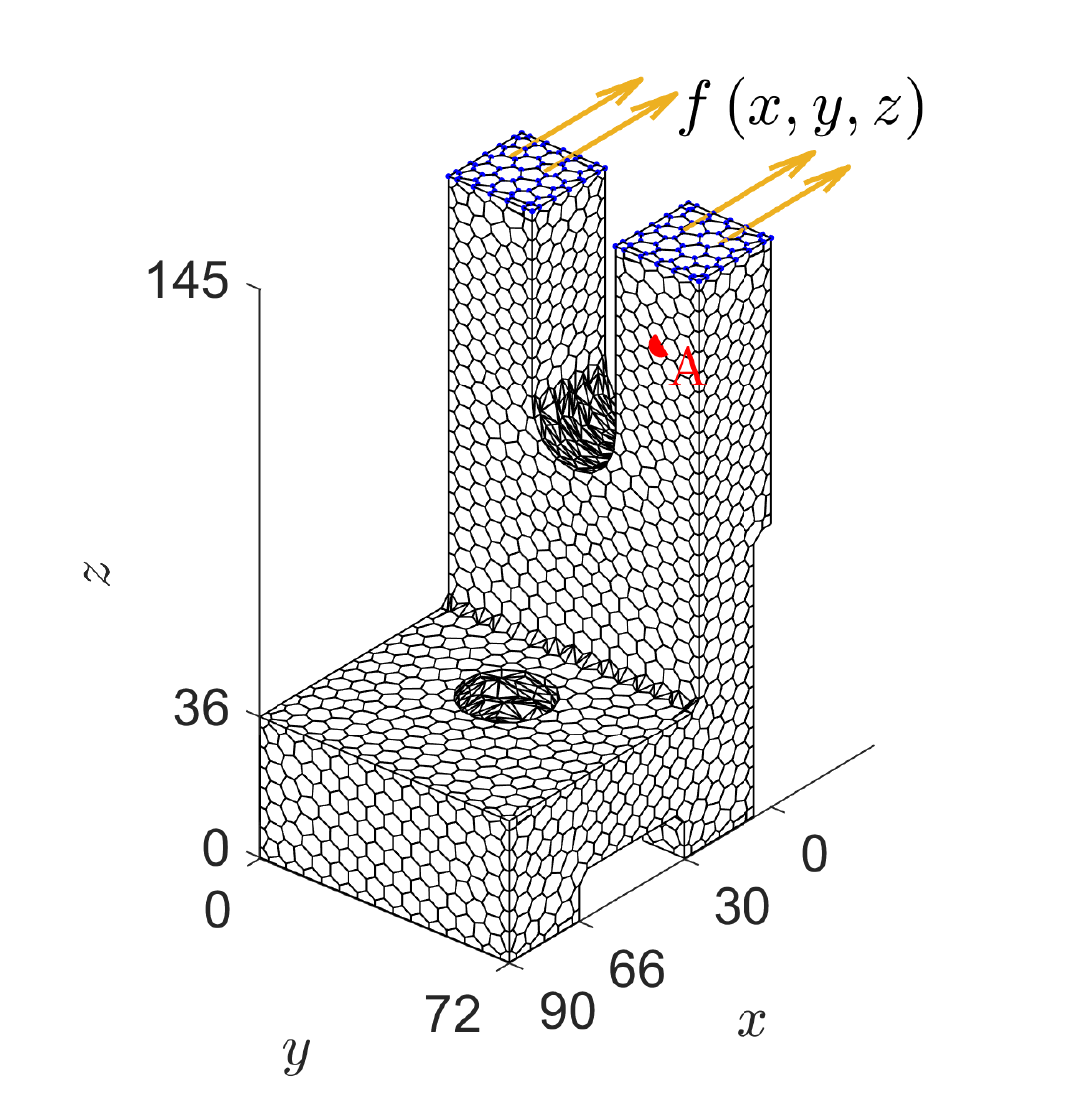}}
    	\end{minipage}
    	\caption{Geometry of the 3D mechanical part (left) and its Voronoi mesh (right).} \label{fig_e2_model_mesh}
    \end{figure}

    In this case, we consider the SVEM analysis of a 3D mechanical part as shown in \figref{fig_e2_model}, where Dirichlet boundary conditions ${\bf u}_x\left( \theta \right) = {\bf u}_y\left( \theta \right) = {\bf u}_z\left( \theta \right) = {\bf 0}$ are imposed on the red surface and an external force $f\left( x,y,z \right) = -500~{\rm N/mm^2}$ is applied along the $x$ direction on the blue surface.
    The model is discretized by use of the Voronoi mesh depicted in \figref{fig_e2_mesh}, including a total of 28232 vertices, 4389 elements and 84696 DoFs.
    In this example, the stochastic material matrix ${\bf G}\left( x,y,z, \theta \right)$ is given by

	\vspace{-0.5cm}
	\begin{equation} \label{eq:G_Ex2}
        {\bf G}\left( x,y,z, \theta \right) = \frac{E\left( x,y,z, \theta \right)}{\left( 1 + \nu \right)\left( 1 - 2\nu \right)}
        \left[ {\begin{array}{*{20}{c}}
        1 - \nu & \nu & \nu & 0 & 0 & 0 \\
        \nu & 1 - \nu & \nu & 0 & 0 & 0 \\
        \nu & \nu & 1 - \nu & 0 & 0 & 0 \\
        0 & 0 & 0 & \frac{1}{2} -\nu & 0 & 0 \\
        0 & 0 & 0 & 0 & \frac{1}{2} -\nu & 0 \\
        0 & 0  & 0 & 0 & 0 & \frac{1}{2} -\nu
        \end{array}} \right] \in \mathbb{R}^{6 \times 6},
	\end{equation}
    where the Poisson ratio $\nu = 0.3$ and the Young's modulus $E\left( x,y,z, \theta \right)$ is a three-dimensional random field with the covariance function
    
    \vspace{-0.5cm}
    \begin{equation} \label{eq:Ex2_Cov}
        {\rm Cov}_{EE}\left(x_1,y_1,z_1;x_2,y_2,z_2\right) = \sigma_E^2 \exp \left( { - \frac{{\left| x_1 - x_2 \right|}}{l_x} - \frac{{\left| y_1 - y_2 \right|}}{l_y} - \frac{{\left| z_1 - z_2 \right|}}{l_z}} \right),
    \end{equation}
    where the standard deviation $\sigma_E = 41.8$~GPa, and the correlation lengths are given by $l_x = \max(x) - \min(x)$, $l_y = \max(y) - \min(y)$ and $l_z = \max(z) - \min(z)$.
    The random field $E\left( x,y,z, \theta \right)$ has a \eqref{eq:Ex1_E}-like series expansion

    \vspace{-0.5cm}
    \begin{equation}\label{eq:Ex2_E}
        E\left( x,y,z, \theta \right) = E_0\left( x,y,z \right) + \sum\limits_{i=1}^m \xi_i\left( \theta \right) \sqrt{\kappa_i} E_i\left( x,y,z \right),
    \end{equation}
    where the function $E_0\left( x,y,z \right) = 208$~GPa, $\left\{ \xi_i\left( \theta \right) \right\}_{i=1}^m$ are mutually independent uniform random variables on $\left[ 0,1 \right]$.
    It is noted that the mean value of the random field $E\left( x,y,z, \theta \right)$ is $E_0\left( x,y,z \right) + 0.5 \sum_{i=1}^m \sqrt{\kappa_i} E_i\left( x,y,z \right)$ instead of $E_0\left( x,y,z \right)$.
    Similarly, $\left\{ \kappa_i, E_i\left( x,y \right) \right\}_{i=1}^m$ are eigenvalues and eigenvectors of the covariance function ${\rm Cov}_{EE}\left(x_1,y_1,z_1; x_2,y_2,z_2\right)$ and can be solved by the \eqref{eq:Ex1_Cov_sol}-like integral equation.
    The truncated number is set as $m = 13$ in this case to achieve the truncated error ${\kappa_{13}} \mathord{\left/ {\vphantom{ \kappa_{13} \sum_{i=1}^{13} \kappa_i }} \right. \kern-\nulldelimiterspace} \sum_{i=1}^{13} \kappa_i < 1 \times 10^{-3}$.

 	\begin{figure}[ht]
		\centering
		\includegraphics[width=0.5\linewidth]{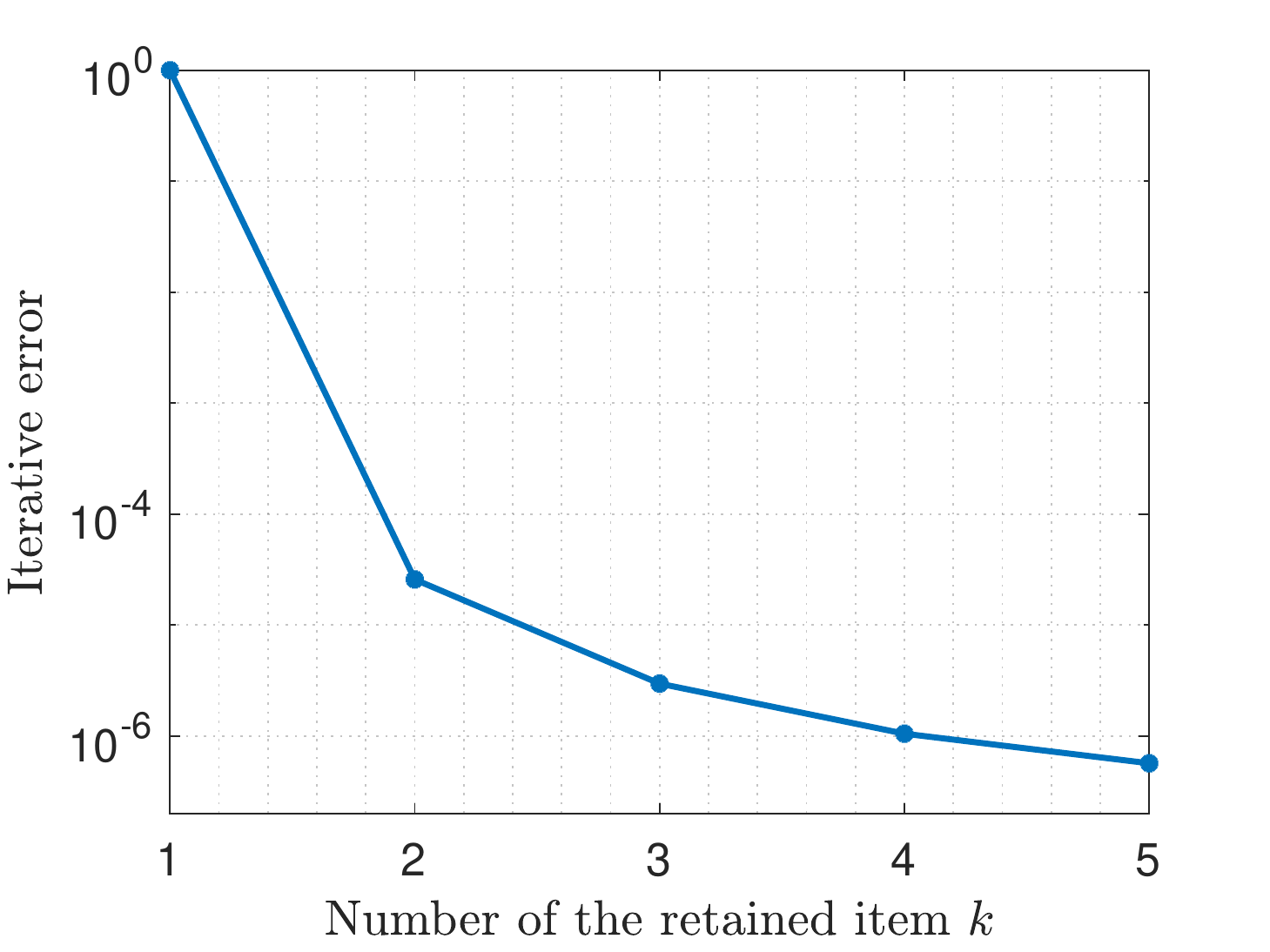}
		\caption{Iterative errors of different retained terms.} \label{fig_e2_It_err}
	\end{figure} 

    \begin{figure}[!ht]
    	\begin{minipage}{0.50\textwidth}
    		\centering
    		\subfloat[][Comparison of PDFs of the stochastic displacement $u_{A,x}\left( \theta \right)$ in the $x$ direction of the point A obtained by WIN-SVEM and MCS, respectively.]{\label{fig_e2_u_x_PDF} \includegraphics[width=1.0\textwidth]{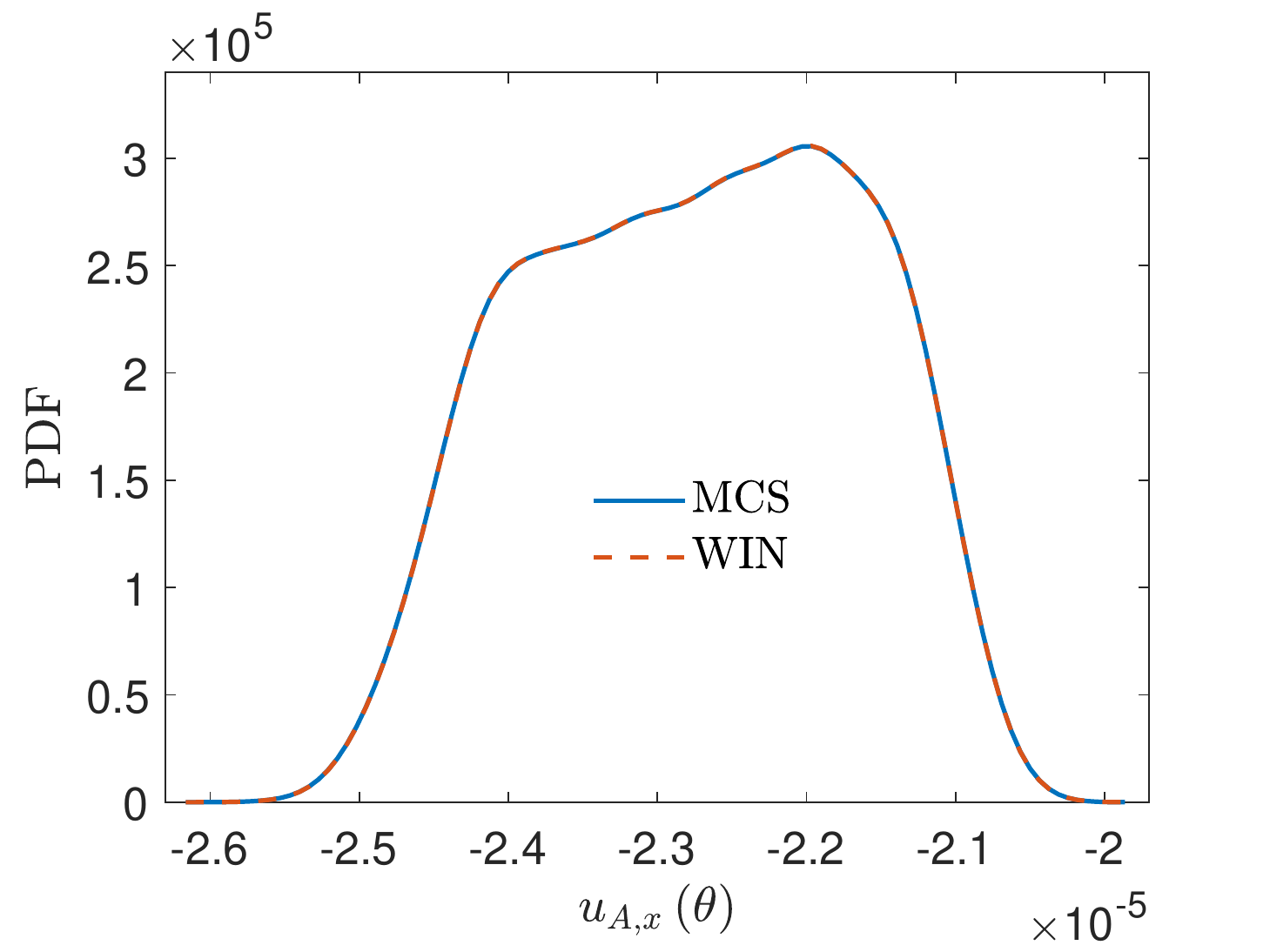}}
    	\end{minipage}
        \hspace{0.12cm}
        \begin{minipage}{0.50\textwidth}
            \centering
            \subfloat[][Comparison of PDFs of the stochastic displacement $u_{A,y}\left( \theta \right)$ in the $y$ direction of the point A obtained by WIN-SVEM and MCS, respectively.]{\label{fig_e2_u_y_PDF} \includegraphics[width=1.0\textwidth]{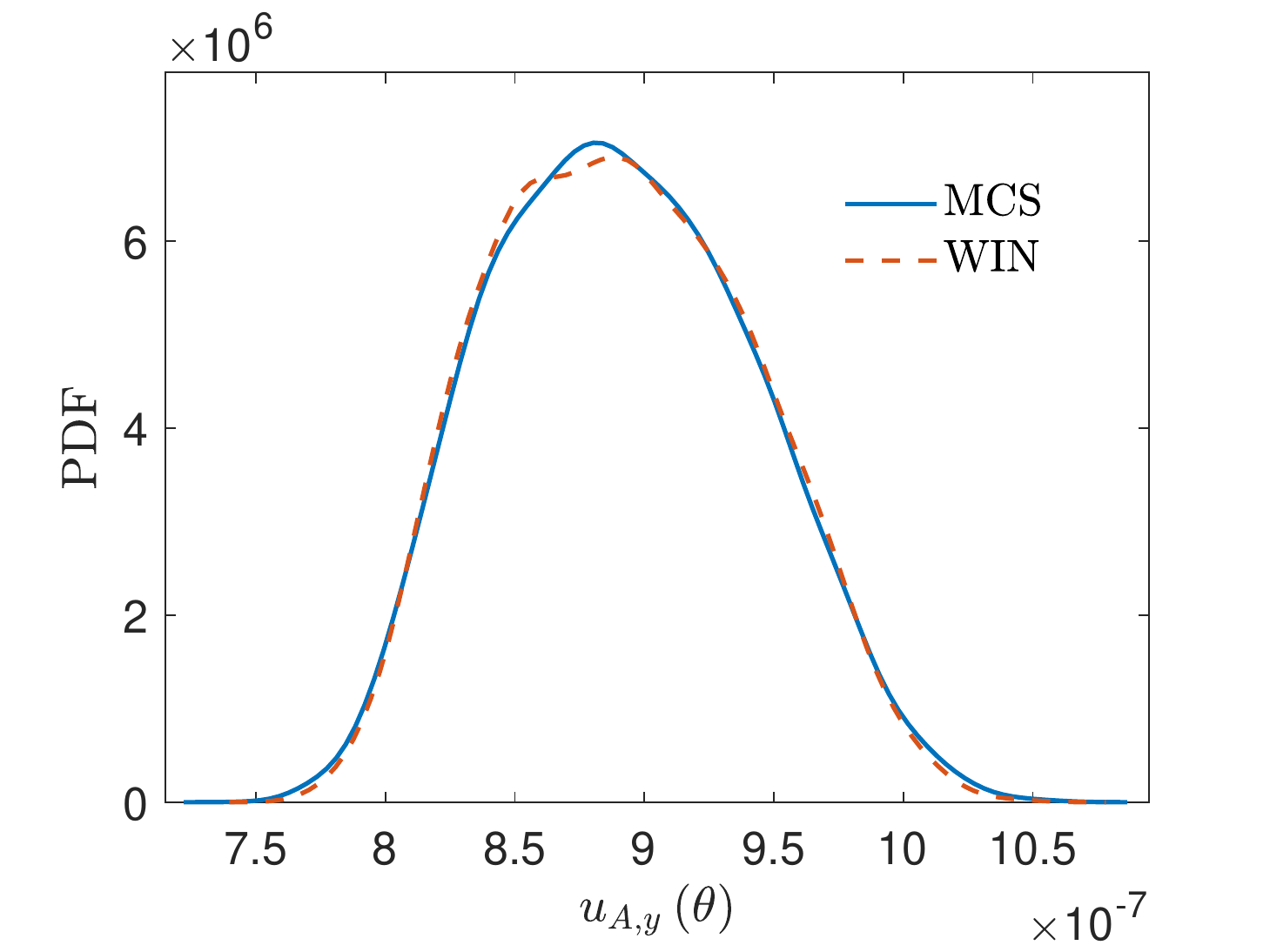}}
    	\end{minipage}
    	\begin{minipage}{1.0\textwidth}
            \centering
            \subfloat[][Comparison of PDFs of the stochastic displacement $u_{A,z}\left( \theta \right)$ in the $z$ direction of the point A obtained by WIN-SVEM and MCS, respectively.]{\label{fig_e2_u_z_PDF} \includegraphics[width=0.5\textwidth]{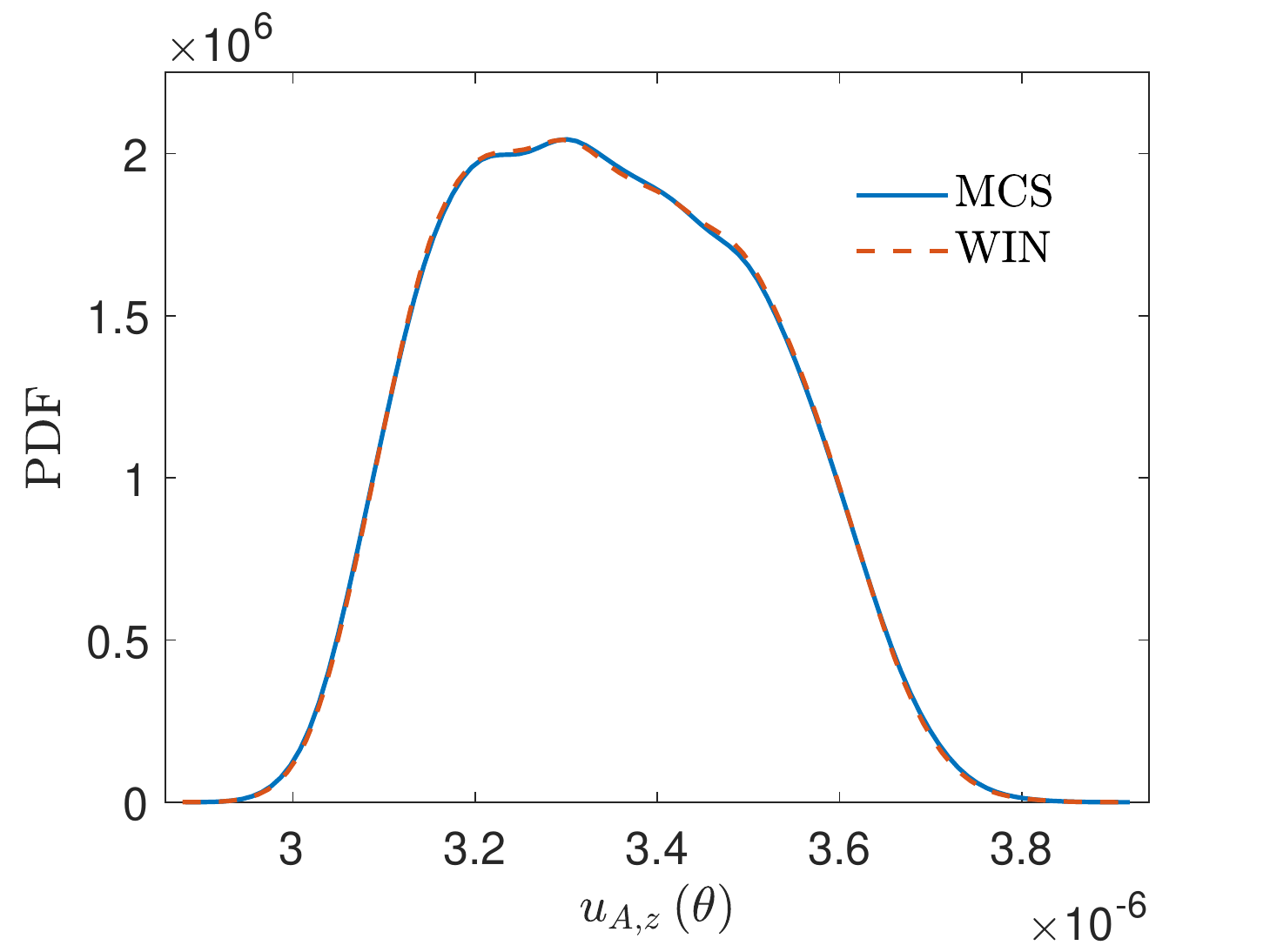}}
    	\end{minipage}
    	\caption{PDFs of the stochastic displacements $u_{A,x}\left( \theta \right)$ (top left), $u_{A,y}\left( \theta \right)$ (top right) and $u_{A,z}\left( \theta \right)$ (bottom) in the $x$, $y$ and $z$ directions of the point A obtained by WIN-SVEM and MCS, respectively.} \label{fig_e2_u_PDF}
    \end{figure}

    Similar to the example in the previous section, PC-SVEM suffers from the curse of dimensionality in this example since the size of the derived deterministic equation (\ref{eq:Sys_PC}) is about $8.89 \times 10^6$, even with only the second order PC basis.
    Thus, only WIN-SVEM is used to solve this problem.
    Iterative errors of different retained terms are seen from \figref{fig_e2_It_err}.
    Five terms are retained to meet the specified convergence error, which verifies the good convergence of WIN-SVEM for 3D stochastic problems.
    Regarding the computational accuracy, PDFs of the stochastic displacements $u_{A,x}\left( \theta \right)$, $u_{A,y}\left( \theta \right)$ and $u_{A,z}\left( \theta \right)$ in the $x$, $y$ and $z$ directions of the point A (shown in \figref{fig_e2_mesh}) obtained by WIN-SVEM and MCS are compared in \figref{fig_e2_u_PDF}.  
    For the PDFs of stochastic displacements $u_{A,x}\left( \theta \right)$ and $u_{A,z}\left( \theta \right)$, WIN-SVEM is in very good agreement with MCS.
    The PDF of the stochastic displacement $u_{A,y}\left( \theta \right)$ is slightly less accurate than those of $u_{A,x}\left( \theta \right)$ and $u_{A,z}\left( \theta \right)$, but acceptable accuracy is still achieved.
    We can simply retain more terms in the stochastic solution approximation to improve the accuracy if higher accuracy is required in practice.
    Further, the computational times of WIN-SVEM and MCS in this example are 213.45s (including the solving time 212.89s and the recalculating time 0.56s) and $2.74 \times 10^4$s, respectively.
    The proposed WIN-SVEM is still much cheaper than MCS and a speedup of more than one hundred times is achieved.

 	\begin{figure}[ht] 
		\centering
		\includegraphics[width=0.8\linewidth]{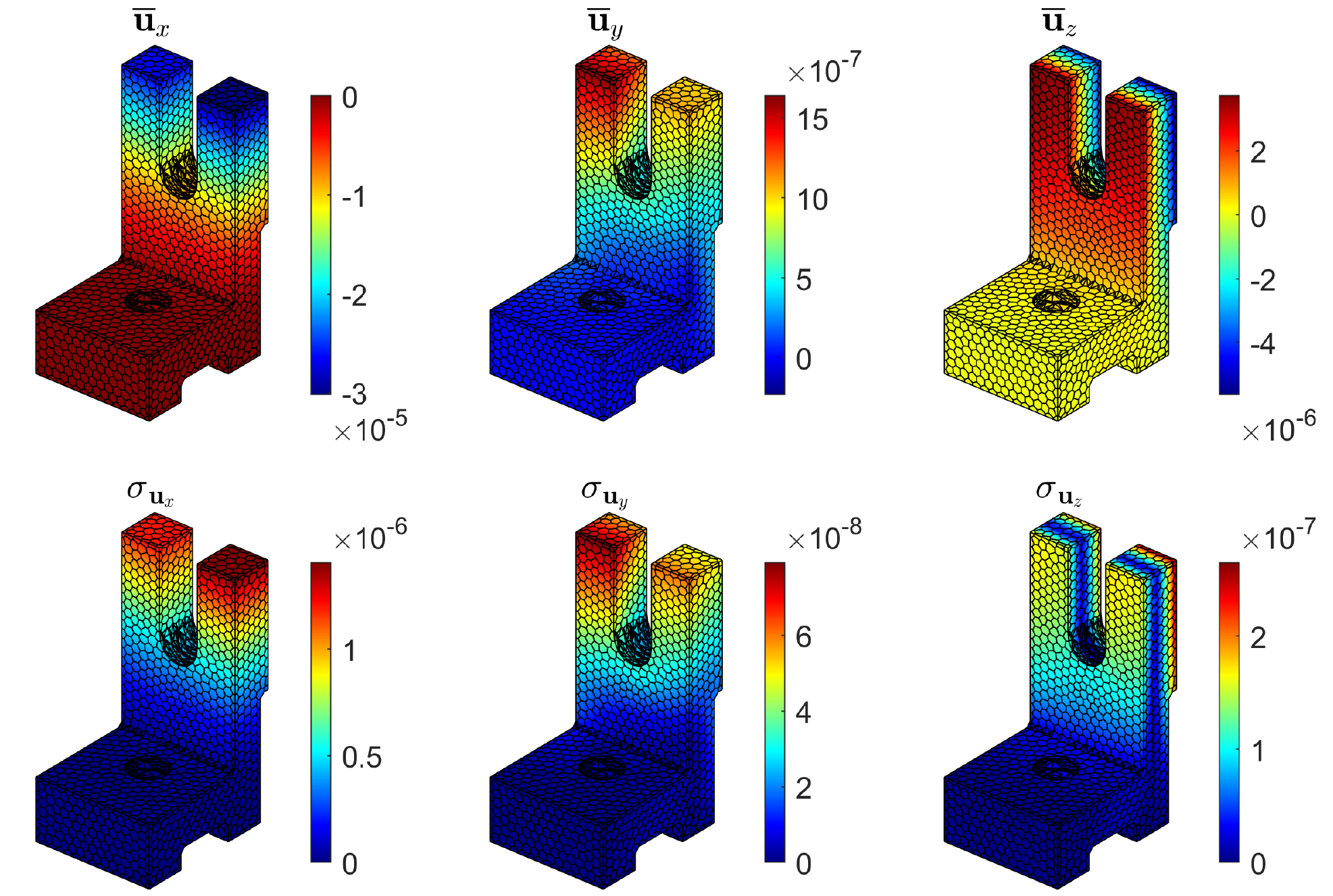}
		\caption{The mean functions $\overline{{\bf u}}_x$, $\overline{{\bf u}}_y$ and $\overline{{\bf u}}_z$ (the first line) in the $x$, $y$ and $z$ directions, and the standard derivation functions ${\bm \sigma}_{{\bf u}_x}$, ${\bm \sigma}_{{\bf u}_y}$ and ${\bm \sigma}_{{\bf u}_z}$ (the second line) in the $x$, $y$ and $z$ directions.} \label{fig_e2_mean_std}
	\end{figure}

    We also highlight that statistical properties of the stochastic solution are easily computed based on WIN-SVEM.
    Here we focus on the first and second order global statistical moments, that is, the mean value and the standard deviation dependent on spatial positions.
    According to \eqref{eq:sol_WIN}, the mean value vector $\overline{{\bf u}}$ is computed as

    \vspace{-0.5cm}
    \begin{equation} \label{eq:mean_fun}
        \overline{{\bf u}} = \sum\limits_{i=1}^k {\widehat{\mathbb{E}}}\left\{ \lambda_i \left( {\widehat{\bm \theta}} \right) \right\} {\bf d}_{{\rm WIN},i},
    \end{equation}
    which only involves the expectations of sample vectors $\left\{ \lambda_i \left( {\widehat{\bm \theta}} \right) \in \mathbb{R}^{10^4} \right\}_{i=1}^k$ and has very low computational effort.
    The mean value components $\overline{{\bf u}}_x$, $\overline{{\bf u}}_y$ and $\overline{{\bf u}}_z$ of $\overline{{\bf u}}$ in the $x$, $y$ and $z$ directions are seen from the first line of \figref{fig_e2_mean_std}.
    Further, the standard deviation vector ${\bm \sigma}_{\bf u}$ is calculated via

    \vspace{-0.5cm}
    \begin{align} \label{eq:std_fun}
        {\bm \sigma}_{\bf u} &= \sqrt{ \sum\limits_{i,j=1}^k {\widehat{\mathbb{E}}}\left\{ \left[ \lambda_i \left( {\widehat{\bm \theta}} \right) - {\widehat{\mathbb{E}}}\left\{ \lambda_i \left( {\widehat{\bm \theta}} \right) \right\} \right] \odot \left[ \lambda_j \left( {\widehat{\bm \theta}} \right) - {\widehat{\mathbb{E}}}\left\{ \lambda_j \left( {\widehat{\bm \theta}} \right) \right\}  \right] \right\} {\bf d}_{{\rm WIN},i} \odot {\bf d}_{{\rm WIN},j} } \nonumber \\
        &= \sqrt{ \sum\limits_{i,j=1}^k \left[ {\widehat{\mathbb{E}}}\left\{ \lambda_i \left( {\widehat{\bm \theta}} \right) \odot \lambda_j \left( {\widehat{\bm \theta}} \right) \right\} - {\widehat{\mathbb{E}}}\left\{ \lambda_i \left( {\widehat{\bm \theta}} \right) \right\} {\widehat{\mathbb{E}}}\left\{ \lambda_j \left( {\widehat{\bm \theta}} \right) \right\}  \right] {\bf d}_{{\rm WIN},i} \odot {\bf d}_{{\rm WIN},j} },
    \end{align}
    whose components ${\bm \sigma}_{{\bf u}_x}$, ${\bm \sigma}_{{\bf u}_y}$ and ${\bm \sigma}_{{\bf u}_z}$ in the $x$, $y$ and $z$ directions can be found in the second line of \figref{fig_e2_mean_std}.

\section{Conclusions}\label{sec:Con}
    We presented two numerical approaches, PC-SVEM and WIN-SVEM, for solving stochastic systems derived from the stochastic virtual element discretization of 2D and 3D linear elastic stochastic problems.
    The deterministic virtual element method is first extended to SVEM that involves stochastic material properties and stochastic external forces, etc.
    Several key calculations of SVEM can be inherited from the deterministic virtual element method, e.g. the gradient computations of virtual basis functions in \eqref{eq:B_i} and the stabilizing element stiffness matrix in \eqref{eq:k_S_ele} (except for the random coefficient $\gamma_S ^{\left( {\rm e} \right)}\left( \theta \right)$).
    Numerical results demonstrate that both PC-SVEM and WIN-SVEM have comparable accuracy to MCS.
    However, PC-SVEM suffers from the curse of dimensionality and cannot be applied to high-dimensional stochastic problems.
    As a comparison, WIN-SVEM can efficiently solve both low- and high-dimensional stochastic problems without any modification, which has been verified by a numerical example of up to 35 stochastic dimensions.
    Further, although only linear elastic stochastic problems are concerned in this paper, both the proposed PC-SVEM and WIN-SVEM can be applied to more general cases, which will be further investigated in subsequent research.

\section*{Acknowledgments}
	The authors are grateful to the Alexander von Humboldt Foundation and the International Research Training Group 2657 (IRTG 2657) funded by the German Research Foundation (DFG) (Grant number 433082294).


	\nocite{*}
	\bibliography{References}

\begin{thebibliography}{42}
\expandafter\ifx\csname natexlab\endcsname\relax\def\natexlab#1{#1}\fi
\providecommand{\url}[1]{\texttt{#1}}
\providecommand{\href}[2]{#2}
\providecommand{\path}[1]{#1}
\providecommand{\DOIprefix}{doi:}
\providecommand{\ArXivprefix}{arXiv:}
\providecommand{\URLprefix}{URL: }
\providecommand{\Pubmedprefix}{pmid:}
\providecommand{\doi}[1]{\href{http://dx.doi.org/#1}{\path{#1}}}
\providecommand{\Pubmed}[1]{\href{pmid:#1}{\path{#1}}}
\providecommand{\bibinfo}[2]{#2}
\ifx\xfnm\relax \def\xfnm[#1]{\unskip,\space#1}\fi
\bibitem[{Quarteroni and Valli(2008)}]{quarteroni2008numerical}
\bibinfo{author}{A.~Quarteroni}, \bibinfo{author}{A.~Valli},
  \bibinfo{title}{Numerical approximation of partial differential equations},
  \bibinfo{publisher}{Springer Science \& Business Media},
  \bibinfo{year}{2008}.
\bibitem[{Beir{\~a}o~da Veiga et~al.(2013)Beir{\~a}o~da Veiga, Brezzi,
  Cangiani, Manzini, Marini, and Russo}]{beirao2013basic}
\bibinfo{author}{L.~Beir{\~a}o~da Veiga}, \bibinfo{author}{F.~Brezzi},
  \bibinfo{author}{A.~Cangiani}, \bibinfo{author}{G.~Manzini},
  \bibinfo{author}{L.~D. Marini}, \bibinfo{author}{A.~Russo},
\newblock \bibinfo{title}{Basic principles of virtual element methods},
\newblock \bibinfo{journal}{Mathematical Models and Methods in Applied
  Sciences} \bibinfo{volume}{23} (\bibinfo{year}{2013})
  \bibinfo{pages}{199--214}.
\bibitem[{Beir{\~a}o~da Veiga et~al.(2014)Beir{\~a}o~da Veiga, Brezzi, Marini,
  and Russo}]{beirao2014hitchhiker}
\bibinfo{author}{L.~Beir{\~a}o~da Veiga}, \bibinfo{author}{F.~Brezzi},
  \bibinfo{author}{L.~D. Marini}, \bibinfo{author}{A.~Russo},
\newblock \bibinfo{title}{The hitchhiker's guide to the virtual element
  method},
\newblock \bibinfo{journal}{Mathematical Models and Methods in Applied
  Sciences} \bibinfo{volume}{24} (\bibinfo{year}{2014})
  \bibinfo{pages}{1541--1573}.
\bibitem[{Mengolini et~al.(2019)Mengolini, Benedetto, and
  Arag{\'o}n}]{mengolini2019engineering}
\bibinfo{author}{M.~Mengolini}, \bibinfo{author}{M.~F. Benedetto},
  \bibinfo{author}{A.~M. Arag{\'o}n},
\newblock \bibinfo{title}{An engineering perspective to the virtual element
  method and its interplay with the standard finite element method},
\newblock \bibinfo{journal}{Computer Methods in Applied Mechanics and
  Engineering} \bibinfo{volume}{350} (\bibinfo{year}{2019})
  \bibinfo{pages}{995--1023}.
\bibitem[{Antonietti et~al.(2022)Antonietti, da~Veiga, and
  Manzini}]{antonietti2022virtual}
\bibinfo{author}{P.~F. Antonietti}, \bibinfo{author}{L.~B. da~Veiga},
  \bibinfo{author}{G.~Manzini}, \bibinfo{title}{The virtual element method and
  its applications}, \bibinfo{publisher}{Springer Science \& Business Media},
  \bibinfo{year}{2022}.
\bibitem[{Da~Veiga et~al.(2013)Da~Veiga, Brezzi, and Marini}]{da2013virtual}
\bibinfo{author}{L.~B. Da~Veiga}, \bibinfo{author}{F.~Brezzi},
  \bibinfo{author}{L.~D. Marini},
\newblock \bibinfo{title}{Virtual elements for linear elasticity problems},
\newblock \bibinfo{journal}{SIAM Journal on Numerical Analysis}
  \bibinfo{volume}{51} (\bibinfo{year}{2013}) \bibinfo{pages}{794--812}.
\bibitem[{Gain et~al.(2014)Gain, Talischi, and Paulino}]{gain2014virtual}
\bibinfo{author}{A.~L. Gain}, \bibinfo{author}{C.~Talischi},
  \bibinfo{author}{G.~H. Paulino},
\newblock \bibinfo{title}{On the virtual element method for three-dimensional
  linear elasticity problems on arbitrary polyhedral meshes},
\newblock \bibinfo{journal}{Computer Methods in Applied Mechanics and
  Engineering} \bibinfo{volume}{282} (\bibinfo{year}{2014})
  \bibinfo{pages}{132--160}.
\bibitem[{Artioli et~al.(2017)Artioli, Beir{\~a}o~da Veiga, Lovadina, and
  Sacco}]{artioli2017arbitrary}
\bibinfo{author}{E.~Artioli}, \bibinfo{author}{L.~Beir{\~a}o~da Veiga},
  \bibinfo{author}{C.~Lovadina}, \bibinfo{author}{E.~Sacco},
\newblock \bibinfo{title}{Arbitrary order {2D} virtual elements for polygonal
  meshes: part {I}, elastic problem},
\newblock \bibinfo{journal}{Computational Mechanics} \bibinfo{volume}{60}
  (\bibinfo{year}{2017}) \bibinfo{pages}{355--377}.
\bibitem[{Chi et~al.(2017)Chi, Da~Veiga, and Paulino}]{chi2017some}
\bibinfo{author}{H.~Chi}, \bibinfo{author}{L.~B. Da~Veiga},
  \bibinfo{author}{G.~Paulino},
\newblock \bibinfo{title}{Some basic formulations of the virtual element method
  ({VEM}) for finite deformations},
\newblock \bibinfo{journal}{Computer Methods in Applied Mechanics and
  Engineering} \bibinfo{volume}{318} (\bibinfo{year}{2017})
  \bibinfo{pages}{148--192}.
\bibitem[{Wriggers et~al.(2021)Wriggers, De~Bellis, and
  Hudobivnik}]{wriggers2021taylor}
\bibinfo{author}{P.~Wriggers}, \bibinfo{author}{M.~De~Bellis},
  \bibinfo{author}{B.~Hudobivnik},
\newblock \bibinfo{title}{A taylor--hood type virtual element formulations for
  large incompressible strains},
\newblock \bibinfo{journal}{Computer Methods in Applied Mechanics and
  Engineering} \bibinfo{volume}{385} (\bibinfo{year}{2021})
  \bibinfo{pages}{114021}.
\bibitem[{Wriggers et~al.(2016)Wriggers, Rust, and Reddy}]{wriggers2016virtual}
\bibinfo{author}{P.~Wriggers}, \bibinfo{author}{W.~T. Rust},
  \bibinfo{author}{B.~Reddy},
\newblock \bibinfo{title}{A virtual element method for contact},
\newblock \bibinfo{journal}{Computational Mechanics} \bibinfo{volume}{58}
  (\bibinfo{year}{2016}) \bibinfo{pages}{1039--1050}.
\bibitem[{Wriggers and Rust(2019)}]{wriggers2019virtual}
\bibinfo{author}{P.~Wriggers}, \bibinfo{author}{W.~T. Rust},
\newblock \bibinfo{title}{A virtual element method for frictional contact
  including large deformations},
\newblock \bibinfo{journal}{Engineering Computations} \bibinfo{volume}{36}
  (\bibinfo{year}{2019}) \bibinfo{pages}{2133--2161}.
\bibitem[{Aldakheel et~al.(2020)Aldakheel, Hudobivnik, Artioli, da~Veiga, and
  Wriggers}]{aldakheel2020curvilinear}
\bibinfo{author}{F.~Aldakheel}, \bibinfo{author}{B.~Hudobivnik},
  \bibinfo{author}{E.~Artioli}, \bibinfo{author}{L.~B. da~Veiga},
  \bibinfo{author}{P.~Wriggers},
\newblock \bibinfo{title}{Curvilinear virtual elements for contact mechanics},
\newblock \bibinfo{journal}{Computer Methods in Applied Mechanics and
  Engineering} \bibinfo{volume}{372} (\bibinfo{year}{2020})
  \bibinfo{pages}{113394}.
\bibitem[{Benedetto et~al.(2014)Benedetto, Berrone, Pieraccini, and
  Scial{\`o}}]{benedetto2014virtual}
\bibinfo{author}{M.~F. Benedetto}, \bibinfo{author}{S.~Berrone},
  \bibinfo{author}{S.~Pieraccini}, \bibinfo{author}{S.~Scial{\`o}},
\newblock \bibinfo{title}{The virtual element method for discrete fracture
  network simulations},
\newblock \bibinfo{journal}{Computer Methods in Applied Mechanics and
  Engineering} \bibinfo{volume}{280} (\bibinfo{year}{2014})
  \bibinfo{pages}{135--156}.
\bibitem[{Aldakheel et~al.(2018)Aldakheel, Hudobivnik, Hussein, and
  Wriggers}]{aldakheel2018phase}
\bibinfo{author}{F.~Aldakheel}, \bibinfo{author}{B.~Hudobivnik},
  \bibinfo{author}{A.~Hussein}, \bibinfo{author}{P.~Wriggers},
\newblock \bibinfo{title}{Phase-field modeling of brittle fracture using an
  efficient virtual element scheme},
\newblock \bibinfo{journal}{Computer Methods in Applied Mechanics and
  Engineering} \bibinfo{volume}{341} (\bibinfo{year}{2018})
  \bibinfo{pages}{443--466}.
\bibitem[{Hussein et~al.(2020)Hussein, Hudobivnik, and
  Wriggers}]{hussein2020combined}
\bibinfo{author}{A.~Hussein}, \bibinfo{author}{B.~Hudobivnik},
  \bibinfo{author}{P.~Wriggers},
\newblock \bibinfo{title}{A combined adaptive phase field and discrete cutting
  method for the prediction of crack paths},
\newblock \bibinfo{journal}{Computer Methods in Applied Mechanics and
  Engineering} \bibinfo{volume}{372} (\bibinfo{year}{2020})
  \bibinfo{pages}{113329}.
\bibitem[{Antonietti et~al.(2017)Antonietti, Bruggi, Scacchi, and
  Verani}]{antonietti2017virtual}
\bibinfo{author}{P.~F. Antonietti}, \bibinfo{author}{M.~Bruggi},
  \bibinfo{author}{S.~Scacchi}, \bibinfo{author}{M.~Verani},
\newblock \bibinfo{title}{On the virtual element method for topology
  optimization on polygonal meshes: {A} numerical study},
\newblock \bibinfo{journal}{Computers \& Mathematics with Applications}
  \bibinfo{volume}{74} (\bibinfo{year}{2017}) \bibinfo{pages}{1091--1109}.
\bibitem[{Chi et~al.(2020)Chi, Pereira, Menezes, and Paulino}]{chi2020virtual}
\bibinfo{author}{H.~Chi}, \bibinfo{author}{A.~Pereira}, \bibinfo{author}{I.~F.
  Menezes}, \bibinfo{author}{G.~H. Paulino},
\newblock \bibinfo{title}{Virtual element method ({VEM})-based topology
  optimization: an integrated framework},
\newblock \bibinfo{journal}{Structural and Multidisciplinary Optimization}
  \bibinfo{volume}{62} (\bibinfo{year}{2020}) \bibinfo{pages}{1089--1114}.
\bibitem[{Smith(2013)}]{smith2013uncertainty}
\bibinfo{author}{R.~C. Smith}, \bibinfo{title}{Uncertainty quantification:
  theory, implementation, and applications}, volume~\bibinfo{volume}{12},
  \bibinfo{publisher}{SIAM}, \bibinfo{year}{2013}.
\bibitem[{Stefanou(2009)}]{stefanou2009stochastic}
\bibinfo{author}{G.~Stefanou},
\newblock \bibinfo{title}{The stochastic finite element method: past, present
  and future},
\newblock \bibinfo{journal}{Computer Methods in Applied Mechanics and
  Engineering} \bibinfo{volume}{198} (\bibinfo{year}{2009})
  \bibinfo{pages}{1031--1051}.
\bibitem[{Papadrakakis and Papadopoulos(1996)}]{papadrakakis1996robust}
\bibinfo{author}{M.~Papadrakakis}, \bibinfo{author}{V.~Papadopoulos},
\newblock \bibinfo{title}{Robust and efficient methods for stochastic finite
  element analysis using {Monte Carlo} simulation},
\newblock \bibinfo{journal}{Computer Methods in Applied Mechanics and
  Engineering} \bibinfo{volume}{134} (\bibinfo{year}{1996})
  \bibinfo{pages}{325--340}.
\bibitem[{Graham et~al.(2011)Graham, Kuo, Nuyens, Scheichl, and
  Sloan}]{graham2011quasi}
\bibinfo{author}{I.~G. Graham}, \bibinfo{author}{F.~Y. Kuo},
  \bibinfo{author}{D.~Nuyens}, \bibinfo{author}{R.~Scheichl},
  \bibinfo{author}{I.~H. Sloan},
\newblock \bibinfo{title}{{Quasi-Monte Carlo} methods for elliptic {PDE}s with
  random coefficients and applications},
\newblock \bibinfo{journal}{Journal of Computational Physics}
  \bibinfo{volume}{230} (\bibinfo{year}{2011}) \bibinfo{pages}{3668--3694}.
\bibitem[{Ghanem and Spanos(2003)}]{ghanem2003stochastic}
\bibinfo{author}{R.~G. Ghanem}, \bibinfo{author}{P.~D. Spanos},
  \bibinfo{title}{Stochastic finite elements: a spectral approach},
  \bibinfo{publisher}{Courier Corporation}, \bibinfo{year}{2003}.
\bibitem[{Xiu and Karniadakis(2002)}]{xiu2002wiener}
\bibinfo{author}{D.~Xiu}, \bibinfo{author}{G.~E. Karniadakis},
\newblock \bibinfo{title}{The {Wiener--Askey} polynomial chaos for stochastic
  differential equations},
\newblock \bibinfo{journal}{SIAM Journal on Scientific Computing}
  \bibinfo{volume}{24} (\bibinfo{year}{2002}) \bibinfo{pages}{619--644}.
\bibitem[{Babu{\v{s}}ka et~al.(2007)Babu{\v{s}}ka, Nobile, and
  Tempone}]{babuvska2007stochastic}
\bibinfo{author}{I.~Babu{\v{s}}ka}, \bibinfo{author}{F.~Nobile},
  \bibinfo{author}{R.~Tempone},
\newblock \bibinfo{title}{A stochastic collocation method for elliptic partial
  differential equations with random input data},
\newblock \bibinfo{journal}{SIAM Journal on Numerical Analysis}
  \bibinfo{volume}{45} (\bibinfo{year}{2007}) \bibinfo{pages}{1005--1034}.
\bibitem[{Xiu(2010)}]{xiu2010numerical}
\bibinfo{author}{D.~Xiu}, \bibinfo{title}{Numerical methods for stochastic
  computations: a spectral method approach}, \bibinfo{publisher}{Princeton
  University Press}, \bibinfo{year}{2010}.
\bibitem[{Khuri and Mukhopadhyay(2010)}]{khuri2010response}
\bibinfo{author}{A.~I. Khuri}, \bibinfo{author}{S.~Mukhopadhyay},
\newblock \bibinfo{title}{Response surface methodology},
\newblock \bibinfo{journal}{Wiley Interdisciplinary Reviews: Computational
  Statistics} \bibinfo{volume}{2} (\bibinfo{year}{2010})
  \bibinfo{pages}{128--149}.
\bibitem[{Fuhg et~al.(2021)Fuhg, Fau, and Nackenhorst}]{fuhg2021state}
\bibinfo{author}{J.~N. Fuhg}, \bibinfo{author}{A.~Fau},
  \bibinfo{author}{U.~Nackenhorst},
\newblock \bibinfo{title}{State-of-the-art and comparative review of adaptive
  sampling methods for kriging},
\newblock \bibinfo{journal}{Archives of Computational Methods in Engineering}
  \bibinfo{volume}{28} (\bibinfo{year}{2021}) \bibinfo{pages}{2689--2747}.
\bibitem[{Zheng et~al.(2022)Zheng, Beer, Dai, and Nackenhorst}]{zheng2022weak}
\bibinfo{author}{Z.~Zheng}, \bibinfo{author}{M.~Beer},
  \bibinfo{author}{H.~Dai}, \bibinfo{author}{U.~Nackenhorst},
\newblock \bibinfo{title}{A weak-intrusive stochastic finite element method for
  stochastic structural dynamics analysis},
\newblock \bibinfo{journal}{Computer Methods in Applied Mechanics and
  Engineering} \bibinfo{volume}{399} (\bibinfo{year}{2022})
  \bibinfo{pages}{115360}.
\bibitem[{Zheng et~al.(2023)Zheng, Valdebenito, Beer, and
  Nackenhorst}]{zheng2023stochastic}
\bibinfo{author}{Z.~Zheng}, \bibinfo{author}{M.~Valdebenito},
  \bibinfo{author}{M.~Beer}, \bibinfo{author}{U.~Nackenhorst},
\newblock \bibinfo{title}{A stochastic finite element scheme for solving
  partial differential equations defined on random domains},
\newblock \bibinfo{journal}{Computer Methods in Applied Mechanics and
  Engineering} \bibinfo{volume}{405} (\bibinfo{year}{2023})
  \bibinfo{pages}{115860}.
\bibitem[{Yosida(2012)}]{yosida2012functional}
\bibinfo{author}{K.~Yosida}, \bibinfo{title}{Functional analysis},
  \bibinfo{publisher}{Springer Science \& Business Media},
  \bibinfo{year}{2012}.
\bibitem[{Da~Veiga et~al.(2017)Da~Veiga, Dassi, and Russo}]{da2017high}
\bibinfo{author}{L.~B. Da~Veiga}, \bibinfo{author}{F.~Dassi},
  \bibinfo{author}{A.~Russo},
\newblock \bibinfo{title}{High-order virtual element method on polyhedral
  meshes},
\newblock \bibinfo{journal}{Computers \& Mathematics with Applications}
  \bibinfo{volume}{74} (\bibinfo{year}{2017}) \bibinfo{pages}{1110--1122}.
\bibitem[{Sakamoto and Ghanem(2002)}]{sakamoto2002polynomial}
\bibinfo{author}{S.~Sakamoto}, \bibinfo{author}{R.~Ghanem},
\newblock \bibinfo{title}{Polynomial chaos decomposition for the simulation of
  {non-Gaussian} nonstationary stochastic processes},
\newblock \bibinfo{journal}{Journal of Engineering Mechanics}
  \bibinfo{volume}{128} (\bibinfo{year}{2002}) \bibinfo{pages}{190--201}.
\bibitem[{Zheng and Dai(2017)}]{zheng2017simulation}
\bibinfo{author}{Z.~Zheng}, \bibinfo{author}{H.~Dai},
\newblock \bibinfo{title}{Simulation of multi-dimensional random fields by
  {Karhunen--Lo{\`e}ve} expansion},
\newblock \bibinfo{journal}{Computer Methods in Applied Mechanics and
  Engineering} \bibinfo{volume}{324} (\bibinfo{year}{2017})
  \bibinfo{pages}{221--247}.
\bibitem[{Zheng et~al.(2021)Zheng, Dai, Wang, and Wang}]{zheng2021sample}
\bibinfo{author}{Z.~Zheng}, \bibinfo{author}{H.~Dai},
  \bibinfo{author}{Y.~Wang}, \bibinfo{author}{W.~Wang},
\newblock \bibinfo{title}{A sample-based iterative scheme for simulating
  non-stationary {non-Gaussian} stochastic processes},
\newblock \bibinfo{journal}{Mechanical Systems and Signal Processing}
  \bibinfo{volume}{151} (\bibinfo{year}{2021}) \bibinfo{pages}{107420}.
\bibitem[{Pellissetti and Ghanem(2000)}]{pellissetti2000iterative}
\bibinfo{author}{M.~F. Pellissetti}, \bibinfo{author}{R.~G. Ghanem},
\newblock \bibinfo{title}{Iterative solution of systems of linear equations
  arising in the context of stochastic finite elements},
\newblock \bibinfo{journal}{Advances in Engineering Software}
  \bibinfo{volume}{31} (\bibinfo{year}{2000}) \bibinfo{pages}{607--616}.
\bibitem[{Keese and Matthies(2005)}]{keese2005hierarchical}
\bibinfo{author}{A.~Keese}, \bibinfo{author}{H.~G. Matthies},
\newblock \bibinfo{title}{Hierarchical parallelisation for the solution of
  stochastic finite element equations},
\newblock \bibinfo{journal}{Computers \& Structures} \bibinfo{volume}{83}
  (\bibinfo{year}{2005}) \bibinfo{pages}{1033--1047}.
\bibitem[{Blatman and Sudret(2010)}]{blatman2010adaptive}
\bibinfo{author}{G.~Blatman}, \bibinfo{author}{B.~Sudret},
\newblock \bibinfo{title}{An adaptive algorithm to build up sparse polynomial
  chaos expansions for stochastic finite element analysis},
\newblock \bibinfo{journal}{Probabilistic Engineering Mechanics}
  \bibinfo{volume}{25} (\bibinfo{year}{2010}) \bibinfo{pages}{183--197}.
\bibitem[{Young(2014)}]{young2014iterative}
\bibinfo{author}{D.~M. Young}, \bibinfo{title}{Iterative solution of large
  linear systems}, \bibinfo{publisher}{Elsevier}, \bibinfo{year}{2014}.
\bibitem[{Talischi et~al.(2012)Talischi, Paulino, Pereira, and
  Menezes}]{talischi2012polymesher}
\bibinfo{author}{C.~Talischi}, \bibinfo{author}{G.~H. Paulino},
  \bibinfo{author}{A.~Pereira}, \bibinfo{author}{I.~F. Menezes},
\newblock \bibinfo{title}{Polymesher: a general-purpose mesh generator for
  polygonal elements written in matlab},
\newblock \bibinfo{journal}{Structural and Multidisciplinary Optimization}
  \bibinfo{volume}{45} (\bibinfo{year}{2012}) \bibinfo{pages}{309--328}.
\bibitem[{Spanos et~al.(2007)Spanos, Beer, and Red-Horse}]{spanos2007karhunen}
\bibinfo{author}{P.~D. Spanos}, \bibinfo{author}{M.~Beer},
  \bibinfo{author}{J.~Red-Horse},
\newblock \bibinfo{title}{Karhunen--{Lo{\`e}ve} expansion of stochastic
  processes with a modified exponential covariance kernel},
\newblock \bibinfo{journal}{Journal of Engineering Mechanics}
  \bibinfo{volume}{133} (\bibinfo{year}{2007}) \bibinfo{pages}{773--779}.
\bibitem[{Saad(2011)}]{saad2011numerical}
\bibinfo{author}{Y.~Saad}, \bibinfo{title}{Numerical methods for large
  eigenvalue problems}, \bibinfo{publisher}{SIAM}, \bibinfo{year}{2011}.

\end{thebibliography}
 
\end{document}